\documentclass{scrarticle}

\usepackage[colorlinks=true,linkcolor=teal,citecolor=teal,hypertexnames=false]{hyperref}
\usepackage[margin=2.5cm]{geometry}
\usepackage{amsmath,amssymb,amsthm}
\usepackage{enumitem}
\usepackage{comment}
\usepackage{microtype}
\usepackage[english]{babel}
\usepackage[utf8]{inputenc}
\usepackage{mathrsfs}
\usepackage{bbm}

\usepackage{tikz}
\usetikzlibrary{matrix,arrows}
\usepackage{tikz-cd}
\usetikzlibrary{cd}

\mathcode`A="7041 \mathcode`B="7042 \mathcode`C="7043 \mathcode`D="7044
\mathcode`E="7045 \mathcode`F="7046 \mathcode`G="7047 \mathcode`H="7048
\mathcode`I="7049 \mathcode`J="704A \mathcode`K="704B \mathcode`L="704C
\mathcode`M="704D \mathcode`N="704E \mathcode`O="704F \mathcode`P="7050
\mathcode`Q="7051 \mathcode`R="7052 \mathcode`S="7053 \mathcode`T="7054
\mathcode`U="7055 \mathcode`V="7056 \mathcode`W="7057 \mathcode`X="7058
\mathcode`Y="7059 \mathcode`Z="705A

\newtheorem{Theorem}{Theorem}[section]
\newtheorem*{Theorem*}{Theorem}
\newtheorem{Lemma}[Theorem]{Lemma}

\newtheorem{Proposition}[Theorem]{Proposition}
\newtheorem*{Proposition*}{Proposition}
\newtheorem{Corollary}[Theorem]{Corollary}
\theoremstyle{definition}
\newtheorem{Definition}[Theorem]{Definition}
\newtheorem{Construction}[Theorem]{Construction}
\newtheorem*{Definition*}{Definition}
\newtheorem{Remark}[Theorem]{Remark}
\newtheorem{Notation}[Theorem]{Notation}
\newtheorem{Reduction}[Theorem]{Reduction}
\newtheorem{Example}[Theorem]{Example}
\newtheorem*{rem*}{Remark}

\newtheorem*{TheoremA}{Theorem A}
\newtheorem*{TheoremB}{Theorem B}

\newcommand{\eg}{\emph{e.g.~}}
\newcommand{\cf}{\emph{cf.~}}
\newcommand{\ie}{\emph{i.e.~}}
\newcommand{\loccit}{\emph{loc.~cit.~}}

\usepackage{thmtools}
\usepackage{thm-restate}

\newcommand{\Betti}[1]{#1_{B}}
\newcommand{\AMot}{A\mathrm{Mot}} %\mathbf{} used before
\newcommand{\HPk}{\mathrm{HdgPk}}

\newcommand{\Hdg}{\mathrm{Hdg}}
\newcommand{\Reg}{\mathrm{Reg}}

\newcommand{\Fil}{\mathrm{Fil}}
\newcommand{\cofib}{\mathrm{cofib}}
\newcommand{\coker}{\operatorname{coker}}

\newcommand{\im}{\operatorname{im}}

\newcommand{\qc}{\mathrm{qc}}
\newcommand{\Cl}{\mathrm{Cl}}

\newcommand{\Spec}{\operatorname{Spec}}
\newcommand{\Spm}{\operatorname{Spm}}
\newcommand{\Spf}{\operatorname{Spf}}
\newcommand{\reg}{\mathrm{reg}}

\newcommand{\Hom}{\mathrm{Hom}}
\newcommand{\RHom}{\mathrm{RHom}}
\newcommand{\Ext}{\mathrm{Ext}}
\newcommand{\End}{\mathrm{End}}
\newcommand{\Aut}{\mathrm{Aut}}
\newcommand{\id}{\mathrm{id}}

\newcommand{\Gal}{\mathrm{Gal}}
\newcommand{\Frob}{\mathrm{Frob}}

\newcommand{\bC}{\mathbb{C}}

\newcommand{\bF}{\mathbb{F}}

\newcommand{\bP}{\mathbb{P}}
\newcommand{\bQ}{\mathbb{Q}}
\newcommand{\bR}{\mathbb{R}}

\newcommand{\bZ}{\mathbb{Z}}

\usepackage{calrsfs}
\DeclareMathAlphabet{\pazocal}{OMS}{zplm}{m}{n}
\newcommand{\cA}{\mathcal{A}}
\newcommand{\cB}{\mathcal{B}}

\newcommand{\cF}{\mathcal{F}}

\newcommand{\cI}{I}

\newcommand{\cL}{\mathcal{L}}
\newcommand{\cM}{\mathcal{M}}
\newcommand{\cN}{\mathcal{N}}
\newcommand{\cO}{\mathcal{O}}

\newcommand{\sD}{\mathcal{D}}

\newcommand{\sH}{\mathcal{H}}

\newcommand{\KI}{K_{\infty}}
\newcommand{\FI}{F_{\infty}}
\newcommand{\CI}{\mathbb{C}_{\infty}}
\newcommand{\OI}{\cO_{\infty}}

\newcommand{\fh}{\mathfrak{h}}
\newcommand{\fj}{\mathfrak{j}}

\newcommand{\fl}{\mathfrak{l}}
\newcommand{\fm}{\mathfrak{m}}

\newcommand{\fp}{\mathfrak{p}}

\newcommand{\fq}{\mathfrak{q}}

\makeatletter
\newsavebox{\@brx}
\newcommand{\llangle}[1][]{\savebox{\@brx}{\(\m@th{#1\langle}\)}%
  \mathopen{\copy\@brx\kern-0.5\wd\@brx\usebox{\@brx}}}
\newcommand{\rrangle}[1][]{\savebox{\@brx}{\(\m@th{#1\rangle}\)}%
  \mathclose{\copy\@brx\kern-0.5\wd\@brx\usebox{\@brx}}}
\makeatother

\numberwithin{equation}{section}

\title{Regulators in the Arithmetic of Function Fields}
\author{Quentin Gazda\footnote{Sorbonne Universit\'e and Universit\'e Paris Cit\'e, CNRS, IMJ-PRG, 75005 Paris, \emph{Email:} \texttt{quentin@gazda.fr}}}
\date{}

\begin{document}

\maketitle

\begin{abstract}
As a natural sequel to the study of $A$-motivic cohomology initiated in \cite{gazda}, we develop a notion of regulator for rigid analytically trivial Anderson $A$\nobreakdash-motives. In accordance with the conjectural picture over number fields, we define it as the morphism on extension modules induced by the exact Hodge\nobreakdash--Pink realization functor. The purpose of this article is twofold: first, we prove a finiteness result for $A$-motivic cohomology; second, under the assumption that all the weights are negative, we show that the $\KI$-linearized regulator is an isomorphism. This proves the analogue of Beilinson’s celebrated conjecture in our setting.
\end{abstract}

{\footnotesize
\tableofcontents
}

\section{Introduction}\label{sec:introduction}
Recent years have witnessed the establishment of very general class formulas for Anderson $A$-modules. The most recent breakthrough was achieved by Angl\`es--Ngo Dac--Tavares Ribeiro \cite{angles}, who obtained the most general form to date (see also the forthcoming work of Caruso, Lucas, and the author \cite{caruso-gazda-lucas}, which allows for even greater generality). This achievement crowns a long line of works dating back to Carlitz \cite{carlitz} in 1935 and culminating in the breakthroughs of Taelman \cite{taelman-L} and V. Lafforgue \cite{lafforgue}. 

In the classical setting of number fields, class formulas at this level of generality are instead formulated in terms of mixed motives, in the form of Beilinson's conjectures. These are far\nobreakdash-reaching conjectures on special $L$-values, formulated by Beilinson in the 1980s in two celebrated papers \cite{beilinson-L,beilinson}. The entire picture is rooted in the notion of a \emph{Beilinson regulator}, conjecturally defined as follows. Let $M$ be a mixed motive over the field of rational numbers $\bQ$. Here, the term ``motive'' is understood in the sense of Deligne (\eg \cite[\S 1]{deligne-droite}). Consider the $\bQ$-vector space $\Ext^1_{\bQ}(\mathbbm{1},M)$ of $1$-fold extensions of the neutral motive $\mathbbm{1}$ by $M$ in the category of mixed motives over $\bQ$. It is expected that one can define a natural subspace $\Ext^1_{\bZ}(\mathbbm{1},M)$ consisting of \emph{extensions having everywhere good reduction} (see, for example, \cite{scholl}). The Hodge realization functor $\sH^+$, from the category of mixed motives to the category of mixed Hodge structures enriched with infinite Frobenii, is expected to be exact and should therefore induce a morphism between extension spaces:
\begin{equation}\label{eq:classical-reg}
\Ext^1_{\bZ}(\mathbbm{1},M)\longrightarrow \Ext^1_{\Hdg}(\mathbbm{1}^+,\sH^+(M)).
\end{equation}
Observe that the right-hand side is a finite-dimensional $\bR$-vector space. The above map is conjecturally the \emph{Beilinson regulator of $M$}. We denote it by $\Reg(M)$. The following statements are expected, although they remain far from being proved.
\begin{enumerate}[label=(\Roman*)]
\item\label{item:conj1} The space $\Ext^1_{\bZ}(\mathbbm{1},M)$ has finite dimension over $\bQ$;
\item\label{item:conj2} If $M$ is pure of weight $<-2$, then $\Reg(M)$ has dense image;
\item\label{item:conj3} There is a $\bQ$-structure $V(M)$ on the target of $\Reg(M)$, natural in $M$, such that we have $\det(\im\Reg(M))=L^*(M,0)\cdot \det V(M)$ as $\bQ$-structures on $\det \Ext^1_{\Hdg}(\mathbbm{1}^+,\sH^+(M))$.
\end{enumerate}
Conjectures \ref{item:conj2} and \ref{item:conj3} are referred to as Beilinson's conjectures. This text is concerned with the function field analogues of \ref{item:conj1} and \ref{item:conj2}, so we do not further comment on the (conjectural) definition of the special $L$-value $L^*(M,0)$ nor on the $\bQ$-structure $V(M)$. We refer the reader instead to the survey \cite{nekovar} for a comprehensive account of Beilinson's conjectures and their history. 

The present work grew out of an attempt to bridge these two pictures and to understand how Beilinson's conjectures could be formulated in the language of Anderson $A$-motives. The study of motivic cohomology for $A$-motives was initiated in \cite{gazda}, and the present text constitutes a natural sequel. Our primary focus is the definition of a \emph{Beilinson regulator} in this context and the study of analogues of conjectures \ref{item:conj1} and \ref{item:conj2}. We hope to address conjecture \ref{item:conj3} in a subsequent work. The reader is invited to consult \cite{caruso-gazda} for the definition of $L$-series of $A$-motives.

\subsection*{Anderson $A$-motives}
Let $(C,\cO_C)$ be a geometrically irreducible smooth projective curve over a finite field $\bF$ with $q$ elements and characteristic $p$, and fix a closed point $\infty$ of $C$. The $\bF$-algebra
\[
A=\cO_C(C\setminus \{\infty\})
\]
consists of functions on $C$ that are regular away from $\infty$. We denote by $K$ its fraction field. Let also $F$ be a finite field extension of $K$.

The notion of $A$-motives dates back to the pioneering work of Anderson \cite{anderson} and generalizes earlier ideas of Drinfeld \cite{drinfeld}. Let us state the definition of an \emph{Anderson $A$-motive}, leaving the details to Section \ref{sec:rat-mix-A-mot}. Throughout this text, unlabeled tensor products or fiber products are taken over $\bF$.
\begin{Definition*}[\ref{def:A-motive}, $A$-motives]
An \emph{Anderson $A$-motive $\underline{M}$ over $F$} consists of a finite locally free $A\otimes F$-module $M$ together with an isomorphism of $A\otimes F$-modules
\[
\tau_M:(\tau^*M)|_{(\Spec A\otimes F)\setminus V(\fj)} \stackrel{\sim}{\longrightarrow} M|_{(\Spec A\otimes F)\setminus V(\fj)},
\]
where $\tau=\Frob_{F}$ denotes the $A$-linear $|\bF|$-Frobenius on $A\otimes F$, and where $V(\fj)$ is the effective Cartier divisor on $\Spec A\otimes F$ associated with the locally free ideal $\fj:=\ker(A\otimes F\to F,~a\otimes r\mapsto ar)$.
\end{Definition*}

We may assume without loss of generality\footnote{In fact, even that $F$ is purely inseparable over $K$, see Reduction \ref{red:F-totally-ramified}.} that there is a unique place of $F$ above $\infty$, which we also denote by $\infty$. We write $F_{\infty}$ for the completion of $F$ with respect to $\infty$. We will focus our attention on $A$-motives that are \emph{rigid analytically trivial} (see Subsection~\ref{subsec:betti} for details). Given an $A$-motive $\underline{M}$ over $F$, we consider the $A$-module
\[
\Betti{\underline{M}}:=\{f\in M\otimes_{A\otimes K^{\natural}}\CI\langle A \rangle ~|~f=\tau_M(\tau^*f)\}
\]
where $K^{\natural}$ stands for the perfection of $K$ inside $F$, and $\CI$ is the completion of an algebraic closure of $\KI$, in which $K^{\natural}$ embeds canonically. Here, $\CI\langle A \rangle$ denotes a certain affinoid algebra over the completion $\CI$ (see \cite[\S 2.3.3]{hartl-juschka}, \cite[\S 3]{gazda-maurischat2}, or Section~\ref{sec:rat-mix-A-mot} below for details). The module $\Betti{\underline{M}}$ is called the \emph{Betti realization of $\underline{M}$}.  To ensure good properties of the functor $\underline{M}\mapsto \Betti{\underline{M}}$ (\eg exactness), we restrict to the subcategory of \emph{rigid analytically trivial} $A$-motives, \ie those $\underline{M}$ for which $\Betti{\underline{M}}$ generates $M\otimes_{A\otimes K^{\natural}}\CI\langle A \rangle$ over $\CI\langle A \rangle$ (Definition \ref{def:rigid analytically trivial}). 

The category of \emph{rigid analytically trivial} $A$-motives over $F$, denoted throughout this text by $\AMot_F^{\operatorname{rat}}$, is an exact $A$-linear category which plays the role of the classical category of mixed motives over $\bQ$. In this paradigm, the functor $\underline{M}\mapsto \Betti{\underline{M}}$ plays the role of the singular realization; the reader is referred to \cite{motif} for surveys on this analogy. In the same way that complex conjugation acts via an involution on the Betti realization of a motive over $\bQ$, we show that, given $\underline{M}$ an object of $\AMot_F^{\operatorname{rat}}$, the $A$-module $\Betti{\underline{M}}$ carries a natural action of the profinite Galois group $G_\infty:=\Gal(\KI^s|\KI)$, where $\KI^s$ denotes the separable closure of $\KI$ inside $\CI$ (\cf Proposition~\ref{prop:finite-extension-betti} and Corollary \ref{cor:betti-exact}). 

Note that $A$-motives do not always carry a weight filtration; those that do are called \emph{mixed} (see \cite[\S 3]{gazda}). However, since the main results of this text do not require mixedness—aside from the evident analogy with the number field theory—we deliberately disregard the notion of mixedness. Nevertheless, there remains a notion of \emph{weights} that plays a major role.

In \cite{gazda}, we established the definition of $\Ext^1_{\cO_F}(\mathbbm{1},\underline{M})$, a natural $A$-module consisting of \emph{integral} extensions of the neutral $A$-motive $\mathbbm{1}$ by $\underline{M}$ in the category $\AMot_F^{\operatorname{rat}}$. However, there are two obstructions to the finite generation of this module, preventing a naive analogue of conjecture \ref{item:conj1} from holding in our context:
\begin{enumerate}
\item The first reason, discussed in detail in \cite[\S 5]{gazda}, is related to the fact that taking Hodge filtrations is not an exact operation on the full class of exact sequences. This issue is addressed by the notion of \emph{regulated extensions}, introduced in Definition~5.1 of \loccit The resulting $A$-submodule $\Ext^{1,\text{reg}}_{\cO_F}(\mathbbm{1},\underline{M})$ of regulated extensions, however, may still fail to be finitely generated.

\item The second reason is due to the profinite nature of the absolute Galois group $G_\infty$. We address this issue by introducing the notion of \emph{analytic reduction at $\infty$}, defined as follows.
\end{enumerate}
The exactness of the Betti realization functor induces a morphism at the level of extension modules:
\begin{equation}\label{eq:rlambda}
r_{B}:\Ext^{1,\text{reg}}_{\cO_F}(\mathbbm{1},\underline{M})\longrightarrow \operatorname{H}^1(G_\infty,\Betti{\underline{M}})
\end{equation}
where the target denotes continuous Galois cohomology. We say that an extension $[\underline{E}]$ in $\Ext^{1,\text{reg}}_{\cO_F}(\mathbbm{1},\underline{M})$ has \emph{analytic reduction at $\infty$} if it lies in the kernel of $r_{B}$; equivalently, if the induced extension of continuous $G_\infty$-modules splits. We denote by $\Ext^{1,\text{reg}}_{\cO_F,\infty}(\mathbbm{1},\underline{M})$ the kernel of $r_B$, and by $\operatorname{Cl}(\underline{M})$ its cokernel. Our first main theorem is the following (repeated from Theorem \ref{thm:finiteness-motcoh}).
%\begin{Theorem*}[\ref{thm:finiteness-motcoh}]
\begin{TheoremA}
The $A$-modules $\Ext^{1,\operatorname{reg}}_{\cO_F,\infty}(\mathbbm{1},\underline{M})$ and $\operatorname{Cl}(\underline{M})$ are finitely generated. If, in addition, the weights of $\underline{M}$ are all negative, then $\operatorname{Cl}(\underline{M})$ is finite. 
\end{TheoremA}
%\end{Theorem*}
Let us make a few comments on Theorem A.
\begin{enumerate}[label=$-$]
    \item The above theorem should be understood as the analogue of conjecture \ref{item:conj1} for rigid analytically trivial Anderson $A$-motives. Indeed, all extensions of classical mixed motives are \emph{regulated} in the obvious sense, and hence \eqref{eq:rlambda} corresponds, in the classical setting, to the morphism of $\bQ$-vector spaces
\[
r_B:\Ext^1_{\bZ}(\mathbbm{1},M)\longrightarrow \operatorname{H}^1(\Gal(\bC|\bR),M_B)
\]
induced by the exactness of the Betti realization (here, $M$ is a mixed motive over $\bQ$, and $M_B$ denotes its Betti realization). Yet, in the $\bQ$-linear category of mixed motives, the right-hand side vanishes, which amounts to saying that \emph{all} extensions in $\Ext^1_{\bZ}(\mathbbm{1},M)$ have \emph{analytic reduction at $\infty$}.
\item As a second remark, Theorem A also shows that $\Ext^{1,\text{reg}}_{\cO_F}(\mathbbm{1},\underline{M})$ is almost \emph{never} finitely generated. Indeed, its size is roughly comparable to that of $\operatorname{H}^1(G_\infty,\Betti{\underline{M}})$ which need not be finitely generated as $G_\infty$ is not topologically finitely generated: via local reciprocity, its maximal abelian quotient has the principal-unit group of $\mathcal O_\infty$ as an open subgroup (up to the usual reciprocity identification), and that group is a countable product of copies of $\mathbb Z_p$.
\item If the $A$-motive $\underline{M}$ arises as (the dual of) the motive of a uniformizable and abelian Anderson module~$E$ with everywhere good reduction over $F$, we prove in the companion paper \cite{gazda-companion} that, up to a twist by the module of K\"ahler differentials, $\Ext^{1,\text{reg}}_{\cO_F,\infty}(\mathbbm{1},\underline{M})$ recovers Taelman’s unit module $U(E/\cO_F)$ and $\operatorname{Cl}(\underline{M})$ recovers its class module $H(E/\cO_F)$ (see \cite{taelman-dirichlet}). In particular, $A$-motivic cohomology naturally generalizes Taelman's theory to general $A$-motives, without assuming effectivity or everywhere good reduction. 
\end{enumerate}

This discussion suggests that the module $\Ext^{1,\operatorname{reg}}_{\cO_F,\infty}(\mathbbm{1},\underline{M})$ is the \emph{appropriate} source of a regulator. This expectation is in agreement with recent computations obtained jointly with Maurischat \cite{gazda-maurischat-ext} for Carlitz tensor powers. The Hodge--Pink side of the theory, which we describe next, appears to confirm this insight.

\subsection*{Hodge-Pink structures}
In an innovative unpublished monograph \cite{pink}, Pink defined and studied a general theory of Hodge structures in function field arithmetic. The appropriate objects of study, highlighted in \loccit and today called \emph{Hodge-Pink structures (over $K_{\infty}^s$)}, consist of pairs $\underline{H}=(H,\fq_H)$, where
\begin{itemize}
\item $H$ is a finite-dimensional $\KI$-vector space,
\item $\fq_{H}$ is a $\KI^s[\![\fj]\!]$-lattice in the $\KI^s(\!(\fj)\!)$-vector space $H_{\KI^s (\!(\fj)\!)}:=H\otimes_{\KI,\nu}\KI^s(\!(\fj)\!)$.
\end{itemize}
We refer to Section \ref{sec:ext-of-MHPS} for details. Here $\KI^s[\![\fj]\!]$ denotes the completion of $A\otimes \KI^s$ with respect to the $\fj$-adic topology and $\KI^s(\!(\fj)\!)$ its fraction field. The map $a\mapsto a\otimes 1$ extends to a continuous map $\nu:\KI\to \KI^s(\!(\fj)\!)$ (Proposition \ref{prop:nu-extends-by-continuity}).

Classically, Hodge structures arising from motives naturally carry an involution, called \emph{the infinite Frobenius}, obtained functorially from the action of complex conjugation. We propose the following function field counterpart: 
\begin{Definition*}[\ref{def:infinite-frob}, Infinite Frobenius]
An \emph{infinite Frobenius} for $\underline{H}$ is a $\KI$-linear continuous representation $\phi:G_\infty\to \Aut_{\KI}(H)$, where $H$ carries the discrete topology, such that for all $\sigma\in G_\infty$,
\[
\phi(\sigma)\otimes_{\KI} \sigma:H_{\KI^s(\!(\fj)\!)}\longrightarrow H_{\KI^s(\!(\fj)\!)}
\]
preserves the Hodge-Pink lattice $\fq_H$. We denote by $\HPk^+$ the category of pairs $(\underline{H},\phi_H)$ where $\underline{H}$ is a Hodge-Pink structure and $\phi_H$ is an infinite Frobenius for $\underline{H}$.
\end{Definition*}

Let $\underline{H}^+$ be a Hodge-Pink structure equipped with an infinite Frobenius, and let $\mathbbm{1}^+$ denote the neutral object in $\HPk^+$. As for $A$-motives, the extension space $\Ext^1_{\HPk^+}(\mathbbm{1}^+,\underline{H}^+)$ need not be finite-dimensional over~$\KI$. This differs from the number field setting. The reasons are almost identical to those for $A$-motives: taking Hodge filtrations is not an exact operation, and this space of extensions is intertwined with the Galois cohomology of the profinite group $G_\infty$. Using Pink's notion of Hodge additivity---which inspired that of \emph{regulated extensions} in \cite{gazda}---we consider the subspace
\begin{equation}\label{eq:subsub}
\Ext^{1,\operatorname{ha}}_{\HPk^+}(\mathbbm{1}^+,\underline{H}^+)
\end{equation}
of Hodge-additive extensions. We denote by $\Ext^{1,\operatorname{ha}}_{\infty}(\mathbbm{1}^+,\underline{H}^+)$ the subspace of \eqref{eq:subsub} consisting of classes of extensions whose infinite Frobenius splits, called extensions \emph{with analytic reduction at $\infty$}.

To a rigid analytically trivial $A$-motive $\underline{M}$ over $F$, we associate a Hodge-Pink structure with an infinite Frobenius, whose underlying space is $\underline{M}_B\otimes_{A} \KI$ (Section \ref{sec:rat-mix-A-mot}). Our construction generalizes that of Hartl\nobreakdash--Juschka \cite{hartl-juschka} in the case where the point $\infty$ has degree one. This assignment turns out to be an exact functor:
\[
\sH^+:\AMot_F^{\operatorname{rat}}\longrightarrow \HPk^+
\]
and it will play the role of the Hodge realization functor (see Definition \ref{def:hodge-pink-real-functor}). Its construction requires treating elements of $\underline{M}_B$ as meromorphic functions on the rigid analytic space associated with $\Spec A\otimes \CI$ and studying their behaviour at the point associated with $\fj$.

This paves the way for the definition of a Hodge-Pink regulator. The exactness of $\sH^+$ induces an $A$-linear morphism at the level of extension modules which, by design, sends regulated, analytically reduced extensions of $A$-motives to Hodge-additive, analytically reduced extensions of Hodge-Pink structures with infinite Frobenii. This motivates the next definition.
\begin{Definition*}[\ref{def:regulator}, Hodge-Pink regulator]
We call the \emph{Hodge-Pink regulator of $\underline{M}$}, and denote it by $\Reg(\underline{M})$, the $A$-linear morphism
\[
\Reg(\underline{M}):\Ext^{1,\operatorname{reg}}_{\cO_F,\infty}(\mathbbm{1},\underline{M})\longrightarrow \Ext^{1,\operatorname{ha}}_{\infty}(\mathbbm{1}^+,\sH^+(\underline{M}))
\]
induced by the exact functor $\sH^+$.
\end{Definition*}

It is the counterpart of Beilinson's regulator in this setting. If $\underline{M}$ arises from an abelian Anderson module~$E$ as above, we prove in the forthcoming work \cite{gazda-companion} that, up to a twist by the module of K\"ahler differentials, it coincides with \emph{Taelman's regulator}
\[
\log_E:U(E/\cO_F)\longrightarrow \mathrm{Lie}_E(F_{\infty})/(K_{\infty}\Lambda_E\cap \mathrm{Lie}_E(F_{\infty}))
\]
where $\Lambda_E$ denotes the \emph{period lattice of $E$} and $K_{\infty}\Lambda_E$ the $\KI$-subspace of $\mathrm{Lie}_E(K^s_{\infty})$ it generates. 
Our second main result therefore generalizes Taelman's theory to a more general class of $A$-motives.

\begin{TheoremB}
Assume that the weights of $\underline{M}$ are all negative. Then the kernel of $\Reg(\underline{M})$ is the torsion submodule of $\Ext^{1,\operatorname{reg}}_{\cO_F,\infty}(\mathbbm{1},\underline{M})$ and its image forms an $A$-lattice of full rank in the $\KI$-vector space $\Ext^{1,\operatorname{ha}}_{\infty}(\mathbbm{1}^+,\sH^+(\underline{M}))$. Equivalently, the $\KI$-linearized regulator is an isomorphism.
\end{TheoremB}

This is the counterpart of Beilinson's conjecture \ref{item:conj2}. Theorem B immediately determines the rank of the $A$-motivic cohomology module. Combining it with the explicit dimension formula for the Hodge-Pink extension space gives the following corollary. 
\begin{Corollary}
Assume that the weights of $\underline{M}$ are all negative. Then
\[
\mathrm{rk}_A~\Ext^{1,\operatorname{reg}}_{\cO_F,\infty}(\mathbbm{1},\underline{M})=(\text{sum~of~the~positive~Hodge-Pink~weights~of~}\underline{M})-\mathrm{rk}_A~ \underline{M}_B^{G_{\infty}}.
\] 
\end{Corollary}

\begin{Example}[Carlitz twists]
Let $C=\mathbb{P}^1_{\bF}$ with affine coordinate $t$, and let $n>0$. Let $K=\mathbb{F}(\theta)$ be (isomorphic to) the function field of $\mathbb{P}^1_{\bF}$ and let $F/K$ be a finite extension. The \emph{Carlitz $n$-twist over $F$} is the $t$-motive $\underline{A}(n):=(F[t],(t-\theta)^{-n})$ (it is the function field analogue of Tate twists). In this case, the corollary gives:
\begin{equation}\label{eq:K-theory-rank}
\mathrm{rk}_{\mathbb{F}[t]}~\Ext^{1,\operatorname{reg}}_{\cO_F,\infty}(\mathbbm{1},\underline{A}(n))=n[F:K]-d_n
\end{equation}
where $d_n$ is the number of infinite places of $F$ at which $(-\theta)^n$ admits a $(q-1)$th root. As expected, when $n=1$ this matches\footnote{The reference erroneously states the formula with the separable degree instead of the degree $[F:K]$, although the proof given there shows that it is the latter that holds true.} the rank of $U(C/\cO_F)$ computed by Taelman \cite[Proposition~1]{taelman-dirichlet}.
\end{Example}

If one pushes the analogy further with the number field situation, \eqref{eq:K-theory-rank} should correspond to the rank of the odd $K$-group $K_{2n-1}(\cO_F)$ for a number field $F$; however, whereas \eqref{eq:K-theory-rank} grows with $n$, the rank of $K_{2n-1}(\cO_F)$ is bounded with respect to $n$. A curious discrepancy! This is perhaps where function field Hodge structures (rather than Hodge-Pink structures; see Definition \ref{def:hodge-structure}) should be rehabilitated in the definition of a ``Hodge regulator'': 
\begin{equation}\label{eq:hodge-regulator-intro}
\Ext^{1,\reg}_{\cO_F,\infty}(\mathbbm{1},\underline{M}) \longrightarrow \Ext^1_{\Hdg,\infty}(\mathbbm{1},(H,\Fil^{\bullet} H_{\KI^s})).
\end{equation}
Note that this map is well-defined (once one clarifies what the right-hand side means, \eg Definition \ref{def:hodge-structure}), since Hodge-additive extensions induce extensions of function field Hodge structures. Although the target has smaller $\KI$-dimension than the $A$-rank of the source, we believe that there exists a natural submodule $Z$ of $\Ext^{1,\reg}_{\cO_F,\infty}(\mathbbm{1},\underline{M})$ such that the image of the restriction of \eqref{eq:hodge-regulator-intro} to $Z$ is a full-rank $A$-lattice. We expect this submodule $Z$ to be related to the eponym module in Taelman's conjecture, studied by Anglès--Ngo Dac--Tavares Ribeiro in \cite{ANT2}. This is the subject of ongoing investigations by the author.

\subsection*{Methods}
Our proofs of the main theorems draw strong inspiration from the work of Lafforgue \cite{lafforgue}, Pink \cite{pink-isogenies}, and Mornev \cite{mornev-shtuka}, and hinge on the concept of \emph{shtuka models}. We associate non\nobreakdash-canonically to $\underline{M}$---which lives over the affine curve $\Spec A\otimes F$---a \emph{shtuka model} $\underline{\cM}_0=(\cM_0,\cN_0,\tau_0)$ over the surface $(\Spec A)\times X$ (Proposition \ref{prop:existence-of-X-models}) where $X$ is a smooth projective curve over $\bF$ whose function field is $F$. Its coherent cohomology defines an object in $\sD(A)$, the derived category of $A$-modules\footnote{$G_{\underline{M}}$ is in fact defined differently in Definition \ref{def:Gm} below, but is quasi-isomorphic to the given formula by Proposition \ref{prop:ge-in-terms-of-cohomology}.}, whose cohomology is concentrated in degrees $0$ and $1$:
\[
G_{\underline{M}}:=\operatorname{R}\!\Gamma\left(\Spec A\times X,\cM_0\xrightarrow{\iota-\tau_0}\cN_0\right).
\]
Although there are several possible choices for $\underline{\cM}_0$, the above object only depends on $\underline{M}$ up to quasi\nobreakdash-isomorphism. It plays a role analogous to that of \emph{shtuka cohomology} as defined and studied extensively in \cite{mornev-shtuka}. We show in Proposition \ref{prop:motivic-exact-sequence} that there is a natural long exact sequence of $A$-modules:
\[
0\to \Hom_{\AMot_F}(\mathbbm{1},\underline{M})\to \Betti{\underline{M}}^{G_{\infty}} \to \operatorname{H}^0(G_{\underline{M}}) \to \Ext^{1,\operatorname{reg}}_{\cO_F}(\mathbbm{1},\underline{M}) \xrightarrow{\Betti{r}} \operatorname{H}^1(G_{\infty},\Betti{\underline{M}}) \to \operatorname{H}^1(G_{\underline{M}}) \to 0
\]
where $r_B$ is the map defined in \eqref{eq:rlambda}. The first part of Theorem A is deduced from the perfectness of the complex $G_{\underline{M}}$, which itself follows from the finiteness of coherent cohomology for the proper $A$-scheme $(\Spec A)\times X\to \Spec A$.

To some extent, the second part of Theorem A follows a similar pattern but is more involved. When the weights of $\underline{M}$ are negative, we further associate to $\underline{M}$ a \emph{shtuka model} $\underline{\cM}$ on $C\times X$ (Proposition \ref{prop:existence-CxX-shtuka-models}). An incarnation of this model in the context of Drinfeld modules with everywhere good reduction already appears in \loccit under the name of a \emph{global model}. The key ingredient of our proof is the relation between the behaviour of $\underline{\cM}$ in the formal neighbourhood of $\{\infty\}\times \{\infty\}\hookrightarrow C\times X$ and the space of Hodge-additive extensions of Hodge-Pink structures having analytic reduction at $\infty$ (Corollary \ref{cor:shtuka-model-to-hodge-pink-additive}): this link is materialized by a fiber sequence in~$\sD(K_{\infty})$:
\[
\Betti{\underline{M}}^{G_{\infty}}[0]\otimes_A \KI \longrightarrow \operatorname{R}\!\Gamma\left(\Spf \cO_{\infty}\hat{\times} X,\cM\xrightarrow{\iota-\tau}\cN\right)\otimes_{\OI}\KI \longrightarrow \Ext^{1,\operatorname{ha}}_{\infty}(\mathbbm{1}^+,\sH^+(\underline{M}))[0].
\]
The middle term is isomorphic to $G_{\underline{M}}\otimes_A \KI$, as we show using Grothendieck's theorem on formal functions (see Section \ref{sec:proof}). In particular, $G_{\underline{M}}\otimes_A \KI$ sits in degree zero, and hence $\operatorname{H}^1(G_{\underline{M}})\cong \operatorname{Cl}(\underline{M})$ is finite, as desired.

This also shows that the rank of $\Ext^{1,\operatorname{reg}}_{\cO_F,\infty}(\mathbbm{1},\underline{M})$ matches the dimension of the vector space $\Ext^{1,\operatorname{ha}}_{\infty}(\mathbbm{1}^+,\sH^+(\underline{M}))$. To prove Theorem B, it therefore remains to establish the injectivity of the Hodge-Pink regulator after extending scalars to $\KI$. We achieve this by first expressing motivic extension modules in terms of meromorphic solutions to $\tau$-difference equations in the style of Sinha \cite{sinha}, and then by using an explicit description of the regulator in terms of the principal parts at $\fj$ of these solutions (Proposition \ref{prop:regulator-injective}).

\paragraph{Acknowledgments.} I learned \emph{a posteriori} that the problem of constructing natural finitely generated extension modules from the category of $A$-motives---such as $\Ext^{1,\operatorname{reg}}_{\cO_F,\infty}(\mathbbm{1},\underline{M})$---had remained open and had been raised by Dinesh Thakur. I am much indebted to him for his interest in my work. While preparing earlier versions of the manuscript, I benefited greatly from many exchanges and discussions with the following people, to whom I wish to reiterate my gratitude: Gebhard B\"ockle, Bhargav Bhatt, Christopher Deninger, Urs Hartl, Annette Huber-Klawitter, Maxim Mornev, and Federico Pellarin. I am also grateful to the Max Planck Institute for Mathematics in Bonn for its hospitality and financial support.

\section{Hodge-Pink structures and their extensions}\label{sec:ext-of-MHPS}
In the inspiring monograph \cite{pink}, Pink introduced the appropriate function field analogue of Hodge structures, called Hodge-Pink structures. In this section, we give the definition and turn without further ado to the computation of extension spaces. The reader interested in additional motivation is invited to consult \loccit or the survey of Hartl--Juschka \cite{hartl-juschka}. One innovation of ours is the consideration of \emph{infinite Frobenii} in this context.

\subsection{Hodge-Pink structures}
\subsubsection*{Definitions}
Let $(C,\cO_C)$ be a geometrically irreducible smooth projective curve over $\bF$, and let $\infty$ be a closed point of $C$. The $\bF$-algebra
\[
A=\cO_C(C\setminus \{\infty\})
\]
consists of functions on $C$ that are regular away from $\infty$. We denote by $K$ the fraction field of $A$ and by $\KI$ the completed local field of $C$ at $\infty$. Let $\FI$ be a finite field extension of $\KI$ and fix a separable closure $\FI^s$ of $\FI$. We denote by $\fj\subset A\otimes A$ the kernel of the multiplication map $A\otimes A\to A$.  

By naive analogy with the classical setting, one could formulate the following definition of a ``(pre-)Hodge structure'':
\begin{Definition}[Hodge structure]\label{def:hodge-structure}
A \emph{function field Hodge structure (over $\FI^s$)} consists of a pair $(H,\Fil^{\bullet} H_{\FI^s})$ where
\begin{itemize}
\item $H$ is a finite-dimensional $\KI$-vector space,
\item $\Fil^{\bullet} H_{\FI^s}$ is a decreasing, separated, and exhaustive $\bZ$-filtration of $H_{\FI^s}:=H\otimes_{\KI} \FI^s$. 
\end{itemize}
\end{Definition}

However, Pink shows in \cite{pink} that this definition (or rather its mixed version) does not satisfy the Tannakian requirements one would expect for a category of (mixed) Hodge structures. He instead suggests replacing the Hodge filtration by the data of a lattice. The definition of \emph{Hodge-Pink structures} is rooted in the ring $\FI[\![\fj]\!]$: for $L$ a separable extension of $\FI$, $L[\![\fj]\!]$ is defined as the completion of the ring $A\otimes L$ along its maximal ideal $\fj$:
\[
L[\![\fj]\!]:= \varprojlim_m A\otimes L/\fj^m.
\]
Note that $L[\![\fj]\!]$ is a discrete valuation ring with maximal ideal $\fj$ and residue field $L$. Since each of the terms $A\otimes L/\fj^m$ in the limit is a finite-dimensional $L$-vector space, they carry a canonical topology (and even a canonical structure of a Banach $L$-space if $L$ is complete); we endow $L[\![\fj]\!]$ with the inverse-limit topology. We denote by $L(\!(\fj)\!)$ its field of fractions.

Given an $\KI$-vector space $V$, there are (at least) two ways to obtain an $\KI[\![\fj]\!]$-module from $V$: the \emph{base}-wise way would be $V\otimes_{\KI} \KI[\![\fj]\!]$, viewing $\KI$ as a subring of $\KI[\![\fj]\!]$. The \emph{coefficient}-wise way is made possible by the following statement (which generalizes \cite[Proposition 3.1]{pink}):
\begin{Proposition}\label{prop:nu-extends-by-continuity}
The map $A\to \KI[\![\fj]\!]$ induced by $a\mapsto a\otimes 1$ extends uniquely to a continuous map $\nu:\KI\to \KI[\![\fj]\!]$.
\end{Proposition}
\begin{proof}
For $m\geq 1$, set $B_m:=A\otimes \KI/\fj^m$ and let $\nu_m:A\to B_m$ be the map $a\mapsto a\otimes 1$. Let us first prove the statement at the level of $B_m$. Given a nonzero $a\in A$, we have $\nu_m(a)\equiv 1\otimes a\pmod{\fj}$, which is a unit in $B_m/\fj\cong \KI$. Since elements of $\fj$ are nilpotent in $B_m$, this implies that $\nu_m(a)$ is itself a unit in $B_m$. In particular, $\nu_m$ factors uniquely through $A\to K$. 

The composite of the extended map $\nu_{K,m}:K\to B_m$ with reduction modulo $\fj$ is the canonical inclusion $K\hookrightarrow\KI$. We deduce that if $\pi_{\infty}\in K$ is a uniformizer at $\infty$, then there exists $n\in \fj$ such that $\nu_{K,m}(\pi_\infty)=(1\otimes \pi_\infty)+n$. Since $n$ is a nilpotent element of $B_m$, we have for $N>m$:
\[
\nu_{K,m}(\pi_\infty)^N=(1\otimes \pi_\infty)^{N-m}\left(\sum_{k=0}^{m-1}{\binom{N}{k} (1\otimes \pi_\infty)^{m-k}n^k}\right).
\] 
In particular, $\nu_{K,m}$ maps a basis of open neighbourhoods of zero in $K$ (for the $\infty$-adic topology) to subsets of open neighbourhoods in $B_m$. Since $B_m$ is complete, $\nu_{K,m}$ extends uniquely to a continuous map $\nu_{\KI,m}:\KI\to B_m$.

This shows that $\nu:A\to \KI[\![\fj]\!]$ extends to a continuous map $\KI\to \KI[\![\fj]\!]$ by the formula $x\mapsto (\nu_{\KI,m}(x))_{m\geq 1}$. To show that it is the unique continuous extension, note that a continuous map $\KI\to \KI[\![\fj]\!]$ extending $\nu$ would induce, by the product topology on $\KI[\![\fj]\!]$, a continuous map $\KI\to B_m$ extending $\nu_m$; hence it must be the map $\nu_{\KI,m}$ constructed above.
\end{proof}

We are now in a position to define Hodge-Pink structures, following \cite[Definition 3.2]{pink}. By an \emph{$L[\![\fj]\!]$-lattice $\fh$ in a finite-dimensional $L(\!(\fj)\!)$-vector space $V$} we mean a finitely generated $L[\![\fj]\!]$-submodule of $V$ that contains a basis.
\begin{Definition}[Hodge-Pink structure]\label{def:pMHPS}
A \emph{Hodge-Pink structure $\underline{H}$ (over $\FI^s$)} consists of the data of $(H,\fq)$ where 
\begin{itemize}
\item $H$ is a finite-dimensional $\KI$-vector space,
\item $\fq=\fq_{H}$ is a $\FI^s[\![\fj]\!]$-lattice in the $\FI^s(\!(\fj)\!)$-vector space $H_{\FI^s(\!(\fj)\!)}:=H\otimes_{\KI,\nu}\FI^s(\!(\fj)\!)$.
\end{itemize}
We call $\fq$ the \emph{Hodge-Pink lattice of $\underline{H}$}. We let $\fp_{H}:=H_{\FI^s[\![\fj]\!]}=H\otimes_{\KI,\nu}\FI^s[\![\fj]\!]$ and call it the \emph{tautological lattice of $\underline{H}$}.
\end{Definition}

We gather Hodge-Pink structures into a $\KI$-linear category $\HPk$. A \emph{morphism from $\underline{H}=(H,\fq)$ to $\underline{H}'=(H',\fq')$ in $\HPk$} is a $\KI$-linear morphism $f:H\to H'$ such that its $\FI^s(\!(\fj)\!)$\nobreakdash-linear extension $f_{\FI^s(\!(\fj)\!)}:H_{\FI^s(\!(\fj)\!)}\to H_{\FI^s(\!(\fj)\!)}'$ preserves the Hodge-Pink lattices; \ie $f_{\FI^s(\!(\fj)\!)}(\fq)\subset \fq'$. We call $f$ \emph{strict} whenever we have:
\[
f_{\FI^s(\!(\fj)\!)}(\fq)=\fq'\cap \im f_{\FI^s(\!(\fj)\!)}.
\]

\subsubsection*{Extensions}
We are interested in computing extension modules in the category $\HPk$. However, $\HPk$ is not abelian, and we must clarify the notion of exact sequences we want to consider.
\begin{Definition}[Strict exact sequence]\label{def:strict-sequence-Hdg}
A sequence $S:0\to \underline{H}'\to \underline{H}\to \underline{H}''\to 0$ in $\HPk$ is called \emph{exact} if its underlying sequence of $\KI$-vector spaces is exact, and \emph{strict} if the morphisms of $S$ are.
\end{Definition}

\begin{Remark}\label{rem:strict-exact-sequences}
It is formal to verify that an exact sequence $S$ in $\HPk$ is strict if and only if the induced sequence $S_{\fq}$ of Hodge-Pink lattices is exact in the category of $\FI^s[\![\fj]\!]$-modules.
\end{Remark}

The category $\HPk$, together with the notion of strict exact sequences, forms an exact category. Given two Hodge-Pink structures $\underline{X}$ and $\underline{Y}$, we denote by $\Ext^1_{\HPk}(\underline{X},\underline{Y})$ the set of equivalence classes of strict exact sequences of the form $0\to \underline{Y}\to \underline{H}\to \underline{X}\to 0$. This is canonically a $\KI$-vector space, with addition given by the \emph{Baer sum} and scalar multiplication given by pulling back sequences along $a\cdot \id:\underline{X}\to \underline{X}$ for $a\in \KI$ (equivalently, pushing out along $a\cdot \id:\underline{Y}\to \underline{Y}$). 

These extension modules can be described as follows. We write $\underline{X}=(X,\fq_{X})$ and $\underline{Y}=(Y,\fq_{Y})$. Given a $\FI^s(\!(\fj)\!)$-linear morphism $f:X_{\FI^s(\!(\fj)\!)}\to Y_{\FI^s(\!(\fj)\!)}$, we consider the object $\underline{E}_f$ whose underlying space is the direct sum $Y\oplus X$ and whose Hodge-Pink lattice $\fq_f$ is given by 
\begin{equation}\label{eq:def-Ef}
\fq_f=\{(q_y+f(q_x),q_x)~|~(q_y,q_x)\in \fq_{Y}\oplus \fq_{X}\}.
\end{equation}
The lattice $\fq_f$ is an extension of $\fq_X$ by $\fq_Y$, hence the short exact sequence $S_f:0\to \underline{Y}\to \underline{E}_f\to \underline{X}\to 0$ is strict (\cf Remark \ref{rem:strict-exact-sequences}). Consequently, the assignment $f\mapsto [S_f]$ defines a map
\begin{equation}\label{eq:varphi}
\varphi:\Hom_{\FI^s(\!(\fj)\!)}\left(X_{\FI^s(\!(\fj)\!)},Y_{\FI^s(\!(\fj)\!)}\right)\longrightarrow \Ext^1_{\HPk}(\underline{X},\underline{Y}).
\end{equation}
It is classical to check that the Baer sum of $[S_f]$ and $[S_g]$ is $[S_{f+g}]$, and that the pullback of $[S_f]$ by $a\cdot \id:\underline{X}\to \underline{X}$, for $a\in \KI$, is $[S_{af}]$. In particular, $\varphi$ is a $\KI$-linear map. The following proposition suffices to describe the extension module:

\begin{Proposition}\label{prop:explicit-extensions-pre-Hodge-Pink}
The morphism $\varphi$ is surjective, with kernel
\[
\Hom_{\KI}(X,Y)+\Hom_{\FI^s[\![\fj]\!]}(\fq_X,\fq_Y).
\] 
\end{Proposition}
\begin{proof}
To show that $\varphi$ is surjective, fix an extension $[S]$ of $\underline{X}$ by $\underline{Y}$. The underlying sequence of $\KI$-vector spaces splits, and consequently $S$ is equivalent to an exact sequence of the form
\begin{equation}\label{eq:exactsequence-form1}
0\longrightarrow \underline{Y}\longrightarrow \left(Y\oplus X, \fq\right)\longrightarrow \underline{X}\longrightarrow 0
\end{equation}
where $\fq\subset Y_{\FI^s(\!(\fj)\!)}\oplus X_{\FI^s(\!(\fj)\!)}$ is a $\FI^s[\![\fj]\!]$-lattice. Let us define a linear map $f:X_{\FI^s(\!(\fj)\!)}\to Y_{\FI^s(\!(\fj)\!)}$ as follows. Since the morphisms of \eqref{eq:exactsequence-form1} are strict, the induced sequence of lattices is exact:
\[
0\longrightarrow \fq_{Y} \longrightarrow \fq \longrightarrow \fq_{X} \longrightarrow 0.
\]
For $q_x\in \fq_{X}$, choose a lift $(\tilde{y}_x,q_x)\in\fq$. Its $Y$-coordinate $\tilde{y}_x$ is unique modulo $\fq_Y$. The assignment $\bar{f}(q_x):=\tilde{y}_x+\fq_Y$ defines a $\FI^s[\![\fj]\!]$-linear morphism $\bar{f}:\fq_{X}\to Y_{\FI^s(\!(\fj)\!)}/\fq_{Y}$. Since $\fq_{X}$ is a projective $\FI^s[\![\fj]\!]$-module, $\bar{f}$ lifts to a $\FI^s(\!(\fj)\!)$-linear morphism $f:X_{\FI^s(\!(\fj)\!)}\to Y_{\FI^s(\!(\fj)\!)}$. It is clear from the construction that $S$ is isomorphic to $S_f$, hence $[S]=\varphi(f)$.

It remains to describe the kernel of $\varphi$. Observe that $S_f$ splits if and only if $f$ admits a representative that preserves the Hodge-Pink lattices; \ie if and only if there exists a $\KI$-linear map $u:X\to Y$ such that $f+u$ preserves the Hodge-Pink lattices. This proves that
\[
\ker \varphi=\Hom_{\KI}(X,Y)+\Hom_{\FI^s[\![\fj]\!]}(\fq_X,\fq_Y)
\]
as desired.
\end{proof}

Let $\mathbbm{1}$ be the Hodge-Pink structure whose underlying space is $\KI$ itself and whose Hodge-Pink lattice is $\fq_{\mathbbm{1}}=\fp_{\mathbbm{1}}=\FI^s[\![\fj]\!]$. We call $\mathbbm{1}$ the \emph{neutral Hodge-Pink structure}. We end this paragraph with the following corollary:

\begin{Corollary}\label{cor:thm-absolute-Hodge-Pink}
Let $\underline{H}$ be a Hodge-Pink structure. We have a natural isomorphism of $\KI$-vector spaces
\[
\varphi:\frac{H_{\FI^s(\!(\fj)\!)}}{H+\fq_H} \stackrel{\sim}{\longrightarrow} \Ext^1_{\HPk}(\mathbbm{1},\underline{H}).
\]
\end{Corollary}

\subsubsection*{Hodge-additive extensions}
In contrast to the number field setting, the space $\Ext^1_{\HPk}(\underline{X},\underline{Y})$ is almost never finite-dimensional. This issue is closely related to the problem of defining the \emph{right} regulator, which in the classical setting is a morphism between finite-dimensional vector spaces. Following Pink, we now discuss the notion of \emph{Hodge additivity} for extensions, which is a step toward addressing this difficulty. 

Let $\underline{H}$ be a Hodge-Pink structure. We first recall how to associate a finite decreasing filtration---the \emph{Hodge filtration}---to $H_{\FI^s}$. For $p\in \bZ$, let $\operatorname{Fil}^pH_{\FI^s}$ denote the image of $\fp_H\cap \fj^p \fq_H$ under the composition:
\begin{equation}
\begin{tikzcd}[column sep=3em]
\fp_H=H\otimes_{\KI,\nu} \FI^s[\![\fj]\!] \arrow[r,"\operatorname{mod}~\fj"] & H\otimes_{\KI}\FI^s=:H_{\FI^s}.
\end{tikzcd}\nonumber
\end{equation}
We call $\Fil^{\bullet}H_{\FI^s}=(\operatorname{Fil}^pH_{\FI^s})_p$ the \emph{Hodge filtration} of $\underline{H}$. Note that the data of $(H,\Fil^{\bullet}H_{\FI^s})$ defines a function field Hodge structure in the sense of Definition \ref{def:hodge-structure}, and that the assignment $\underline{H}\mapsto (H,\Fil^{\bullet}H_{\FI^s})$ is functorial. 

\begin{Definition}[{\cite[\S 8]{pink}}, Hodge additivity]\label{def:hodge-additive}
A strict exact sequence $0\to \underline{H}'\to \underline{H}\to \underline{H}''\to 0$ is said to be \emph{Hodge additive} if, for all $p\in \bZ$, it induces an exact sequence of $\FI^s$-vector spaces
\[
0\longrightarrow \Fil^p H'_{\FI^s}\longrightarrow \Fil^p H_{\FI^s}\longrightarrow \Fil^p H''_{\FI^s}\longrightarrow 0.
\]
\end{Definition}

Note that the property of being Hodge additive is invariant under equivalences of extensions. We make the following definition.
\begin{Definition}[$\Ext^{1,\operatorname{ha}}_{\HPk}$]
Let $\underline{X}$, $\underline{Y}$ be objects in $\HPk$. We denote by $\Ext^{1,\operatorname{ha}}_{\HPk}(\underline{X},\underline{Y})$ the subset of Hodge-additive extensions of $\underline{X}$ by $\underline{Y}$ in $\HPk$.
\end{Definition} 

The group $\Ext^{1,\operatorname{ha}}_{\HPk}(\underline{X},\underline{Y})$ carries the structure of a $\KI$-vector space. The following criterion for Hodge additivity is due to Pink \cite[Proposition 8.7]{pink}.
\begin{Proposition}
Let $\underline{X}$, $\underline{Y}$ be objects in $\HPk$. Let $f:X_{\FI^s(\!(\fj)\!)} \to Y_{\FI^s(\!(\fj)\!)}$ be a $\FI^s(\!(\fj)\!)$-linear map. Then the sequence $S_f$ is Hodge additive if and only if
\[
f\in \Hom_{\FI^s[\![\fj]\!]}(\fp_X,\fp_Y)+\Hom_{\FI^s[\![\fj]\!]}(\fq_X,\fq_Y).
\]
\end{Proposition}

In the particular case where $\underline{X}$ is $\mathbbm{1}$, we get:
\begin{Corollary}\label{cor:hodge-additive-ext-neutral}
Let $\underline{H}$ be a Hodge-Pink structure. The morphism $\varphi$ of Corollary \ref{cor:thm-absolute-Hodge-Pink} induces:
\[
\varphi:\frac{\fp_H}{H+(\fq_H\cap \fp_H)} \stackrel{\sim}{\longrightarrow} \Ext^{1,\operatorname{ha}}_{\HPk}(\mathbbm{1},\underline{H}).
\]
\end{Corollary}

\subsection{Infinite Frobenii}\label{subsec:inf-frob}
\subsubsection*{The classical picture}
Before introducing \emph{infinite Frobenii} for Hodge-Pink structures, let us briefly recall the classical story. According to Nekov\'a\v{r} \cite[(2.4)]{nekovar} and Deligne \cite[\S 1.4 (M7)]{deligne-droite}, an \emph{infinite Frobenius} $\phi_{\infty}$ for a Hodge structure $(H,\operatorname{Fil}H_\bC)$ (with coefficients in $\bR$ and base $\bC$) is an involution of the $\bR$-vector space $H$ such that $\phi_{\infty}\otimes_{\bR} c$ preserves $\operatorname{Fil}H_\bC$, where $c$ denotes complex conjugation. Hodge structures arising from the singular cohomology groups of a variety $X$ over $\bR$ are naturally equipped with an infinite Frobenius, induced functorially by the action of complex conjugation on the complex points $X(\bC)$. 

We let $\Hdg_{\bR}^+$ be the category whose objects are pairs $(\underline{H},\phi_{\infty})$, where $\underline{H}$ is a Hodge structure and $\phi_{\infty}$ is an infinite Frobenius for $\underline{H}$. Morphisms in $\Hdg_{\bR}^+$ are morphisms of mixed Hodge structures over $\bC$ which commute with the infinite Frobenii. 

Extension modules in the category $\Hdg_{\bC}$ of mixed Hodge structures over $\bC$ are well known. Given an object $\underline{H}$ of $\Hdg_{\bC}$, the complex of $\bR$-vector spaces
\begin{equation}
\left[H\oplus F^0H_{\bC}\stackrel{\left(\begin{smallmatrix} 1 \\ \text{-}1 \end{smallmatrix}\right)}{\longrightarrow}H_{\bC}\right] \nonumber
\end{equation}
represents the cohomology of $\RHom_{\Hdg_{\bC}}(\mathbbm{1},\underline{H})$ (e.g. \cite[\S 1]{beilinson}, \cite[Proposition 2]{carlson}, \cite[Theorem 3.31]{peters}). We obtain an $\bR$-linear morphism
\begin{equation}\label{eq:isom-classical-exthodge}
\frac{H_{\bC}}{H+F^0H_{\bC}}\stackrel{\sim}{\longrightarrow} \Ext^1_{\Hdg_{\bC}}(\mathbbm{1},\underline{H}).
\end{equation}
If now $\underline{H}^+$ denotes an object in the category $\Hdg^+_{\bR}$ with infinite Frobenius $\phi_{\infty}$, the complex $\RHom_{\Hdg_{\bR}^+}(\mathbbm{1}^+,\underline{H}^+)$ is instead represented by 
\begin{equation}
\left[H^+\oplus (F^0H_{\bC})^+\stackrel{\left(\begin{smallmatrix} 1 \\ \text{-}1 \end{smallmatrix}\right)}{\longrightarrow}(H_{\bC})^+\right] \nonumber
\end{equation}
where the superscript $+$ denotes the corresponding $\bR$-subspace fixed by $\phi_{\infty}\otimes c$ (e.g. \cite[\S 1]{beilinson}, \cite[(2.5)]{nekovar}). We obtain an $\bR$-linear morphism
\begin{equation}\label{eq:isom-classical-exthodge+}
\frac{(H_{\bC})^+}{(H)^++(F^0H_{\bC})^+}\stackrel{\sim}{\longrightarrow} \Ext^1_{\Hdg_{\bR}^+}(\mathbbm{1}^+,\underline{H}^+). 
\end{equation}

\subsubsection*{Infinite Frobenii for Hodge-Pink structures}
We keep the notation from the previous subsection. In this subsection, we enrich Hodge-Pink structures with a compatible continuous action of the profinite Galois group $G_\infty:=\Gal(\FI^s|\FI)$. In favourable cases, we compute the corresponding extension spaces. 

Let $\underline{H}=(H,\fq_H)$ be a Hodge-Pink structure.
\begin{Definition}[Infinite Frobenius]\label{def:infinite-frob}
An \emph{infinite Frobenius} for $\underline{H}$ is a $\KI$-linear continuous representation $\phi:G_\infty\to \Aut_{\KI}(H)$, where $G_\infty$ carries the profinite topology and $H$ the discrete topology, such that, for all $\sigma\in G_\infty$,
\[
\phi(\sigma)\otimes_{\KI} \sigma:H_{\FI^s(\!(\fj)\!)}\longrightarrow H_{\FI^s(\!(\fj)\!)}
\]
preserves the Hodge-Pink lattice $\fq_H$. 

Above, we denoted by $\sigma$ its extension to $\FI^s(\!(\fj)\!)$, \ie the one obtained by functoriality of the assignment $L\mapsto L(\!(\fj)\!)$ on separable extensions $L$ of $\FI$, applied to $\sigma:\FI^s\to \FI^s$.
\end{Definition}

We let $\HPk^+$ denote the category whose objects are pairs $(\underline{H},\phi_{H})$, where $\underline{H}$ is a Hodge-Pink structure and $\phi_{H}$ is an infinite Frobenius for $\underline{H}$. We extend Definition \ref{def:strict-sequence-Hdg} to this setting. 
\begin{Definition}[Strict exact sequence]\label{def:strict-sequence-HdgFrob}
A sequence $S:0\to \underline{H}'\to \underline{H}\to \underline{H}''\to 0$ in $\HPk^+$ is called \emph{exact} (resp. \emph{strict}) if the underlying sequence in $\HPk$ is exact (resp. strict) as in Definition \ref{def:strict-sequence-Hdg}.
\end{Definition}

The category $\HPk^+$, equipped with the above class of strict exact sequences, is an exact category, and we now study its extension modules. While the analogue of \eqref{eq:isom-classical-exthodge} holds for Hodge\nobreakdash-Pink structures (this was Corollary \ref{cor:hodge-additive-ext-neutral}), a description analogous to \eqref{eq:isom-classical-exthodge+} does not hold in our setting, since the action of complex conjugation is replaced by that of the (infinite) profinite group $G_\infty$. Consequently, the extension modules are intertwined with the Galois cohomology of $G_\infty$, preventing an isomorphism as simple as \eqref{eq:isom-classical-exthodge+} from existing. We next clarify how Galois cohomology interferes with the computation of extension spaces. 

By definition, we have a forgetful functor from $\HPk^+$ to the category of $\KI$-linear continuous representations of $G_\infty$, which assigns to $\underline{H}^+=(\underline{H},\phi_H)$ the representation $\phi_H$. Being exact, it induces a natural $\KI$-linear morphism on extensions:
\begin{equation}\label{eq:forgetful-on-ext}
\Ext^1_{\HPk^+}(\underline{G}^+,\underline{H}^+)\longrightarrow \Ext^1_{G_\infty}(\phi_G,\phi_H).
\end{equation}
The target is the module of extensions of continuous $\KI$-linear $G_\infty$-representations. 
\begin{Definition}[Analytic reduction]\label{def:analytic-reduction-HP}
An extension $[\underline{E}]\in \Ext^1_{\HPk^+}(\underline{G}^+,\underline{H}^+)$ is said to have \emph{analytic reduction} if the induced extension of continuous $G_\infty$-modules splits; \ie if it lies in the kernel of \eqref{eq:forgetful-on-ext}. 

We denote by $\Ext^1_{\HPk^+,\infty}(\underline{G}^+,\underline{H}^+)$ (or simply by $\Ext^1_{\infty}(\underline{G}^+,\underline{H}^+)$) the module they form modulo equivalences.
\end{Definition}

Let $\mathbbm{1}^+$ denote the object of $\HPk^+$ given by the pair $(\mathbbm{1},\phi_{\mathbbm{1}})$, where $\phi_{\mathbbm{1}}:G_\infty\to \KI^{\times}$, $\sigma\mapsto 1$, is the trivial representation.

\begin{Definition}\label{def:dH}
Given an object $\underline{H}^+$ of $\HPk^+$, we denote by 
\[
d_{\underline{H}^+}:\Ext^1_{\HPk^+}(\mathbbm{1}^+,\underline{H}^+)\longrightarrow \operatorname{H}^1(G_\infty,H)
\]
the $\KI$-linear morphism \eqref{eq:forgetful-on-ext} with $\underline{G}^+=\mathbbm{1}^+$, landing in continuous group cohomology. 
\end{Definition}

Recall that $H_{\FI^s(\!(\fj)\!)}=H\otimes_{\KI,\nu} \FI^s(\!(\fj)\!)$ is endowed with a continuous action of $G_\infty$, where $\sigma\in G_\infty$ acts by $\phi_H(\sigma)\otimes \sigma$. For $S$ a subset of $H_{\FI^s(\!(\fj)\!)}$, we denote by $S^+\subseteq S$ the subset of elements fixed under this action of $G_\infty$. Let $\varphi$ be the isomorphism of Corollary \ref{cor:thm-absolute-Hodge-Pink}. There is a $\KI$-linear morphism
\[
\varphi^+: \frac{(H_{\FI^s(\!(\fj)\!)})^+}{H^{+}+\fq_H^+}\longrightarrow \Ext^1_{\HPk^+}(\mathbbm{1}^+,\underline{H}^+),
\]
sending the class of $h\in (H_{\FI^s(\!(\fj)\!)})^+$ to the extension $(\varphi(h),\left(\begin{smallmatrix}\phi_H & 0 \\ 0 & 1 \end{smallmatrix}\right))$ with split infinite Frobenius. By definition, the image of $\varphi^+$ lands in $\Ext^1_{\infty}(\mathbbm{1}^+,\underline{H}^+)$. Under an additional assumption, we can say more.
\begin{Proposition}\label{prop:ext-+}
Suppose that $\operatorname{H}^1(G_\infty,\fq_H)$ is trivial. Then the sequence of $\KI$-vector spaces
\[
0\longrightarrow \frac{(H_{\FI^s(\!(\fj)\!)})^+}{H^{+}+\fq_H^+} \xrightarrow{\varphi^+} \Ext^1_{\HPk^+}(\mathbbm{1}^+,\underline{H}^+) \xrightarrow{d_{\underline{H}^+}} \operatorname{H}^1(G_\infty,H)
\]
is exact. In particular, under the same assumption, $\varphi^+$ induces an isomorphism of $\KI$-vector spaces:
\[
\varphi^+:\frac{(H_{\FI^s(\!(\fj)\!)})^+}{H^{+}+\fq_H^+} \stackrel{\sim}{\longrightarrow} \Ext^1_{\infty}(\mathbbm{1}^+,\underline{H}^+).
\]
\end{Proposition}

\begin{Remark}
It will appear that the condition $\operatorname{H}^1(G_\infty,\fq_H)=(0)$ is always satisfied for Hodge\nobreakdash-Pink structures arising from rigid analytically trivial $A$-motives. We refer to Lemma \ref{lem:equivariant-lattices-hodge} below. 
\end{Remark}

\begin{proof}
We begin with some notation and an observation. Recall that $\varphi(h)$ is represented by the extension of $\mathbbm{1}$ by $\underline{H}$ whose middle term is
\[
\underline{E}_h:=\left(H\oplus \KI, \begin{pmatrix}
1 & h \\ 0 & 1
\end{pmatrix} \fq_H\oplus \FI^s[\![\fj]\!] \right).
\]
Given a continuous cocycle $c:G_\infty\to H$, we denote by $[c]$ the continuous $\KI$-linear $G_\infty$\nobreakdash-representation of $H\oplus \KI$ given by
\begin{equation}
[c]:G_\infty\longrightarrow \End_{\KI}(H\oplus \KI), \quad \sigma \longmapsto \begin{pmatrix}\phi_{H}(\sigma) & c(\sigma) \\ 0 & 1 \end{pmatrix}. \nonumber
\end{equation}

Now assume that $[c]$ defines an infinite Frobenius for the Hodge-Pink structure $\underline{E}_h$. Then $(\underline{E}_h,[c])$ defines an extension of $\mathbbm{1}^+$ by $\underline{H}^+$ in $\HPk^+$. For $m\in H$, the diagram
\[
\begin{tikzcd}[row sep=2.5em,ampersand replacement = \&]
0\arrow[r] \& {\underline{H}^+}  \arrow[r]\arrow[d,"{\id_{\underline{H}}}"] \& \left(\underline{E}_h,[c]\right) \arrow[d,"{\left(\begin{smallmatrix}\id_H & m \\ 0 & 1 \end{smallmatrix}\right)}"]\arrow[r] \& \mathbbm{1}^+ \arrow[d,"1"]\arrow[r] \& 0 \\
0\arrow[r] \& \underline{H}^+ \arrow[r] \& \left(\underline{E}_{h+m},[\sigma \mapsto c(\sigma)+m-\phi_{H}(\sigma)(m)]\right) \arrow[r] \& \mathbbm{1}^+ \arrow[r] \& 0
\end{tikzcd}
\]
defines an equivalence in $\HPk^+$ between the extensions
\begin{equation}\label{eq:congruence-hodge}
\left(\underline{E}_h,[c]\right)\quad\text{and}\quad \left(\underline{E}_{h+m},[\sigma\mapsto c(\sigma)+m-\phi_{H}(\sigma)(m)]\right). 
\end{equation}

If $[\underline{E}^+]$ is an element of $\ker d_{\underline{H}^+}$, there exist $h'\in H_{\FI^s(\!(\fj)\!)}$ and $m\in H$ such that $[\underline{E}^+]$ is equivalent to an extension of the form $(\underline{E}_{h'},[\sigma\mapsto m-\phi_{H}(\sigma)(m)])$. By \eqref{eq:congruence-hodge}, we may assume without loss of generality that $[\underline{E}^+]$ is of the form $\left(\underline{E}_h,[0]\right)$. The condition that the infinite Frobenius of $\underline{E}^+$ preserves the Hodge-Pink lattice reads
\begin{equation}
\text{for~all}~ \sigma\in G_\infty,\quad (\phi_{H}(\sigma)\otimes \sigma)(h)-h\in \fq_H. \nonumber
\end{equation}
In particular, $h+\fq_H$ is invariant under $G_\infty$ as an element of the quotient module $H_{\FI^s(\!(\fj)\!)}/\fq_H$. Under our assumption that $\operatorname{H}^1(G_\infty,\fq_H)$ is trivial, we obtain
\[
\left(\frac{H_{\FI^s(\!(\fj)\!)}}{\fq_H}\right)^+=\frac{(H_{\FI^s(\!(\fj)\!)})^+}{\fq_H^+}.
\]
Hence, $h$ is congruent to some $h_0\in (H_{\FI^s(\!(\fj)\!)})^+$ modulo $\fq_H$. We have $\underline{E}_{h}=\underline{E}_{h_0}$ since their Hodge\nobreakdash-Pink lattices coincide. To conclude, it suffices to note that any equivalence in $\HPk^+$ between the extensions $[\underline{E}^+]=\left(\underline{E}_{h_0},[0]\right)$ and $\left(\underline{E}_{k_0},[0]\right)$ for some $k_0\in  H_{\FI^s(\!(\fj)\!)}$ implies $h_0-k_0\in H^++\mathfrak q_H^+$.
\end{proof}

We now combine Hodge additivity and analytic reduction. We denote by
\[
\Ext^{1,\operatorname{ha}}_{\infty}(\mathbbm{1}^+,\underline{H}^+):=\Ext^{1}_{\infty}(\mathbbm{1}^+,\underline{H}^+) \cap \Ext^{1,\operatorname{ha}}_{\HPk^+}(\mathbbm{1}^+,\underline{H}^+)
\]
the subspace of Hodge-additive extensions having analytic reduction. The same argument as in the proof of Proposition \ref{prop:ext-+} applies to show:
\begin{Proposition}\label{prop:ext-+-ha}
Suppose that $\operatorname{H}^1(G_\infty,\fq_H\cap \fp_H)$ is trivial. Then the sequence of $\KI$-vector spaces
\[
0\longrightarrow \frac{\fp_H^+}{H^{+}+(\fq_H\cap\fp_H)^+} \xrightarrow{\varphi^+} \Ext^{1,\operatorname{ha}}_{\HPk^+}(\mathbbm{1}^+,\underline{H}^+) \xrightarrow{d_{\underline{H}^+}} \operatorname{H}^1(G_\infty,H)
\]
is exact. In particular, under the same assumption, $\varphi^+$ induces an isomorphism of $\KI$-vector spaces:
\[
\varphi^+:\frac{\fp_H^+}{H^{+}+(\fq_H\cap \fp_H)^+} \stackrel{\sim}{\longrightarrow} \Ext^{1,\operatorname{ha}}_{\infty}(\mathbbm{1}^+,\underline{H}^+).
\]
\end{Proposition}

\section{Rigid analytically trivial $A$-motives}\label{sec:rat-mix-A-mot}
The chief aim of this section is to associate a Hodge-Pink structure to a rigid analytically trivial $A$-motive (Subsection \ref{subsec:asso-MHPS}). Note that this was already carried out by Hartl--Juschka in \cite{hartl-juschka}, under the assumption that $\infty$ is a degree-one point of $C$.

In Subsection \ref{subsec:betti}, we review the construction of the \emph{Betti realization} $\Betti{\underline{M}}$ of an $A$-motive $\underline{M}$ and the notion of rigid analytic triviality. We show that the module $\underline{M}_B$ carries a canonical continuous action of the profinite group $G_{\infty}$. This will be used later to define the infinite Frobenii on the associated Hodge-Pink structure. In Subsection \ref{subsection:analytic continuation}, we show that elements of the Betti realization of a rigid analytically trivial $A$-motive can be interpreted as rigid analytic functions, and we prove that they admit a meromorphic continuation to the rigid analytification of $\Spec A\otimes\CI$. This allows us to take Laurent expansions of elements of $\underline{M}_B$ at $\fj$ and thereby define the Hodge-Pink lattice. 

\subsection{Definitions}\label{subsec:def}
In this subsection we review the usual setup of $A$-motives. Recall that $(C,\cO_C)$ is a geometrically irreducible smooth projective curve over a finite field $\bF$ with $q$ elements, that $\infty$ is a closed point of $C$, and that $A$ is the $\bF$-algebra $\cO_C(C\setminus \{\infty\})$. 

Let $R$ be an $A$-algebra via an $\bF$-algebra morphism $\kappa:A\to R$. We denote by $A\otimes R$ the tensor product over $\bF$, and we let $\fj=\fj_\kappa$ be the kernel of the multiplication map $A\otimes R\to R$, $a\otimes r\mapsto \kappa(a)r$. The following observation appears in \cite{hartl-isogeny}.
\begin{Lemma}\label{lem:j-Cartier}
The ideal $\fj$ is a locally free $A\otimes R$-module of rank $1$. In particular, $V(\fj)$ defines a Cartier divisor on $\Spec A\otimes R$.
\end{Lemma}
\begin{proof}
Denote by $\fj:=\fj_{\id}$ the kernel of the multiplication map $A\otimes A\to A$. Being locally free of rank $1$ is a property stable under base change, and as $\fj_\kappa=\fj_{\id}\otimes_{A\otimes A}(A\otimes R)$, it suffices to show that $\fj:=\fj_{\id}$ is locally free of rank $1$ over $A\otimes A$.

Let $\mathfrak{P}\subset A\otimes A$ be a prime ideal and denote by $\mathfrak{j}_{\mathfrak{P}}$ the localization of $\mathfrak{j}$ at $\mathfrak{P}$. We aim to prove that $\mathfrak{j}_{\mathfrak{P}}$ is a principal ideal. There are two situations:
\begin{itemize}
\item If $\mathfrak{j}\nsubseteq \mathfrak{P}$, then there exists an element $f\in \mathfrak{j}$ with $f\notin \mathfrak{P}$, hence $1=f^{-1}\cdot f\in \mathfrak{j}_{\mathfrak{P}}\subset (A\otimes A)_{\mathfrak{P}}$, and thus $\mathfrak{j}_{\mathfrak{P}}= (A\otimes A)_{\mathfrak{P}}$.
\item If $\mathfrak{j}\subseteq \mathfrak{P}$, let $\mathfrak{p}:=\mathfrak{P}/\mathfrak{j}\subset A\otimes A/\mathfrak{j}\cong A$; this defines a prime ideal of $A$. Let $j\in \mathfrak{j}_{\mathfrak{P}}$ be an element whose class generates the $A_{\mathfrak{p}}$-module
\[
\mathfrak{j}_{\mathfrak{P}}/\mathfrak{j}_{\mathfrak{P}}^2\cong (\Omega_{A/\mathbb{F}}^1)_{\mathfrak{p}}
\]
(because $C$ is smooth over $\bF$, $\Omega_{A/\mathbb{F}}^1$ is locally free of rank $1$ over $A$). Then $\mathfrak{j}_{\mathfrak{P}}=(j)+\mathfrak{j}_{\mathfrak{P}}^2$ and Nakayama's lemma implies that $\mathfrak{j}_{\mathfrak{P}}=(j)$.
\end{itemize}
\end{proof}

Consequently, the ideal $\fj\subset A\otimes R$ is invertible. Given an $A\otimes R$-module $N$, we define
\[
N[\fj^{-1}]:=\mathrm{colimit} \left(N\longrightarrow N\otimes_{A\otimes R} \fj^{-1} \longrightarrow N\otimes_{A\otimes R} \fj^{-2}\longrightarrow \cdots \right).
\]
Let $\tau$ be the ring endomorphism of $A\otimes R$ acting as the identity on $A$ and as raising to the $q$th power on $R$. 
\begin{Definition}[$A$-motives]\label{def:A-motive}
An \emph{$A$-motive of rank $r$ over $R$} is a pair $\underline{M}=(M,\tau_M)$ where $M$ is a locally free module of constant rank $r$ over $A\otimes R$ together with an $A\otimes R$-linear isomorphism:
\[
\tau_M:(\tau^*M)[\fj^{-1}]\stackrel{\sim}{\longrightarrow} M[\fj^{-1}].
\]
We call $\underline{M}$ \emph{effective} whenever $\tau_M$ factors as a morphism $\tau^* M\to M$.

A morphism $(M,\tau_M)\to (N,\tau_N)$ of $A$-motives over $R$ is an $A\otimes R$-linear morphism $f:M\to N$ such that $f\circ \tau_M=\tau_N\circ \tau^*f$. We let $\AMot_R$ denote the $A$-linear category of $A$-motives over $R$.
\end{Definition}

We will be mainly interested in the situation where $R=F$ is a finite field extension of $K$, where $K$ is the field of fractions of $A$ and $\kappa:A\to F$ is the canonical inclusion. It is often convenient to assume that $F$ is totally ramified at $\infty$, as then $F_{\infty}:=F\otimes_{K} \KI$ is a finite field extension of $\KI$. A reduction is possible under the following lemma\footnote{The decomposition $K\subset L\subset F$ where $L$ is separable over $K$ and $F$ purely inseparable over $L$ is well-known, see \href{https://stacks.math.columbia.edu/tag/030K}{[030K]}. The other way round works in our situation of transcendence degree one.}:
\begin{Lemma}\label{lem:insep-sep}
Let $F$ be a finite extension of $K$ and let $K^{\natural}$ denote the perfection of $K$ inside $F$. Then $K^{\natural}$ is purely inseparable over $K$ and $F$ is separable over $K^{\natural}$. In addition $[F:K^{\natural}]$ equals the separable degree of $F$ over $K$ and $[K^{\natural}:K]$ the inseparable degree.
\end{Lemma}
\begin{proof}
Recall that $K^{\natural}$ is the subfield of elements $x$ in $F$ for which there exists $k\geq 0$ such that $x^{p^k}\in K$. Clearly, $K^{\natural}$ is purely inseparable over $K$, so it suffices to show that $F$ is separable over~$K^{\natural}$. 

We write the cotangent exact sequence \href{https://stacks.math.columbia.edu/tag/00RS}{[00RS]} with respect to the inclusions $\bF\subset K^{\natural}\subset F$:
\[
F\otimes_{K^{\natural}} \Omega_{K^{\natural}/\bF}^1\longrightarrow \Omega_{F/\bF}^1 \longrightarrow \Omega_{F/K^{\natural}}^1\longrightarrow 0.
\]
It suffices to show that $\Omega_{F/K^{\natural}}^1=(0)$, hence that the first map, between one-dimensional $F$-vector spaces, is nonzero. Assume otherwise: let $a\in K^{\natural}$ be a separating element; \ie $K^{\natural}/\bF(a)$ is finite separable. Then $da\neq 0$ but $da=0$ in $\Omega_{F/\bF}^1$ by assumption. By \href{https://stacks.math.columbia.edu/tag/031W}{[031W]}, there exists $b\in F$ such that $a=b^p$. This implies $b\in K^{\natural}$ and then $da=0$ in $\Omega_{K^{\natural}/\bF}^1$, a contradiction. 

For the assertion on the degree, let $\Omega$ be an algebraically closed field containing $K$. Note that any $K$-linear embedding $\sigma:F\to \Omega$ fixes $K^{\natural}$. In particular,
\[
\Hom_K(F,\Omega)=\Hom_{K^{\natural}}(F,\Omega).
\]
Because $F$ is separable over $K^{\natural}$, comparing cardinalities yields $[F:K]_{\mathrm{sep}}=[F:K^{\natural}]$ and hence $[K^{\natural}:K]=[F:K]_{\mathrm{insep}}$.
\end{proof}

The following ``trick'' allows us to reduce to the situation where $F$ is totally ramified at $\infty$ over $K$.
\begin{Reduction}[to the totally ramified case]\label{red:F-totally-ramified}
If $R\to S$ is a finite \'etale $A$-algebra morphism, there is an exact \emph{restriction functor} $\mathrm{Res}_{S/R}$ from the category of $A$-motives over $S$ to that over $R$ which is right-adjoint to the base-change functor $\underline{M}\mapsto \underline{M}_S$ (\eg \cite[\S 2]{gazda}).

Let $\underline{M}$ be an $A$-motive over $F$. Consider the inclusions $K\subset K^{\natural} \subset F$ given by Lemma \ref{lem:insep-sep}. Because $F$ is finite separable over $K^{\natural}$, $\mathrm{Res}_{F/K^{\natural}} \underline{M}$ is an $A$-motive over $K^{\natural}$ (note that the naive ``$\mathrm{Res}_{F/K}$'' does not preserve $A$-motives in the sense of Definition \ref{def:A-motive}). We have canonical isomorphisms
\[
\Hom_{\AMot_F}(\mathbbm{1}_F,\underline{M}) \cong \Hom_{\AMot_{K^{\natural}}}(\mathbbm{1}_{K^{\natural}},\mathrm{Res}_{F/K^{\natural}}\underline{M}),\quad \Ext^1_{\AMot_F}(\mathbbm{1}_F,\underline{M}) \cong \Ext^1_{\AMot_{K^{\natural}}}(\mathbbm{1}_{K^{\natural}},\mathrm{Res}_{F/K^{\natural}}\underline{M})
\]
which allow us to reduce the study of morphisms and extension modules to the case where $F=K^{\natural}$. 

From now on, we will therefore assume that $F$ is a finite field extension of $K$ which is totally ramified at $\infty$ (we could also assume that $F$ is purely inseparable over $K$, although this does not seem to simplify our arguments further).
\end{Reduction}

\begin{Notation}[The $t$-setting]\label{not:t-setting}
The $t$-setting refers to the following situation: $C$ is the curve $\bP^1_{\bF}$ and $\infty$ denotes the closed point $[0:1]$. The variable ``$t$'' denotes the affine coordinate on $\bP^1_{\bF}\setminus\{\infty\}=\Spec\bF[t]$, so that $A=\bF[t]$. We write $F[t]$ instead of the tensor product $A\otimes F$, $t$ for $t\otimes 1$ and $\theta$ for $1\otimes t$. The ideal $\fj$ of $F[t]$ is principal, generated by $t-\theta$. 

In this situation, the typical examples of $A$-motives (also referred to as \emph{$t$-motives} in Anderson's terminology) are the Carlitz twists $\underline{A}(n)$, $n\in \bZ$. Their underlying $F[t]$-module is $F[t]$ itself, and $\tau_{A(n)}$ acts by semilinear multiplication by $(t-\theta)^{-n}$ (the minus-sign convention is chosen to fit the analogy with the classical Tate twists $\bZ(n)$, \eg for weight reasons). 
\end{Notation}

For $A$-motives over $F$, one has a notion of weights that we briefly recall (the reader will find all details in \cite[\S 3]{gazda}). There is an (exact) functor 
\[
\AMot_F\longrightarrow {}_{\infty}\mathbf{Isoc}_F,\quad \underline{M} \longmapsto \cI_{\infty}({\underline{M}})
\]
where the target category is that of $\infty$-isocrystals over $F$. The category ${}_{\infty}\mathbf{Isoc}_F$ admits a \emph{slope function} in the sense of \cite{andre}, hence any object admits a unique slope filtration (the Harder--Narasimhan filtration) \cite[Theorem 3.10]{gazda}. We say that $\underline{M}$ has weights $\nu_1<\nu_2<\cdots<\nu_s$ if $\cI_\infty(\underline{M})$ has slopes $-\nu_1>-\nu_2>\cdots>-\nu_s$.

\subsection{The Betti realization functor}\label{subsec:betti}
We let $K$ be the fraction field of $A$ (equivalently, the function field of $C$). Let $F$ be a finite extension of $K$ which is totally ramified at $\infty$ (see Reduction \ref{red:F-totally-ramified}). We denote by $F_{\infty}$ its completion with respect to the unique place above $\infty$. Here, we introduce the \emph{Betti realization} of an $A$-motive over $F$ (Definition \ref{def:betti-realization}) and discuss \emph{rigid analytic triviality} (Definition~\ref{def:rigid analytically trivial}). One goal is to define the full subcategory $\AMot_F^{\text{rat}}$ of $\AMot_F$ consisting of \emph{rigid analytically trivial $A$-motives}, which will be the source of the Hodge-Pink realization functor to be defined in Subsection \ref{subsec:asso-MHPS}. Historically, the notion of rigid analytic triviality dates back to Anderson \cite[\S 2]{anderson}, and most of this subsection owes to his work. A novelty of our account is the consideration of a natural continuous action of $G_\infty$---the absolute Galois group at $\infty$---on the Betti realization $A$-module. The existence of canonical infinite Frobenii attached to the associated Hodge-Pink structures will follow from this construction. 

\subsubsection*{Tate algebras}
Let $L$ be a field over $\bF$ complete with respect to a non-archimedean norm $|\cdot|$, and let $\cO_L$ be its valuation ring with maximal ideal $\fm_L$. Typically, $L$ will be a subfield of $\CI$ (the completion of a separable closure $\FI^s$ of $\FI$). 

Let $L$ be a non-archimedean Banach algebra containing $\bF$ with norm $|\cdot|$. Typically, $L$ is a subring of $\CI$ (the completion of a separable closure $\FI^s$ of $\FI$) equipped with its induced norm. We define a norm on $A\otimes L$ by setting
\[
\|x\|:=\inf\left(\max_i |l_i|\right) \quad \text{for~}x\in A\otimes L
\]
where the infimum is taken over all the representations of $x$ of the form $\sum_{i}{a_i\otimes l_i}$. This is a (multiplicative) norm \cite[Proposition~2.2]{gazda-maurischat}; we call it \emph{the Gauss norm}. 

\begin{Definition}\label{def:aff-Tate-algebra}
We denote by $L\langle A \rangle$ the $L$-algebra given by the completion of $A\otimes L$ with respect to the Gauss norm. We again denote by $\tau$ the continuous extension of $A\otimes L\to A\otimes L$, $a\otimes c\mapsto a\otimes c^q$, to $L\langle A \rangle$.
\end{Definition}

\begin{Remark}
The notation $L\langle A\rangle$ is meant to stress that it generalizes the classical Tate algebra over $L$. In the $t$-setting (\cf Notation \ref{not:t-setting}), $L\langle A\rangle\cong L\langle t \rangle$. For general $A$ and a complete field $L$, one can show that $L\langle A \rangle$ is an affinoid algebra (see \cite[Proposition 3.2]{gazda-maurischat2}). 
\end{Remark}

The following preliminary lemma will be used later, in the definition of the Betti realization functor.
\begin{Lemma}\label{lem:j-invertible-in-Tate-algebra}
Let $\kappa:A\to L$ be an $\bF$-algebra morphism with discrete unbounded image. We have $\fj_{\kappa} L\langle A\rangle=L\langle A \rangle$.
\end{Lemma}
\begin{proof}
Because $\kappa(A)$ is unbounded in $L$, it contains an element $\alpha$ of norm $|\alpha|>1$. Let $a\in A$ be such that $\alpha=\kappa(a)$. Then $\kappa(a)^{-1}$ has norm $<1$ and the series 
\begin{equation}
-\sum_{n\geq 0}{a^n\otimes \kappa(a)^{-(n+1)}}\nonumber
\end{equation}
converges in $L\langle A\rangle$ to the inverse of $a\otimes 1-1\otimes \kappa(a)$. 
\end{proof}

\subsubsection*{The Betti realization of an $A$-motive}
We denote by $\cO_{\FI}$ the valuation ring of $\FI$ with maximal ideal $\fm_{\FI}$. We fix a separable closure $\FI^s$ of $\FI$, and denote by $\CI$ its completion (which is algebraically closed and complete, by Krasner's lemma). The canonical norm on $\KI$ extends uniquely to $\FI$ and to a norm $|\cdot|$ on $\CI$. The action of $G_{\infty}=\Gal(\FI^s|\FI)$ extends continuously to $\CI$. 

Let $\underline{M}=(M,\tau_M)$ be an $A$-motive over $F$. By Lemma \ref{lem:j-invertible-in-Tate-algebra}, the ideal $\fj\subset A\otimes F$ becomes the unit ideal in $\CI\langle A\rangle$, and thus $\tau_M$ induces an isomorphism of $\CI\langle A \rangle$-modules:
\begin{equation}\label{eq:tau_M-v}
\tau_M:\tau^*(M\otimes_{A\otimes F}\CI\langle A \rangle)\stackrel{\sim}{\longrightarrow} M\otimes_{A\otimes F}\CI\langle A \rangle
\end{equation}
which commutes with the action of $G_\infty$ on $M\otimes_{A\otimes F}\CI\langle A \rangle$, inherited from its action on the right-hand factor of the tensor product. We still denote by $\tau_M$ the isomorphism \eqref{eq:tau_M-v}. 
\begin{Definition}\label{def:betti-realization}
The \emph{Betti realization of $\underline{M}$} is the $A$-module
\begin{equation}
\Betti{\underline{M}}:=\left\{\omega\in M\otimes_{A\otimes F}\bC_{\infty}\langle A \rangle~|~\omega=\tau_M(\tau^* \omega)\right\}. \nonumber
\end{equation}
It is endowed with the compatible action of $G_{\infty}$ that it inherits as a submodule of $M\otimes_{A\otimes F}\bC_{\infty}\langle A \rangle$. We let $\Betti{\underline{M}}^{+}$ denote the $A$-submodule of $\Betti{\underline{M}}$ consisting of elements fixed by $G_{\infty}$.
\end{Definition}

The next definition is borrowed from \cite[\S 2.3]{anderson}.
\begin{Definition}[Rigid analytically trivial]\label{def:rigid analytically trivial}
The $A$-motive $\underline{M}$ is called \emph{rigid analytically trivial} if the map
\[
\Betti{\underline{M}}\otimes_A \bC_{\infty}\langle A \rangle\longrightarrow M\otimes_{A\otimes F}\bC_{\infty}\langle A \rangle, \quad \omega\otimes c\longmapsto \omega\cdot c
\]
given by multiplication is bijective.
\end{Definition}

\begin{Remark}[Anderson--Thakur generating function]
In the $t$-setting (\cf Notation \ref{not:t-setting}), consider \emph{the Anderson--Thakur generating function}:
\[
\omega(t):=\eta \prod_{i=0}^{\infty}{\left(1-\frac{t}{\theta^{q^i}}\right)^{-1}},
\]
where $\eta$ is a $(q-1)$th root of $-\theta$ in $\CI$. The product converges in the Tate algebra $\CI\langle t \rangle$ to a unit, and $\omega(t)^n$ generates the $\bF[t]$-module $\underline{A}(n)_B$. In particular, Carlitz twists are rigid analytically trivial.
\end{Remark}

\begin{Remark}
Not every $A$-motive is rigid analytically trivial. An example of an $A$-motive that is not rigid analytically trivial is given in \cite[2.2]{anderson} or \cite[Ex. 3.2.10]{taelman}.
\end{Remark}

The following proposition rephrases \cite[Corollary~4.3]{hartl-boeckle}:
\begin{Proposition}\label{prop:rank-betti}
Let $\underline{M}$ be an $A$-motive over $F$ of rank $r$. Then $\Betti{\underline{M}}$ is a finite projective $A$-module of rank $r'$ satisfying $r'\leq r$, with equality if and only if $\underline{M}$ is rigid analytically trivial.
\end{Proposition}

When $\underline{M}$ is rigid analytically trivial, in Definition \ref{def:rigid analytically trivial} the field $\CI$ can be replaced by a much smaller field. This is the subject of the next proposition.
\begin{Proposition}\label{prop:finite-extension-betti}
Let $\underline{M}$ be a rigid analytically trivial $A$-motive over $F$. There exists a (complete) finite separable field extension $L$ of $\FI$ in $\CI$ such that $\Betti{\underline{M}}$ is contained in $M\otimes_{A\otimes F}L\langle A \rangle$. In particular, the action of $G_\infty$, equipped with the profinite topology, on $\Betti{\underline{M}}$, equipped with the discrete topology, is continuous.
\end{Proposition}
\begin{proof}
For effective $t$-motives this is proved in \cite[Theorem~4]{anderson}, so we explain how one may reduce to this case. Let $t$ be a non-constant element of $A$. The inclusion $\bF[t]\subset A$ makes $A$ into a finite flat $\bF[t]$-module, and therefore $\underline{M}$ defines a $t$-motive of rank $\deg(t)\cdot \operatorname{rank}\underline{M}$ over $F$ (\cf Notation \ref{not:t-setting}). Let $n>0$ be an integer such that $(t-\theta)^n\tau_M(\tau^*M)\subset M$. Let $\underline{N}:=\underline{A}(-n)$ denote the Carlitz $(-n)$th twist over $F$. We have 
\begin{equation}
\Betti{\underline{N}}=\omega(t)^{-n} \bF[t] \subset \KI(\eta)\langle t \rangle. \nonumber
\end{equation}
The $\bF[t]$-motive $\underline{N}$ has been chosen so that $\underline{M}\otimes \underline{N}$ is effective (see Definition \ref{def:A-motive}). We are thus in the range of application of \cite[Theorem~4]{anderson}, which states that there exists a finite extension $H$ of $F_\infty$ in $\CI$ such that 
\begin{equation}
\Betti{\underline{M}}\otimes_{\bF[t]} \Betti{\underline{N}}=\Betti{(\underline{M}\otimes \underline{N})}\subset (M\otimes_{F[t]} N)\otimes_{F[t]}H\langle t \rangle= M\otimes_{F[t]}H\langle t \rangle. \nonumber
\end{equation}
It follows that there exists a finite extension $L'$ of $F_\infty$ such that $\Betti{\underline{M}}\subset M\otimes_{F[t]}L'\langle t \rangle$ (\eg one can take $L':=H(\eta)$). 

We now show that one can choose $L'$ separable over $F_\infty$. We let $\FI^s\langle t \rangle$ denote the subring of $\CI\langle t \rangle$ whose coefficients belong to a finite separable extension of $F_{\infty}$ inside $\FI^s$. Note that $M\otimes_{F[t]}\FI^s\langle t \rangle$ is free of finite rank over $\FI^s\langle t \rangle$. Therefore, $(M\otimes_{F[t]}\FI^s\langle t \rangle)/(t^n)$ is a finite-dimensional $\FI^s$-vector space for all positive integers $n$. By Lang's isogeny theorem (e.g. \cite[Proposition 1.1]{katz}), the multiplication map
\begin{equation}
\left\{m\in (M\otimes_{F[t]}\FI^s\langle t \rangle)/(t^n)~|~m=\tau_M(\tau^*m)\right\}\otimes \FI^s\longrightarrow  (M\otimes_{F[t]}\FI^s\langle t \rangle)/(t^n) \nonumber
\end{equation}
is an isomorphism. In particular, the inclusion
\[
\left\{m\in (M\otimes_{F[t]}\FI^s\langle t \rangle)/(t^n)\mid m=\tau_M(\tau^*m)\right\} \subseteq \left\{m\in (M\otimes_{F[t]}\bC_{\infty}\langle t \rangle)/(t^n)\mid m=\tau_M(\tau^*m)\right\}
\]
is an equality. This shows that $\Betti{\underline{M}}$ is both a submodule of $M\otimes_{F[t]}\FI^s \langle t \rangle$ and of $M\otimes_{F[t]}L'\langle t \rangle$. Because $M$ is free over $F[t]$, it follows that $\Betti{\underline{M}}\subset M\otimes_{F[t]}L\langle t \rangle$ where $L=L'\cap \FI^s$ is a finite separable extension of $\FI$ in $\CI$. As $(A\otimes F)\otimes_{F[t]} L\langle t\rangle$ is isomorphic to $L\langle A \rangle$, we deduce that $\Betti{\underline{M}}\subset M\otimes_{A\otimes F}L\langle A \rangle$.
\end{proof}

As a consequence of Proposition \ref{prop:finite-extension-betti}, by the faithful flatness of the inclusion $L\langle A \rangle \to \CI\langle A \rangle$ (\cite[AC I\S 3.5 Proposition 9]{bourbaki}), we have:
\begin{Proposition}
Let $\underline{M}$ be a rigid analytically trivial $A$-motive over $F$. Let $L$ be chosen as in Proposition \ref{prop:finite-extension-betti}. The multiplication map 
\[
\Betti{\underline{M}}\otimes_A L\langle A \rangle \longrightarrow M\otimes_{A\otimes F}L\langle A \rangle
\]
is an isomorphism of $L\langle A \rangle$-modules.
\end{Proposition}

The first part of the next result is inspired by \cite[Proposition~6.1]{hartl-boeckle}. We have adapted its proof to allow the smaller field $\FI^s$ instead of $\CI$. This is needed in order to compute the $A$-module $\operatorname{H}^1(G_{\infty},\Betti{\underline{M}})$ of continuous Galois cohomology.
\begin{Theorem}\label{thm-analytically}
Let $\underline{M}$ be a rigid analytically trivial $A$-motive over $F$. There is an exact sequence of $A[G_\infty]$-modules:
\begin{equation}\label{exact-sequence-over-CI}
0\longrightarrow \Betti{\underline{M}}\longrightarrow M\otimes_{A\otimes F}\bC_{\infty}\langle A \rangle \xrightarrow{\id-\tau_M} M\otimes_{A\otimes F}\bC_{\infty}\langle A \rangle\longrightarrow 0.
\end{equation}
Furthermore, it induces a long exact sequence of $A$-modules 
\begin{equation}\label{exact-sequence-over-KI}
0\longrightarrow \Betti{\underline{M}}^{+}\longrightarrow M\otimes_{A\otimes F}\FI\langle A \rangle \xrightarrow{\id-\tau_M} M\otimes_{A\otimes F}\FI\langle A \rangle \longrightarrow \operatorname{H}^1(G_\infty,\Betti{\underline{M}})\longrightarrow 0.
\end{equation}
\end{Theorem}
\begin{Remark}
The fact that \eqref{exact-sequence-over-CI} implies \eqref{exact-sequence-over-KI} is not immediate. We have to descend from the completion of the perfection of $\FI$---which, by the Ax--Sen--Tate theorem, corresponds to the fixed subfield of $\CI$ under $G_\infty$---to the much smaller field $F_\infty$.
\end{Remark}

\begin{proof}[Proof of Theorem \ref{thm-analytically}.]
Let $\bF[t]\to A$ be a non-constant morphism of rings. We have $\bC_{\infty}\langle A\rangle=A\otimes_{\bF[t]}\bC_{\infty}\langle t\rangle$, where $\bC_{\infty}\langle t\rangle$ is the Tate algebra over $\bC_{\infty}$ in the variable $t$.

The exactness of \eqref{exact-sequence-over-CI} follows from \cite[Proposition~6.1]{hartl-boeckle}. We use the same argument as in \loccit to show that the sequence
\begin{equation}\label{exact-sequence-over-KIs}
0\longrightarrow \Betti{\underline{M}}\longrightarrow M\otimes_{F[t]}\FI^s\langle t \rangle \xrightarrow{\id-\tau_M} M\otimes_{F[t]}\FI^s\langle t \rangle\longrightarrow 0 
\end{equation}
is exact, where the first inclusion is well defined by Proposition \ref{prop:finite-extension-betti}. It suffices to show the surjectivity of $\id-\tau_M$ on $M\otimes_{F[t]}\FI^s\langle t \rangle$. Let $f\in M\otimes_{F[t]}\FI^s\langle t \rangle$. Since $\underline{M}$ is rigid analytically trivial, we may assume without loss of generality that $f=c\cdot \omega$ for $c=\sum_{n\geq 0}{c_nt^n}\in \FI^s\langle t \rangle$ and $\omega\in \Betti{\underline{M}}$. By definition of $\FI^s\langle t \rangle$, there exists a finite separable extension $L$ of $\FI$ such that $c_n\in L$ for all $n\geq 0$. As $(c_n)_n$ is a null sequence, there exists $n_0\geq 0$ large enough such that for all $n>n_0$ the series $b_n:=c_n+c_n^q+c_n^{q^2}+\ldots$ converges in $L$ and satisfies $b_n-b_n^q=c_n$. For $n\leq n_0$, let $b_n\in \FI^s$ be a solution of $x-x^q=c_n$. That $(c_n)_n$ is a null sequence also implies that $(b_n)_n$ is a null sequence. In particular, the series 
\[
b:=\sum_{n=0}^{\infty}{b_nt^n}
\]
belongs to $L(b_0,\ldots,b_{n_0})\langle t \rangle\subset \FI^s\langle t\rangle$. Hence, the element $g:=b \omega$ belongs to $M\otimes_{F[t]}\FI^s\langle t \rangle$ and satisfies $(\id-\tau_M)(g)=f$. Surjectivity follows. 

We turn to the second part of the statement. By Proposition \ref{prop:finite-extension-betti}, $G_\infty$ acts continuously on \eqref{exact-sequence-over-KIs}, where each module is equipped with the discrete topology. Taking invariants yields a long exact sequence of $A$-modules:
\begin{equation}
0\longrightarrow \Betti{\underline{M}}^{+}\longrightarrow  M\otimes_{A\otimes F}\FI\langle A \rangle \xrightarrow{\id-\tau_M} M\otimes_{A\otimes F}\FI\langle A \rangle\longrightarrow \operatorname{H}^1(G_\infty,\Betti{\underline{M}})\longrightarrow \cdots \nonumber
\end{equation}
We prove that the $F[t]$-module
\begin{equation}
\operatorname{H}^1(G_\infty,M\otimes_{F[t]} \FI^s\langle t \rangle),\nonumber
\end{equation}
which appears in the dotted part, is zero. Note that it is isomorphic to $M\otimes_{F[t]}\operatorname{H}^1(G_\infty,\FI^s\langle t \rangle)$, hence it suffices to show that $\operatorname{H}^1(G_\infty,\FI^s\langle t \rangle)$ vanishes. Because $\FI^s\langle t \rangle$ is endowed with the discrete topology, a continuous cocycle $G_{\infty}\to \FI^s\langle t \rangle$ descends to a cocycle $c:H\to L\langle t\rangle$ for some finite Galois extension $L$ of $\FI$, where $H$ is the finite group $\Gal(L|\FI)$. 

Because $L$ is separable over $\FI$, its trace form is non-degenerate; that is, there exists $\alpha\in L$ with $\operatorname{Tr}_{L/\FI}(\alpha)=1$. Then set 
\[
f=-\sum_{\sigma\in H}c(\sigma)\sigma(\alpha)\in L\langle t\rangle.
\]
A computation shows that $c(\rho)=\rho(f)-f$ for all $\rho\in H$. It follows that $c$ is a coboundary, and hence that $\operatorname{H}^1(G_\infty,\FI^s\langle t\rangle)=0$. This concludes the proof.
\end{proof}

We are now ready to introduce the category of rigid analytically trivial $A$-motives over $F$, as mentioned in the introduction.
\begin{Definition}\label{def:MM-F-vrig}
We let $\AMot_F^{\operatorname{rat}}$ be the full subcategory of $\AMot_F$ whose objects are rigid analytically trivial.
\end{Definition}

Recall that a sequence of objects in $\AMot_F$ is called \emph{exact} if the underlying sequence of modules is (see \cite[Definition 2.14]{gazda}), and that the category $\AMot_F$ endowed with the class of exact sequences forms an exact category (\cite[Proposition 2.15]{gazda}). The next proposition, which ensures that extension modules in the category $\AMot_F^{\operatorname{rat}}$ are well defined, is borrowed from \cite[Lemma 2.3.25]{hartl-juschka}.
\begin{Proposition}\label{prop:Mrat-iif-M'rat-M''rat}
Let $0\to \underline{M}'\to \underline{M}\to \underline{M}''\to 0$ be an exact sequence in $\AMot_F$. Then $\underline{M}$ is rigid analytically trivial if and only if $\underline{M}'$ and $\underline{M}''$ are. In particular, the category $\AMot_F^{\operatorname{rat}}$ is exact.
\end{Proposition}

We finally record that the Betti realization functor, with source $\AMot_F^{\operatorname{rat}}$, is exact (this is not true on the full category $\AMot_F$).
\begin{Corollary}\label{cor:betti-exact}
The functor $\underline{M}\mapsto \Betti{\underline{M}}$ from $\AMot_F^{\operatorname{rat}}$ to the category $\operatorname{Rep}_{A}(G_\infty)$ of continuous $A$-linear representations of $G_\infty$ is exact.
\end{Corollary}
\begin{proof}
This follows from Theorem \ref{thm-analytically} together with the snake lemma.
\end{proof}

\subsubsection*{Extensions with analytic reduction}
In view of Corollary \ref{cor:betti-exact}, the following definition makes sense:
\begin{Definition}[Analytic reduction at $\infty$]\label{def:analytic-reduction-mot}
Let $S:0\to \underline{M}'\to \underline{M}\to \underline{M}''\to 0$ be an exact sequence in $\AMot_F^{\operatorname{rat}}$. We say that $S$ \emph{has analytic reduction at $\infty$} if the sequence $S_B$ splits in $\operatorname{Rep}_{A}(G_\infty)$.
\end{Definition}

\begin{Remark}
The above notion is crucial for obtaining finitely generated extension modules from the category $\AMot_F^{\operatorname{rat}}$ (see Section \ref{sec:reg-fin-thm}). It is the counterpart of Definition \ref{def:analytic-reduction-HP} for Hodge-Pink structures.
\end{Remark}

Given two rigid analytically trivial $A$-motives $\underline{N}, \underline{M}$, the exactness of the Betti realization functor induces a morphism of $A$-modules:
\begin{equation}\label{eq:rB(N,M)}
r_B(\underline{N},\underline{M}):\Ext^1_{\AMot_F^{\operatorname{rat}}}(\underline{N},\underline{M}) \longrightarrow \Ext^1_{\operatorname{Rep}_{A}(G_\infty)}(\Betti{\underline{N}},\Betti{\underline{M}})
\end{equation}
which is the counterpart of \eqref{eq:forgetful-on-ext} in the context of Hodge-Pink structures. By definition, the kernel of $r_B(\underline{N},\underline{M})$ consists of extensions having analytic reduction at $\infty$. We set:
\begin{equation}\label{eq:ext-an-red-mot}
\Ext^1_{F,\infty}(\underline{N},\underline{M}):=\ker r_B(\underline{N},\underline{M}).
\end{equation}

\subsection{Analytic continuation}\label{subsection:analytic continuation}
To associate a Hodge-Pink structure to a rigid analytically trivial $A$-motive $\underline{M}$, it is crucial to understand the behaviour of elements of $\Betti{\underline{M}}$---which can be viewed as \emph{functions} on the affinoid subdomain $\Spm~\CI\langle A \rangle$ \emph{with values in} $M\otimes_{F}\CI$---near the point $\fj$. However, the latter is not a point of this domain, as one deduces from Lemma \ref{lem:j-invertible-in-Tate-algebra}. Hence it is necessary to extend elements of $\Betti{\underline{M}}$ to a larger domain. In this subsection, we show that elements of $\Betti{\underline{M}}$ can indeed be uniquely continued \emph{meromorphically} to the whole rigid analytification of the affine curve $\Spec A\otimes \CI$, with poles supported only at $\fj$ and its Frobenius twists $\fj_1$, $\fj_2$, \dots. In the case $\deg(\infty)=1$, this is treated in \cite[\S 2.3.4]{hartl-juschka}; for special functions of Anderson $A$-modules, this is the subject of \cite{gazda-maurischat2}. 

\subsubsection*{Prerequisites of rigid analytic geometry}
Let $L$ be a non-archimedean local field equipped with a continuous map $\iota:F\hookrightarrow L$ with respect to the topology induced by $v_{\infty}$. We denote by $\mathcal{O}_L$ its valuation ring and $k_L$ its residue field. We define $C_L^{\mathrm{an}}$ as the rigid analytification of $C_L:=C\times \Spec L$. We also denote by $\mathfrak{U}=\mathfrak{U}_L$ the analytification of $U_L:=(C\setminus \{\infty\})\times \Spec L$. 

\begin{Construction}[Analytic and meromorphic functions on $\mathfrak{U}$]
We recall some material from \cite[\S 3]{gazda-maurischat2} to determine $\mathcal{O}(\mathfrak{U})$, the algebra of analytic functions on $\mathfrak{U}$. Let $|\cdot|$ be the norm on $L$. We generalize Gauss norms on $A\otimes L$ as follows: let $\rho>1$. For $f\in A\otimes L$, we set
\[
\|f\|_\rho:=\inf \left(\max_i\{|\ell_i|\rho^{\deg a_i}\}\right),
\]
where the infimum is taken over all representations of $f$ as a finite sum $\sum_{i}{a_i\otimes \ell_i}$ in $A\otimes L$. The map $\|\cdot \|_\rho$ is indeed a Gauss norm \cite[Proposition 3.2]{gazda-maurischat2}, and we denote by $L\langle A\rangle_\rho$ the affinoid algebra over $L$ obtained by completing $A\otimes L$ with respect to it. Observe that $L\langle A \rangle$ as in Definition \ref{def:aff-Tate-algebra} coincides with $L\langle A \rangle_1$ (see \cite[Proposition 2.2]{gazda-maurischat}). If $\rho<\rho'$, there is a canonical inclusion $L\langle A\rangle_{\rho'}\hookrightarrow L\langle A\rangle_{\rho}$. We have the following equivalent description of $\mathfrak{U}$:
\[
\mathfrak{U}:=\bigcup_{\rho >1} \mathrm{Spm} L\langle A \rangle_{\rho}
\]
where $\rho$ runs over $\mathbb{R}_{>1}$. The sheaf property for this admissible covering gives the following formula to compute its ring of analytic functions:
\begin{equation}\label{eq:analytic-functions-on-U}
\mathcal{O}(\mathfrak{U}) = \varprojlim_{\rho >1} L\langle A \rangle_{\rho}.
\end{equation}
One may derive a similar formula for the algebra of meromorphic functions on $\mathfrak{U}$:
\begin{equation}\label{eq:meromorphic-functions-on-U}
\mathcal{M}(\mathfrak{U})=\varprojlim_{\rho >1} Q(L\langle A\rangle_\rho).
\end{equation}
The reader unfamiliar with rigid analytic geometry may very well take the formulas \eqref{eq:analytic-functions-on-U} and \eqref{eq:meromorphic-functions-on-U} as definitions for $\mathcal{O}(\mathfrak{U})$ and $\mathcal{M}(\mathfrak{U})$. We note that by standard non-archimedean analysis, $\mathcal{M}(\mathfrak{U})$ coincides with $Q(\mathcal{O}(\mathfrak{U}))$. Since $\|\tau(f)\|_{\rho}=\|f\|_{\rho^{1/q}}^q$ for $f\in A\otimes L$, $\tau$ induces compatible morphisms $L\langle A\rangle_{\rho}\to L\langle A \rangle_{\rho^{1/q}}$. We still denote by $\tau$ the resulting endomorphism of $\mathcal{O}(\mathfrak{U})$ or $\mathcal{M}(\mathfrak{U})$.
\end{Construction}

The closed point $\mathfrak{j}$---corresponding to the eponym ideal---is the $L$-valued point of $C$ given by the composition $\Spec L \to \Spec K\xrightarrow{\eta} C$ where $\eta$ is the generic point. Its iterates by the Frobenius $\mathfrak{j}_n:\Spec L \xrightarrow{\mathrm{Frob}^n_L} \Spec L \xrightarrow{\mathfrak{j}} C$ define points of $\mathfrak{U}$ and we denote by $J$ the set $\{\fj,\fj_1,\fj_2,\ldots\}$. For $I\subset J$ a subset, we define $\mathcal{M}_I(\mathfrak{U})$ as the algebra of meromorphic functions whose poles are supported at $I$:
\[
\mathcal{M}_I(\mathfrak{U}) = \varprojlim_{\rho >1} L\langle A \rangle_{\rho}[\mathfrak{i}^{-1}~|~\mathfrak{i}\in I]. 
\]
In the above formula, we use the same notation for the point $\mathfrak{i}\in I$ and the corresponding invertible ideal of $A\otimes L$. In particular, for fixed $\rho$, only finitely many ideals of $I$ have to be inverted. The following criterion follows from \cite[Lemma~2.4]{gazda-maurischat2}.

\begin{Lemma}\label{lem:global-section-polar-at-J}
For $u\in A\setminus \bF$, the following product converges in $\mathcal{O}(\mathfrak{U})$:
\[
\Pi_u:=\prod_{j=0}^{\infty}{\left( 1-\frac{u\otimes 1}{1\otimes u^{q^j}}\right)}.
\]
In addition, let $f\in \mathcal{M}(\mathfrak{U})$ and assume the existence of an integer $n\geq 0$ such that, for all $u\in A\setminus \bF$, $\Pi_u ^n\cdot f\in \mathcal{O}(\mathfrak{U})$. Then $f\in \mathcal{M}_J(\mathfrak{U})$. 
\end{Lemma} 

\begin{Example}\label{example:algebra in the t-case}
In the $t$-setting (\cf Notation \ref{not:t-setting}), $\fj_i$ corresponds to the principal ideal $(t-\theta^{q^i})$ of $L[t]$. We have 
\begin{equation}
L\langle A\rangle=L\langle t\rangle=\left\{\sum_{k=0}^{\infty}{a_k t^k}~\bigg\vert~a_k\in L;~a_k\xrightarrow[k\to\infty]{}0\right\}, \nonumber
\end{equation}
\begin{equation}\label{eq:entire-series}
\mathcal{O}(\mathfrak{U})=L\llangle t\rrangle=\left\{\sum_{k=0}^{\infty}{a_k t^k}~\bigg\vert~a_k\in L;~\forall \rho>1:~|a_k|\rho^k\xrightarrow[k\to\infty]{}0\right\}. 
\end{equation}
The algebra $L\langle A\rangle$ corresponds to series converging on the \emph{closed} unit disc, whereas $\mathcal{O}(\mathfrak{U})$ consists of entire series. The algebra $\mathcal{M}_J(\mathfrak{U})$ is identified with the subalgebra
\[
\bigcap_{\rho >1}\left(\left\{\sum_{k=0}^{\infty}{a_k t^k}~\bigg\vert~a_k\in L;~a_k\rho^k\xrightarrow[k\to\infty]{}0\right\}\left[\frac{1}{(t-\theta^{q^i})}~\Bigg |~i\geq 0\right]\right).
\]
\end{Example}

\begin{Construction}[Laurent expansion map]\label{cons:Laurent-series-expansion-map}
Let $u\in A\setminus \bF$ denote a separating element as in \cite[\S 2.2]{gazda-maurischat} and choose $\rho>|\iota(u)|$. Then it follows from Weierstrass division that the canonical map
\[
A\otimes L/\mathfrak{j}^n \longrightarrow L\langle A \rangle_{\rho}/\mathfrak{j}^n
\]
is an isomorphism for all $n\geq 0$. Passing to the limit as $n$ grows, we obtain an injective map 
\[
\mathcal{O}(\mathfrak{U})\longrightarrow L\langle A \rangle_{\rho} \longrightarrow L[\![\fj]\!]
\]
that we call \emph{the Laurent expansion map at $\fj$.} Passing to quotient fields immediately extends it to an injective map
\[
\mathcal{M}(\mathfrak{U})\longrightarrow L(\!(\fj)\!).
\]
\end{Construction}

\subsubsection*{Analytic continuation of elements in the Betti realization}
Given $L\subset \CI$ a finite separable extension of $\FI$, we have a canonical injective map $\mathcal{M}_J(\mathfrak{U}_L)\to L\langle A \rangle$. We are now in a position to prove the main result of this subsection (compare with \cite[Proposition 2.3.30]{hartl-juschka}).
\begin{Theorem}\label{thm:analytic-continuation}
Let $\underline{M}$ be a rigid analytically trivial $A$-motive over $F$. There exists a finite separable extension $L$ of $\FI$ such that the canonical inclusion
\[
\left(M\otimes_{A\otimes F}\mathcal{M}_J(\mathfrak{U}_L)\right)^{\tau_M=1}\hookrightarrow \Betti{\underline{M}}
\]
is an isomorphism.
\end{Theorem}

\begin{proof}
Let $n\geq 0$ be such that $\fj^n \tau_M(\tau^*M)\subset M$. By Proposition \ref{prop:finite-extension-betti}, there exists a finite separable extension $L$ of $\FI$ such that
\[
\underline{M}_B\subset M\otimes_{A\otimes F}L\langle A \rangle.
\]
Let $u\in A$ be non-constant. By Lemma \ref{lem:global-section-polar-at-J}, it suffices to show that $\Pi_u^n \cdot \Betti{\underline{M}}\subset  M\otimes_{A\otimes F}\mathcal{O}(\mathfrak{U})$. Choose a $(q-1)$st root $\eta_u$ of $-u$ in $\FI^s$. Then
\begin{equation}\label{eq:omega}
\omega_u:=\eta_u\Pi_u^{-1} \in Q(\mathcal{O}(\mathfrak{U}_{L'})), \quad \text{where}\quad L':=L(\eta_u)
\end{equation}
is such that $\omega_u^{-n}$ generates the Betti realization of the $\bF[u]$-motive $\underline{A}(-n)$. Set $\underline{N}:=\underline{A}(-n)$. Then $\underline{M}\otimes \underline{N}$---where the tensor product is taken in the category of $\bF[u]$-motives---is effective, and from \cite[Proposition~3.4]{hartl-boeckle} we deduce that $\Betti{(\underline{M}\otimes \underline{N})}\cong  \Betti{\underline{M}}\otimes_{\bF[u]}\Betti{\underline{N}}$ is contained in $M\otimes_{F[u]}L'\llangle u\rrangle$. Since\footnote{The map $A\otimes_{\bF[u]}L\langle u \rangle_{\rho} \to L\langle A \rangle_{\rho}$ is an isomorphism (it is a surjective map of finite free $L\langle u \rangle_{\rho}$-modules of rank $\deg(u)$), and since $A$ is a finite free $\bF[u]$-module, tensoring along $\bF[u]\to A$ preserves limits.} $A\otimes_{\bF[u]}L'\llangle u\rrangle \stackrel{\sim}{\to} \mathcal{O}(\mathfrak{U}_{L'})$, we obtain $\omega_u^{-n} \cdot \Betti{\underline{M}}\subset M\otimes_{A\otimes F}\mathcal{O}(\mathfrak{U}_{L'})$, and the desired statement follows by Galois descent.
\end{proof}

\begin{Corollary}\label{cor:mult-continuation}
Under the notation and assumptions of Theorem \ref{thm:analytic-continuation}, the multiplication map
\[
\Betti{\underline{M}}\otimes_A \mathcal{M}_{J}(\mathfrak{U}) \longrightarrow M\otimes_{A\otimes F}\mathcal{M}_{J}(\mathfrak{U})
\]
is bijective.
\end{Corollary}
\begin{proof}
Choose $u\in A$ non-constant. The $\bF[u]$-linear determinant of this multiplication map is, up to a scalar multiple in $L$, a power of $\omega_u$. Hence it is invertible in $\mathcal{M}_{J}(\mathfrak{U})$. 
\end{proof}

\subsection{The associated Hodge-Pink structure}\label{subsec:asso-MHPS}
Let $\underline{M}$ be a rigid analytically trivial $A$-motive over $F$ (Definition \ref{def:rigid analytically trivial}). In this subsection, we associate a Hodge-Pink structure to $\underline{M}$, following unpublished work of Pink. 

Let $\Betti{\underline{M}}$ be the Betti realization of $\underline{M}$ (Definition \ref{def:betti-realization}). By Corollary \ref{cor:mult-continuation}, there exists a finite separable extension $L$ of $\FI$ such that the multiplication map
\begin{equation}\label{eq:mult-mero}
\Betti{\underline{M}}\otimes_A \mathcal{M}_{J}(\mathfrak{U})\longrightarrow M\otimes_{A\otimes F}\mathcal{M}_{J}(\mathfrak{U}),\quad \omega\otimes c\longmapsto \omega \cdot c
\end{equation}
is bijective. Base-changing along the Laurent expansion map of Construction \ref{cons:Laurent-series-expansion-map}, we obtain an isomorphism of $\FI^s(\!(\fj)\!)$-modules:
\begin{equation}\label{eq:gammaMv}
\Betti{\underline{M}}\otimes_{A,\nu} \FI^s(\!(\fj)\!)\stackrel{\sim}{\longrightarrow} M\otimes_{A\otimes F}\FI^s(\!(\fj)\!),
\end{equation}
where $\nu$ denotes the morphism $A\to \FI^s[\![\fj]\!]$, $a\mapsto a\otimes 1$, introduced earlier in Section \ref{sec:ext-of-MHPS} in the context of Hodge-Pink structures.
\begin{Definition}\label{def:gamma}
We denote by $\gamma_{\underline{M}}$ the isomorphism \eqref{eq:gammaMv}.
\end{Definition}
A trivial yet important remark is the following.
\begin{Lemma}
The morphism $\gamma_{\underline{M}}$ is $G_\infty$-equivariant, where $\sigma\in G_\infty$ acts on the left-hand side of \eqref{eq:gammaMv} via $\sigma\otimes \sigma$ and on the right via $\id_M\otimes \sigma$. 
\end{Lemma}

In the next definition, which is due to Pink, we attach a Hodge-Pink structure to $\underline{M}$ (see also \cite[Definition~2.3.32]{hartl-juschka}).
\begin{Definition}[Hodge-Pink structure associated to $\underline{M}$]\label{def:hodge-pink-real-functor}
We let $\sH(\underline{M})$ be the Hodge-Pink structure 
\begin{itemize}
\item whose underlying $\KI$-vector space is $\Betti{\underline{M}}\otimes_A \KI$,
\item whose Hodge-Pink lattice is $\fq_{\underline{M}}=\gamma_{\underline{M}}^{-1}(M\otimes_{A\otimes F}\FI^s[\![\fj]\!])$.
\end{itemize} 
The tautological lattice of $\sH(\underline{M})$ is $\fp_{\underline{M}}=\Betti{\underline{M}}\otimes_{A,\nu} \FI^s[\![\fj]\!]$. The action of $G_\infty$ on $\Betti{\underline{M}}$ is continuous (Proposition \ref{prop:finite-extension-betti}) and defines an infinite Frobenius $\phi_{\underline{M}}$ for $\sH(\underline{M})$. We denote by $\sH^+(\underline{M})$ the pair $(\sH(\underline{M}),\phi_{\underline{M}})$.
\end{Definition}

\begin{Proposition}\label{prop:exactness-HodgePink}
The assignment $\underline{M}\mapsto \sH(\underline{M})$ (resp. $\sH^+(\underline{M})$) defines an exact functor $\sH:\AMot_F^{\operatorname{rat}}\to \HPk$ (resp. $\sH^+:\AMot_F^{\operatorname{rat}}\to \HPk^+$) between exact categories.
\end{Proposition}
\begin{proof}
Let $S:0\to \underline{M}'\to \underline{M}\to \underline{M}''\to 0$ be an exact sequence in $\AMot_F^{\text{rat}}$. By Corollary \ref{cor:betti-exact}, the sequence $\Betti{S}$ is exact, hence $\sH(S)$ and $\sH^+(S)$ are exact. To show that they are strict, observe that the sequence
\[
0\longrightarrow \fq_{\underline{M}'} \longrightarrow \fq_{\underline{M}} \longrightarrow \fq_{\underline{M}''} \longrightarrow 0
\]
is isomorphic to $0\to M'\otimes_{A\otimes F}\FI^s[\![\fj]\!]\to M\otimes_{A\otimes F}\FI^s[\![\fj]\!]\to M''\otimes_{A\otimes F}\FI^s[\![\fj]\!]\to 0$ via $\gamma_S$, and the latter is exact by flatness of $A\otimes F\to \FI^s[\![\fj]\!]$. 
\end{proof}

We recall the notion of \emph{regulated extensions} in $\AMot_F^{\operatorname{rat}}$, as introduced in \cite[Definition~5.1]{gazda}:
\begin{Definition}[Regulated extensions]\label{def:regulated}
Let $S:0\to \underline{M}'\to \underline{M}\to \underline{M}''\to 0$ be an exact sequence in $\AMot_F^{\text{rat}}$. We say that $S$ is \emph{regulated} if $\sH(S)$ is Hodge additive (Definition \ref{def:hodge-additive}).
\end{Definition}

We conclude this section by giving a description of the extension modules of Hodge-Pink structures arising from $A$-motives. It consists mainly of a reformulation of Propositions \ref{prop:ext-+} and \ref{prop:ext-+-ha} in the case $\underline{H}^+=\sH^+(\underline{M})$ for a rigid analytically trivial $A$-motive $\underline{M}$ over $F$.

\begin{Theorem}\label{thm:extensions-of-MHPS-arising-from-motives}
Let $\underline{M}$ be a rigid analytically trivial $A$-motive over $F$ whose weights are all non-positive. Let $\underline{H}^+=\sH^+(\underline{M})$. In the notation of Proposition \ref{prop:ext-+}, we have an exact sequence 
\[
0\longrightarrow \frac{M\otimes_{A\otimes F}\FI(\!(\fj)\!)}{(\Betti{\underline{M}}^+)_{\KI}+M\otimes_{A\otimes F}\FI[\![\fj]\!]}\xrightarrow{\varphi^+} \Ext^1_{\HPk^+}(\mathbbm{1}^+,\underline{H}^+)\xrightarrow{d_{\underline{H}^+}} \operatorname{H}^1(G_\infty,(\Betti{\underline{M}})_{\KI})
\]
where $(\Betti{\underline{M}})_{\KI}:=\Betti{\underline{M}}\otimes_A \KI$. The Hodge-additive version of this exact sequence holds:
\[
0\longrightarrow \frac{(M+\tau_M(\tau^*M))\otimes_{A\otimes F}\FI[\![\fj]\!]}{(\Betti{\underline{M}}^+)_{\KI}+M\otimes_{A\otimes F}\FI[\![\fj]\!]}\xrightarrow{\varphi^+} \Ext^{1,\operatorname{ha}}_{\HPk^+}(\mathbbm{1}^+,\underline{H}^+)\xrightarrow{d_{\underline{H}^+}} \operatorname{H}^1(G_\infty,(\Betti{\underline{M}})_{\KI}).
\]
\end{Theorem}

The theorem follows from the next two lemmas, the first of which specifies the form of $\fp_{\underline{M}}$ when viewed as a submodule of $M\otimes_{A\otimes F}\FI^s(\!(\fj)\!)$.
\begin{Lemma}\label{lem:tautological-lattice-tau}
We have $\gamma_{\underline{M}} (\fp_{\underline{M}})=\tau_M(\tau^*M)\otimes_{A\otimes F}\FI^s[\![\fj]\!]$. 
\end{Lemma}
\begin{proof}
Let $L$ be a finite separable extension of $\FI$ as in Theorem \ref{thm:analytic-continuation}. By Lemma \ref{lem:twist-polar-part}, pulling back \eqref{eq:mult-mero} by $\tau$ produces an isomorphism
\begin{equation}
\Betti{\underline{M}}\otimes_A \mathcal{M}_{J^{(1)}}(\mathfrak{U}) \stackrel{\sim}{\longrightarrow} (\tau^*M)\otimes_{A\otimes F}\mathcal{M}_{J^{(1)}}(\mathfrak{U}) \nonumber
\end{equation}
which, upon localizing at $\fj$, yields
\begin{equation}
\delta_{\underline{M}}:\Betti{\underline{M}}\otimes_{A}\FI^s[\![\fj]\!]\stackrel{\sim}{\longrightarrow} (\tau^*M)\otimes_{A\otimes F}\FI^s[\![\fj]\!]. \nonumber
\end{equation}
The latter fits into a commutative diagram
\begin{equation}
\begin{tikzcd}%[column sep=8em,row sep=4em]
\Betti{\underline{M}}\otimes_A \FI^s(\!(\fj)\!)\arrow[r,"\delta_{\underline{M}}"]\arrow[dr,"\gamma_{\underline{M}}"'] & (\tau^*M)\otimes_{A\otimes F}\FI^s(\!(\fj)\!)\arrow[d,"\tau_M\otimes \id_{\FI^s(\!(\fj)\!)}"] \\
 & M\otimes_{A\otimes F}\FI^s(\!(\fj)\!)
\end{tikzcd} \nonumber
\end{equation}
which already appears in \cite[Proposition~2.3.30]{hartl-juschka} under different notation. The equality $\gamma_{\underline{M}}(\fp_{\underline{M}})=\tau_M(\tau^*M)\otimes_{A\otimes F}\FI^s[\![\fj]\!]$ follows from the commutativity of the above diagram together with the fact that both $\gamma_{\underline{M}}$ and $\delta_{\underline{M}}$ are isomorphisms.
\end{proof}

To apply Proposition \ref{prop:ext-+-ha} and conclude the proof of Theorem \ref{thm:extensions-of-MHPS-arising-from-motives}, we need a vanishing result for Galois cohomology, supplied by the next lemma.

\begin{Lemma}\label{lem:equivariant-lattices-hodge}
Let $\fl$ be a $G_{\infty}$-equivariant $\FI^s[\![\fj]\!]$-lattice in $\Betti{\underline{M}}\otimes_A \FI^s(\!(\fj)\!)$. Then $G_{\infty}$ acts continuously on $\fl$ for the $\fj$-adic topology and $\operatorname{H}^1(G_\infty,\fl)=(0)$ as continuous Galois cohomology.
\end{Lemma}
\begin{proof}
The action of $G_{\infty}$ on $\FI^s(\!(\fj)\!)$ is continuous for the $\fj$-adic topology. Because $G_{\infty}$ also acts continuously on the discrete module $\underline{M}_B$, the induced semilinear action on
\[
\underline{M}_B\otimes_A \FI^s(\!(\fj)\!)
\]
is also continuous for the $\fj$-adic topology. This topology is equivalently induced from $\fl$ (making $(\fj^n \fl)_{n\geq 0}$ a neighborhood basis at zero); since $\fl$ is $G_{\infty}$ stable, we deduce that the action restricted to $\fl$ is also continuous. This proves the first assertion.

Now, let $W$ be any finite-dimensional $\FI^s$-vector space equipped with a continuous semilinear action of $G_{\infty}$ where $W$ is discrete. The additive version of the Hilbert $90$ theorem \cite[x.\S1,~Proposition~1]{serre} implies $H^1(G_{\infty},W)=(0)$.

Let $c:G_{\infty}\to \fl$ be a continuous cocycle and let $c_n$ denote its reduction modulo $\fj^n$ for $n\geq 0$. By the above fact for $W=\fl/\fj^n \fl$, $c_n$ is trivial and there exists $b_n\in \fl/\fj^n \fl$ such that $c_n(\sigma)=\,^{\sigma}b_n-b_n$ for all $\sigma\in G_{\infty}$. 

Applying this to $W=\fj^{n-1}\fl/\fj^n \fl$, we also deduce that $H^1(G_{\infty},\fj^{n-1}\fl/\fj^n \fl)$ vanishes, and in particular that the reduction map $(\fl/\fj^{n}\fl)^{G_{\infty}}\to (\fl/\fj^{n-1}\fl)^{G_{\infty}}$ is surjective. This in turn implies that the sequence $(b_n)_n$ may be chosen compatibly; \ie in $\lim_n \fl/\fj^n \fl=\fl$. This shows the existence of $b\in \fl$ such that $c(\sigma)=\,^{\sigma}b-b$, hence that $c$ is trivial.
\end{proof}

\section{Regulators and finiteness theorems}
Let $F$ be a finite extension of $K$ which is totally ramified at $\infty$ (see Reduction \ref{red:F-totally-ramified}). Let $\underline{M}$ be a rigid analytically trivial $A$-motive over $F$, and let $\underline{H}^+$ denote the Hodge-Pink structure with infinite Frobenii associated to $\underline{M}$ in Definition \ref{def:hodge-pink-real-functor}. In this section, we associate to $\underline{M}$ a \emph{Hodge-Pink regulator map}
\[
\Reg(\underline{M}):\Ext^{1,\operatorname{reg}}_{\cO_F,\infty}(\mathbbm{1},\underline{M})\longrightarrow \Ext^{1,\operatorname{ha}}_{\infty}(\mathbbm{1}^+,\underline{H}^+)
\]
resulting from the exactness of the Hodge-Pink realization functor; see Subsection \ref{subsec:regulator}. The source of the regulator is the $A$-module of integral, regulated extensions having analytic reduction at $\infty$. We define this module in Subsection \ref{subsec:integral-regulated-analytic-extensions} together with \emph{the class module of $\underline{M}$}. The target of the regulator is the space of Hodge-additive extensions with analytic reduction at $\infty$ in the category of Hodge-Pink structures with infinite Frobenii, studied in Subsection \ref{subsec:asso-MHPS}.

\subsection{Extensions of $A$-motives}\label{subsec:integral-regulated-analytic-extensions}
Let $\underline{M}=(M,\tau_M)$ be an $A$-motive over $F$. Let $\cO_F$ denote the ring of integers of $F$; \ie the integral closure of $A$ in $F$. Recall from \cite[\S 4]{gazda} that there exists a unique $A\otimes \cO_F$-submodule $M_{\cO_F}\subset M$ with the following properties:
\begin{enumerate}
    \item\label{item:maxmod-1} $M_{\cO_F}$ is finitely generated and generates $M$ over $A\otimes F$,
    \item\label{item:maxmod-2} $M_{\cO_F}$ is stable under $\tau_M$; namely, $\tau_M(\tau^*M_{\cO_F})\subset M_{\cO_F}[\fj^{-1}]$,
    \item\label{item:maxmod-3} $M_{\cO_F}$ is maximal; \ie it is not strictly contained in another submodule satisfying \ref{item:maxmod-1} and \ref{item:maxmod-2}.
\end{enumerate}
The module $M_{\cO_F}$ is called \emph{the maximal integral model} of $\underline{M}$ over $\cO_F$. One can show that $M_{\cO_F}$ is a projective $A\otimes \cO_F$-module \cite[Theorem 4.32]{gazda}. We also introduce the following submodule of $M[\fj^{-1}]$: 
\begin{equation}\label{def:NA}
N_{\cO_F}:=(M+\tau_M(\tau^*M))\cap M_{\cO_F}[\fj^{-1}].
\end{equation}

In \loccit we introduced the $A$-submodule of \emph{integral} and \emph{regulated} extensions of $\mathbbm{1}$ by $\underline{M}$ in $\AMot_F$, denoted $\Ext^{1,\text{reg}}_{\cO_F}(\mathbbm{1},\underline{M})$, and proved that the map
\begin{equation}\label{eq:iota}
\iota:\frac{N_{\cO_F}}{(\id-\tau_M)(M_{\cO_F})}\longrightarrow \Ext^{1,\text{reg}}_{\cO_F}(\mathbbm{1},\underline{M}),
\end{equation}
assigning to the class of $m\in N_{\cO_F}\subset M[\fj^{-1}]$ the extension whose middle object has underlying module $M\oplus (A\otimes F)$ and $\tau$-morphism $\left(\begin{smallmatrix} \tau_M & m \\ 0 & 1 \end{smallmatrix}\right)$ (with the obvious arrows), is a natural isomorphism of $A$-modules \cite[Theorem~E and Corollary~5.5]{gazda}. We also formulated a conjecture \cite[Conjecture 5.8]{gazda} relating the above to the submodule of regulated extensions having everywhere good reduction. 

Some computations suggested that $\Ext^{1,\text{reg}}_{\cO_F}(\mathbbm{1},\underline{M})$ is generally \emph{not} finitely generated (see below), contrary to what is expected in the number field setting. For instance, if $\underline{M}$ arises as the dual of the motive of an abelian Anderson module $E$ with everywhere good reduction, then one can show, \cf \cite{gazda-companion}, that up to a twist by the module of K\"ahler differentials it is isomorphic to $E(\cO_F)$ (which is not finitely generated over $A$, as Poonen proved in the case of Drinfeld modules \cite{poonen}).

Assuming that $\underline{M}$ is rigid analytically trivial, we will demonstrate that the defect of finite generation is measured by the morphism $r_B$ introduced in \eqref{eq:rB(N,M)}:
\begin{equation}\label{eq:rBetti}
\Betti{r}(\underline{M}):\Ext^{1,\mathrm{reg}}_{\cO_F}(\mathbbm{1},\underline{M}) \longrightarrow \operatorname{H}^1(G_\infty,\Betti{\underline{M}}).
\end{equation}
It assigns to an extension of rigid analytically trivial $A$-motives the class of the continuous cocycle associated with the induced extension of $A$-linear representations of $G_\infty$. 

For later reference, let us explain how $\Betti{r}(\underline{M})$ is computed. Choose $m\in N_{\cO_F}$ and let $[\underline{E}]\in \Ext^{1,\text{reg}}_{\cO_F}(\mathbbm{1},\underline{M})$ be the class of the extension $\iota(m)$ given by \eqref{eq:iota}. The group underlying the Betti realization of $\underline{E}$ consists of pairs $(\xi,a)$, with $\xi \in M\otimes_{A\otimes F}\CI\langle A \rangle$ and $a\in \CI\langle A \rangle$, satisfying
\begin{equation}
\begin{pmatrix} \tau_M & m \\ 0 & 1 \end{pmatrix} \begin{pmatrix}\tau^*\xi \\ \tau^*a \end{pmatrix}=\begin{pmatrix}\xi \\ a \end{pmatrix}. \nonumber
\end{equation}
The bottom row implies $a\in A$, and the top row gives $\xi-\tau_M(\tau^*\xi)=am$. A splitting of the sequence $0\to \Betti{\underline{M}}\to \Betti{\underline{E}}\to A \to 0$ in the category of $A$-modules corresponds to a choice of a solution $\xi_m\in M\otimes_{A\otimes F}\CI\langle A \rangle$ to the equation $\xi-\tau_M(\tau^*\xi)=m$. Such a choice yields a decomposition
\begin{equation}\label{eq:mu}
\mu_B:\Betti{\underline{M}}\oplus \Betti{\mathbbm{1}}\stackrel{\sim}{\longrightarrow} \Betti{\underline{E}}, \quad (\omega,a)\longmapsto (\omega+a\xi_m,a). \nonumber
\end{equation}
Under this description, an element $\sigma\in G_\infty$ acts on $\Betti{\underline{E}}$ by 
\begin{equation}
(\omega+a\xi_m,a)\mapsto (\,^{\sigma}\omega+a\,^{\sigma}\xi_m,a)=(\,^{\sigma}\omega+a(\,^{\sigma}\xi_m-\xi_m)+a\xi_m,a), \nonumber
\end{equation}
where $\,^{\sigma}\xi_m-\xi_m\in \Betti{\underline{M}}$. Hence, $\sigma$ acts as the matrix $\left(\begin{smallmatrix} \sigma & \, \,^{\sigma}\xi_m-\xi_m \\ 0 & 1 \end{smallmatrix}\right)$. We deduce the following.

\begin{Lemma}\label{lem:compute-rBM}
$\Betti{r}(\underline{M})$ maps the class of the extension $\iota(m)$ to the cohomology class of the cocycle $\sigma\mapsto \,^{\sigma}\xi_m-\xi_m$, where $\xi_m\in M\otimes_{A\otimes F} \CI\langle A \rangle$ is any solution to $\xi_m-\tau_M(\tau^*\xi_m)=m$. 
\end{Lemma}

We give names to the kernel and cokernel of $\Betti{r}(\underline{M})$.
\begin{Definition}[Extensions with analytic reduction \& Class module]
We define the following $A$-modules:
\begin{equation}\label{eq:class-module}
\Ext^{1,\mathrm{reg}}_{\cO_F,\infty}(\mathbbm{1},\underline{M}):=\ker \Betti{r}(\underline{M}), \quad \operatorname{Cl}(\underline{M}):=\coker \Betti{r}(\underline{M}).
\end{equation}
We call the former \emph{the module of integral, regulated extensions having analytic reduction at $\infty$} and the latter \emph{the class module of $\underline{M}$}. 
\end{Definition}

\begin{Remark}
If the $A$-motive $\underline{M}$ arises as the dual of the motive of a uniformizable and abelian Anderson module~$E$ with everywhere good reduction over $F$, we prove in the companion paper \cite{gazda-companion} that, up to a twist by the module of K\"ahler differentials, $\Ext^{1,\text{reg}}_{\cO_F,\infty}(\mathbbm{1},\underline{M})$ recovers Taelman’s unit module $U(E/\cO_F)$ and $\operatorname{Cl}(\underline{M})$ recovers its class module $H(E/\cO_F)$ (see \cite{taelman-dirichlet}). In particular, $A$-motivic cohomology naturally generalizes Taelman's theory to general $A$-motives, without assuming effectivity or everywhere good reduction. 
\end{Remark}

We turn to the computation of the modules defined in \eqref{eq:class-module}. To do so, we introduce a cochain complex $G_{\underline{M}}$, which will later be compared to coherent cohomology of shtuka models in Section~\ref{sec:proof}. It deserves its own definition.
\begin{Definition}\label{def:Gm}
Let $G_{\underline{M}}$ denote the complex of $A$-modules sitting in degrees $0$ and $1$:
\[
G_{\underline{M}}=\left[\frac{M\otimes_{A\otimes F} \FI\langle A \rangle}{M_{\cO_F}}\xrightarrow{\id-\tau_M} \frac{M\otimes_{A\otimes F} \FI\langle A \rangle}{N_{\cO_F}} \right]
\]
where the differential is induced by $m\mapsto m-\tau_M(\tau^*m)$ on $M\otimes_{A\otimes F} \FI\langle A \rangle$. 
\end{Definition}

We have a commutative diagram of $A$-modules:
\begin{equation}\label{diag:relation-complex}
\begin{tikzcd}
0 \arrow[r] & M_{\cO_F} \arrow[d,"\id-\tau_M"]\arrow[r] & M\otimes_{A\otimes F} \FI\langle A \rangle \arrow[r]\arrow[d,"\id-\tau_M"] & G_{\underline{M}}^0 \arrow[r]\arrow[d,"d^1"] & 0 \\
0 \arrow[r] & N_{\cO_F} \arrow[r] & M\otimes_{A\otimes F} \FI\langle A \rangle \arrow[r] & G_{\underline{M}}^1 \arrow[r] & 0
\end{tikzcd}
\end{equation}
where we denote by $G_{\underline{M}}^i$ the $i$th term of the complex $G_{\underline{M}}$, and by $d^i$ its differentials. Note that both rows are exact; for the top one, this follows by the flatness of $M_{\cO_F}$. The bottom row requires an additional argument to show that $N_{\cO_F}\to M\otimes_{A\otimes \cO_F} F_{\infty}\langle A \rangle$ is injective, as $N_{\cO_F}$ may not be flat. But this follows by noticing that there exists an integer $n\geq 0$ such that $N_{\cO_F}\subset \fj^{-n}M_{\cO_F}$, and that $\fj^{-n}M_{\cO_F}$ embeds in $M\otimes_{A\otimes \cO_F} F_{\infty}\langle A\rangle$ by flatness of $M_{\cO_F}$.

The following long exact sequence will be \emph{fundamental} for the sequel.
\begin{Proposition}\label{prop:motivic-exact-sequence}
The snake lemma applied to the diagram \eqref{diag:relation-complex} yields a long exact sequence of $A$-modules:
\[
0\to \Hom_{\AMot_F}(\mathbbm{1},\underline{M})\to \Betti{\underline{M}}^+\to \operatorname{H}^0(G_{\underline{M}})\to \Ext^{1,\operatorname{reg}}_{\cO_F}(\mathbbm{1},\underline{M}) \xrightarrow{r_B(\underline{M})} \operatorname{H}^1(G_\infty,\Betti{\underline{M}}) \to \operatorname{H}^1(G_{\underline{M}}) \to 0.
\]
\end{Proposition}
\begin{proof}
The complex $[M_{\cO_F}\xrightarrow{\id-\tau_M} N_{\cO_F}]$ has $\Hom_{\AMot_F}(\mathbbm{1},\underline{M})$ and $\Ext^{1,\operatorname{reg}}_{\cO_F}(\mathbbm{1},\underline{M})$ as its $0$th and $1$st cohomology modules (see \cite[Proposition 4.23]{gazda} and \eqref{eq:iota}). The kernel and cokernel of the middle vertical arrow in the diagram \eqref{diag:relation-complex} are computed by Theorem \ref{thm-analytically}. Hence everything is clear except, perhaps, that the map $r_B(\underline{M})$ is indeed the one that appears as the connecting homomorphism. To see this, first observe that Theorem \ref{thm-analytically} gives an isomorphism of $A$-modules:
\begin{equation}\label{eq:H1-explicit}
\frac{M\otimes_{A\otimes F} \FI\langle A \rangle}{(\id-\tau_M)(M\otimes_{A\otimes F} \FI\langle A \rangle)} \stackrel{\sim}{\longrightarrow} \operatorname{H}^1(G_\infty,\Betti{\underline{M}}).
\end{equation}
The above map sends the class of $f\in M\otimes_{A\otimes F}\FI\langle A \rangle$ to the class of the cocycle $c_f:\sigma \mapsto \,^\sigma\xi_f-\xi_f$, where $\xi_f\in M\otimes_{A\otimes F}\FI^s\langle A \rangle$ is \emph{any} solution $\xi$ to the equation $\xi-\tau_M(\tau^*\xi)=f$. That $r_B(\underline{M})$ agrees with the connecting homomorphism is the content of Lemma \ref{lem:compute-rBM}.
\end{proof}

In particular, we obtain convenient descriptions of the modules of interest in terms of $G_{\underline{M}}$:
\begin{Corollary}\label{cor:class-ext-GM}
We have $\Cl(\underline{M})\cong H^1(G_{\underline{M}})$ and $\Ext^{1,\mathrm{reg}}_{\cO_F,\infty}(\mathbbm{1},\underline{M})\cong H^0(G_{\underline{M}})/\underline{M}_B^+$.
\end{Corollary}

Below, we will show that the complex $G_{\underline{M}}$ can be interpreted in terms of coherent cohomology of shtuka models of $\underline{M}$. In particular, $G_{\underline{M}}$ is a perfect complex and its cohomology groups are finitely generated; by Corollary \ref{cor:class-ext-GM}, this will imply the finite generation of $\Cl(\underline{M})$ and $\Ext^{1,\mathrm{reg}}_{\cO_F,\infty}(\mathbbm{1},\underline{M})$.

\subsection{Hodge-Pink regulator}\label{subsec:regulator}
Assume that $\underline{M}$ is a rigid analytically trivial $A$-motive over $F$. We move to the definition of the regulator of $\underline{M}$. Let $\underline{H}^+$ denote the Hodge-Pink structure $\sH^+(\underline{M})$ attached to $\underline{M}$. In view of Proposition \ref{prop:exactness-HodgePink}, the Hodge-Pink realization functor $\sH^+$ is exact and hence induces an $A$-linear morphism between the corresponding extension groups:
\begin{equation}\label{eq:general-full-regulator}
r_{\sH^+}(\underline{M}):\Ext^{1}_{\AMot_F^{\text{rat}}}(\mathbbm{1},\underline{M}) \longrightarrow \Ext^{1}_{\HPk^+}(\mathbbm{1}^+,\underline{H}^+).
\end{equation}
By Definition \ref{def:regulated}, $r_{\sH^+}(\underline{M})$ maps regulated exact sequences in $\AMot_F^{\text{rat}}$ to Hodge-additive strict exact sequences in $\HPk^+$. Likewise, extensions of $A$-motives having analytic reduction at $\infty$ are sent to extensions having analytic reduction (Definition~\ref{def:analytic-reduction-HP}). In particular, \eqref{eq:general-full-regulator} induces a map
\begin{equation}\label{eq:HP-regulator}
r_{\sH^+}(\underline{M}):\Ext^{1,\mathrm{reg}}_{\cO_F,\infty}(\mathbbm{1},\underline{M}) \longrightarrow \Ext^{1,\mathrm{ha}}_{\infty}(\mathbbm{1}^+,\underline{H}^+).
\end{equation}

Recall that the target is a finite-dimensional $\KI$-vector space, and that the source is a finitely generated $A$-module (as a consequence of Proposition \ref{prop:ge-in-terms-of-cohomology} below). This provides one possible function field counterpart of Beilinson's regulator.

\begin{Definition}[Hodge-Pink regulator]\label{def:regulator}
We call \emph{the Hodge-Pink regulator of $\underline{M}$}, and denote it by $\Reg(\underline{M})$, the map \eqref{eq:HP-regulator}.
\end{Definition}

We now explain how the Hodge-Pink regulator can be computed in terms of solutions to $\tau$-difference equations. Consider the $A$-module
\begin{equation}
\Delta_{\underline{M}}:=\left\{\xi \in M\otimes_{A\otimes F} \FI\langle A \rangle~\big|~\xi-\tau_M(\tau^* \xi)\in N_{\cO_F} \right\}. \nonumber
\end{equation}
Informally, it consists of ``analytic functions'' $\xi$ with values in $M_{\FI}$ which, after applying the $\tau$-difference operator $f\mapsto f-\tau_M(\tau^*f)$, become ``integral''. Note that, by the analytic continuation theorem (Theorem \ref{thm:analytic-continuation}) applied to extensions of $\mathbbm{1}$ by $\underline{M}$, elements of $\Delta_{\underline{M}}$ extend to $M\otimes_{A\otimes F} \mathcal{M}_J(\mathfrak{U}_{\FI})$ and therefore can be regarded as elements of $M\otimes_{A\otimes F} \FI(\!(\fj)\!)$ by the Laurent expansion map of Construction \ref{cons:Laurent-series-expansion-map}. Also note that elements $\xi$ of either $M_{\cO_F}$ or $\underline{M}_B^+$ obviously satisfy $\xi-\tau_M(\tau^* \xi)\in N_{\cO_F}$ and thus belong to $\Delta_{\underline{M}}$; moreover, we have an isomorphism
\begin{equation}\label{eq:ext-mot-iso}
\frac{\Delta_{\underline{M}}}{\underline{M}_B^++M_{\cO_F}} \stackrel{\sim}{\longrightarrow} \Ext^{1,\mathrm{reg}}_{\cO_F,\infty}(\mathbbm{1},\underline{M}), \quad \xi\longmapsto \iota\!\left(\xi-\tau_M(\tau^* \xi)\right).
\end{equation}
Indeed, by definition of $G_{\underline{M}}$ we have $H^0(G_{\underline{M}})\cong \Delta_{\underline{M}}/M_{\cO_F}$, so this follows from Corollary \ref{cor:class-ext-GM}. The powerful idea of computing extension modules in terms of solutions to $\tau$-difference equations goes back to the work of Sinha \cite{sinha}.

On the other hand, the space of Hodge-additive extensions having analytic reduction has been computed in Section \ref{sec:rat-mix-A-mot}, Theorem \ref{thm:extensions-of-MHPS-arising-from-motives}. Under the assumption that the weights of $\underline{M}$ are all negative, we have
\begin{equation}\label{eq:ext-hdg-iso}
\frac{(M+\tau_M(\tau^*M))\otimes_{A\otimes F} \FI[\![\fj]\!]}{(\underline{M}_B^+)_{\KI}+M\otimes_{A\otimes F}\FI[\![\fj]\!]}\stackrel{\sim}{\longrightarrow}\Ext^{1,\mathrm{ha}}_{\infty}(\mathbbm{1}^+,\underline{H}^+), \quad \psi \longmapsto \varphi^+(\psi).
\end{equation}
The Hodge-Pink regulator may then be computed as follows.

\begin{Proposition}\label{prop:explicit-computation-regulator}
Under the isomorphisms \eqref{eq:ext-mot-iso} and \eqref{eq:ext-hdg-iso}, the Hodge-Pink regulator is induced by
\[
\Delta_{\underline{M}}\longrightarrow (M+\tau_M(\tau^* M))\otimes_{A\otimes F} \FI[\![\fj]\!], \quad \xi \longmapsto -\xi.
\]
\end{Proposition}

\begin{proof}
Start with $\xi\in \Delta_{\underline{M}}$ and set $m=\xi-\tau_M(\tau^* \xi)$ in $N_{\cO_F}[\fj^{-1}]$. Let $[\underline{E}]$ be the class of the extension $\iota(m)$. We denote by $\fq_E$ the Hodge-Pink lattice of its associated Hodge-Pink structure. We must show that the inclusion $\fq_E\subset \underline{E}_B\otimes_{A,\nu} \FI^s(\!(\fj)\!)$ is isomorphic, via the map $\mu_B$ of \eqref{eq:mu}, to
\[
\begin{pmatrix} \id_M & -\gamma_{\underline{M}}(\xi) \\ 0 & 1 \end{pmatrix}\bigl(\fq_{\underline{M}}\oplus \fq_{\mathbbm{1}}\bigr)\hookrightarrow
(\underline{M}_B\oplus \mathbbm{1}_B)\otimes_{A,\nu} \FI^s(\!(\fj)\!).
\]
By definition, $\fq_E=\gamma_{\underline{E}}^{-1}(E\otimes_{A\otimes F} \FI^s[\![\fj]\!])$. The claim then follows from the commutativity of the diagram
\[
\begin{tikzcd}[ampersand replacement=\&, column sep=4em]
(\underline{M}_B\oplus \mathbbm{1}_B)\otimes_{A,\nu} \FI^s(\!(\fj)\!) \arrow[r,"\mu_B"]\arrow[d,"\gamma_{\underline{M}}\oplus \gamma_{\mathbbm{1}}"] \&
\underline{E}_B\otimes_{A,\nu} \FI^s(\!(\fj)\!) \arrow[d,"\gamma_{\underline{E}}"] \\  
\bigl(M\otimes_{A\otimes F}\FI^s(\!(\fj)\!)\bigr)\oplus \FI^s(\!(\fj)\!) \arrow[r,"{\left(\begin{smallmatrix} \id_M & \xi \\ 0 & 1 \end{smallmatrix}\right)}"] \&
E\otimes_{A\otimes F} \FI^s(\!(\fj)\!)
\end{tikzcd}
\]
together with the identity
\[
\fq_{\underline{M}}\oplus \fq_{\mathbbm{1}}
=(\gamma_{\underline{M}}\oplus \gamma_{\mathbbm{1}})^{-1}\!\Bigl(\bigl(M\otimes_{A\otimes F}\FI^s[\![\fj]\!]\bigr)\oplus \FI^s[\![\fj]\!]\Bigr).
\qedhere
\]
\end{proof}

\section{Shtuka models}
Let $F$ be a finite extension of $K$ which is totally ramified at $\infty$ (see Reduction \ref{red:F-totally-ramified}); let $X$ be a smooth projective curve over $\bF$ whose function field is $F$. Let $\underline{M}$ be an $A$-motive over $F$ (we will not require rigid analytic triviality here). In this section, we associate non-canonically to $\underline{M}$ a \emph{shtuka model} over the surface $(\Spec A)\times X$ and, when all the weights of $\underline{M}$ are negative, on $C\times X$. Roughly speaking, these are data $(\cM,\cN,\tau_{\cM})$ where $\cM\hookrightarrow \cN$ is an inclusion of coherent sheaves on $\Spec A\times X$ (respectively $C\times X$), and where $\tau_{\cM}:\tau^*\cM\to \cN$ is a morphism which restricts to $\tau_M:\tau^*M_{\cO_F}\to N_{\cO_F}$ on the affine open $\Spec(A\otimes \cO_F)$ of $C\times X$, together with a vanishing condition along the complement.

This powerful technique was, to the author's knowledge, first introduced by Pink in \cite{pink-isogenies} as a key ingredient in order to prove the isogeny conjecture for $A$-motives. It was also reintroduced by Mornev \cite[\S 12]{mornev-shtuka} in the context of Drinfeld modules with everywhere good reduction. An incarnation of this construction also appears in \cite[Definition 1.13]{fang} in the setting of Anderson $t$-modules. Although our motivations owe much to Mornev's work, our definition of $C\times X$-shtuka models differs: we require $\cN$ to be torsion free rather than locally free, since the $A\otimes \cO_F$-module $N_{\cO_F}$ of \eqref{def:NA} need not be locally free. The definition presented below in Definition \ref{def:shtuka-model CC} has the convenient feature that an existence result follows from the sole assumption that the weights of $\underline{M}$ are negative (Proposition \ref{prop:existence-CxX-shtuka-models}). 

An important step towards the proof of the main theorems is the existence of a quasi\nobreakdash-isomorphism
\begin{equation}\label{eq:step1}
G_{\underline{M}}\stackrel{\sim}{\longrightarrow} \operatorname{R}\!\Gamma\!\left(\Spec A \times X, \cM\xrightarrow{\iota-\tau_{\cM}}\cN\right)\quad (\text{Proposition~}\ref{prop:ge-in-terms-of-cohomology}).
\end{equation}
Since $X$ is proper over $\bF$, this implies the perfectness of $G_{\underline{M}}$. Corollary \ref{cor:class-ext-GM} then shows that both modules $\Ext^{1,\operatorname{reg}}_{\cO_F,\infty}(\mathbbm{1},\underline{M})$ and $\operatorname{Cl}(\underline{M})$ are finitely generated, completing the proof of the first part of Theorem \ref{thm:finiteness-motcoh}. Theorem \ref{thm:rank-dim} and the second part of Theorem \ref{thm:finiteness-motcoh} follow from the further observation that $C\times X$-shtuka models, when restricted to the formal neighbourhood of $\{\infty\}\times \{\infty\}$, are related to the extension space $\Ext^{1,\operatorname{ha}}_{\infty}(\mathbbm{1}^+,\sH^+(\underline{M}))$. We study this relation in Subsection \ref{subsec:shtuka-models-MHPS}. 

Let us introduce some notation. We still denote by $\tau:C\times X\to C\times X$ the morphism of $\bF$-schemes which acts as the identity on the left-hand factor $C$ and as the absolute $q$-Frobenius on the right-hand one. Since $C$ is separated over $\bF$, the graph morphism $X\to C\times X$, where $X\to C$ is induced by $K\subset F$, is a closed immersion; its image defines a closed subscheme $\Delta\subset C\times X$ of codimension $1$. We regard $\Delta$ as an effective divisor on $C\times X$. The evaluation of $\cO(-\Delta)$ on the affine open subscheme $\Spec(A\otimes \cO_F)\subset C\times X$ recovers the ideal $\fj\subset A\otimes \cO_F$.

We also borrow notation from \cite{gazda}. Let $\OI$ be the valuation ring of $\KI$ and $\fm_{\infty}$ its maximal ideal. For $R$ a noetherian $\bF$-algebra, $\cA_{\infty}(R)$ denotes the ring
\[
\cA_{\infty}(R)=\varprojlim_n (\OI\otimes R)/(\fm_{\infty}^n\otimes R).
\]
This ring was used in \emph{loc.\ cit.} to define isocrystals and mixedness. Let also $\cB_{\infty}(R)$ be the ring $\KI\otimes_{\OI}\cA_{\infty}(R)$. More generally, for an $A\otimes R$-module $M$ we will write $\cB_{\infty}(M)$ for $M\otimes_{A\otimes R} \cB_{\infty}(R)$. Geometrically, the formal spectrum $\Spf \cA_{\infty}(R)$ corresponds to a formal neighbourhood of the closed subscheme $\{\infty\}\times \Spec R$ in $C\times \Spec R$.

Let $\bF_{\infty}$ denote the residue field of $\KI$ and let $\pi_{\infty}$ be a uniformizer. We have an identification
\begin{equation}\label{eq:identification-Binfty}
\cB_{\infty}(R)=(\mathbb{F}_{\infty}\otimes R)(\!(\pi_{\infty}\otimes 1)\!)
\end{equation}
where the subring $(\mathbb{F}_{\infty}\otimes R)[\![\pi_{\infty}\otimes 1]\!]$ corresponds to $\cA_{\infty}(R)$. In particular, the valuation $v_{\infty}$ on $\KI$ extends to a function $v_{\infty} :\cB_{\infty}(R)\to \bZ\cup \{\infty\}$ (which is only submultiplicative). 

We also denote by $\infty$ the unique infinite place of $F$. The closed subscheme $C\times \{\infty\}$ defines an effective divisor on $C\times X$ which we denote by $(0,\infty)$. Similarly, we let $(\infty,0)$ denote the effective divisor $\{\infty\}\times X$.

\subsection{$X$-shtuka models}\label{subsec:C-shtuka}
Let $\underline{M}$ be an $A$-motive over $F$. Let $M_{\cO_F}$ denote its maximal integral model \cite[\S 4]{gazda}; by Theorem 4.32 in \loccit $M_{\cO_F}$ is a locally free $A\otimes \cO_F$-module. We set 
\[
N:=M+\tau_M(\tau^*M) \quad \text{and} \quad N_{\cO_F}:=(M+\tau_M(\tau^*M))\cap M_{\cO_F}[\fj^{-1}].
\]
Note that there is no particular reason for $N_{\cO_F}$ to be locally free. However, it is finitely generated (as a submodule of $M_{\cO_F}\otimes_{A\otimes \cO_F} \fj^{-e}$ for some large enough $e$). The map $\tau_M:\tau^* M\to M[\fj^{-1}]$ factors as $\tau_M:\tau^*M_{\cO_F}\to N_{\cO_F}$. 

\begin{Definition}[$X$-shtuka models]\label{def:shtuka-model-C}
An \emph{$X$-shtuka model} $\underline{\cM}$ for $\underline{M}$ is the datum $(\cN,\cM,\tau_{\cM})$ of
\begin{enumerate}[label=$(\alph*)$]
\item a torsion free coherent sheaf $\cN$ on $(\Spec A)\times X$ such that $\cN(\Spec(A\otimes \cO_F))=N_{\cO_F}$,
\item a locally free coherent sheaf $\cM\subset \cN$ such that $\cM(\Spec(A\otimes \cO_F))=M_{\cO_F}$ and such that the cokernel of the inclusion $\iota:\cM\to \cN$ is supported on $\Delta$,
\item a morphism $\tau_{\cM}:\tau^*\cM\to \cN(-(0,\infty))$ which coincides with $\tau_M:\tau^*M_{\cO_F}\to N_{\cO_F}$ on the affine open subscheme $\Spec(A\otimes \cO_F)$.
\end{enumerate}
\end{Definition}

\begin{Proposition}\label{prop:existence-of-X-models}
An $X$-shtuka model for $\underline{M}$ exists.
\end{Proposition}

\begin{proof}
Let $U:= \Spec \cO_F$ and $V\subset X$ be an affine subscheme over $\bF$ containing the point $\infty$. Then $\{U,V\}$ forms a cover of $X$ in the Zariski topology. Let $B\subset F$ be the $\bF$-subalgebra such that $V=\Spec B$ and let $D\subset F$ be the subalgebra such that $W:=U\cap V=\Spec D$. Write $\fj_S$ for the ideal $\cO(-\Delta)(\Spec(A\otimes S))$ of $A\otimes S$, where $\Spec S$ is an affine Zariski open subscheme of $X$.

We choose a locally free coherent sheaf $\cM'$ on $(\Spec A)\times X$ which extends $M_{\cO_F}$ on the dense open subscheme $\Spec A\otimes \cO_F$; in practice, let $M_D$ be the localization of $M_{\cO_F}$ on $(\Spec A)\times W$ and let $M_{B}'$ be a locally free $A\otimes B$-lattice in $M_D$. Then $\cM'$ is the unique coherent sheaf extending $M_{\cO_F}$ such that $\cM'(\Spec A\otimes B)=M'_B$.

Since $\tau_M(\tau^*M_{\cO_F})\subset M_{\cO_F}[\fj^{-1}]$, we have $\tau_M(\tau^*M_D)\subset M_D[\fj_D^{-1}]$. However, it need not be true that $\tau_M(\tau^*M_B')\subset M_B'[\fj_B^{-1}]$; in particular, $\tau_M$ does not yet extend to a map $\tau_{\cM'}$ away from the divisor $\Delta$. Still, there exists $d\in B$, invertible in $D$, such that 
\begin{equation}
\tau_M(\tau^*M_B')\subset d^{-1}M_B'[\fj_B^{-1}]. \nonumber
\end{equation}
Choose $r\in B$, invertible in $D$, which vanishes\footnote{Such an $r$ always exists: the divisor $E:=\deg(x)\cdot\infty-\deg(\infty)\cdot x$ has degree zero, hence $nE$ is principal for $n$ some positive integer $n$ (since $\mathrm{Pic}^0(X)$ is finite; see \cite[Lemma~5.6]{rosen}). Choosing $r$ with $(r)=nD$ yields $r\in B$, and $r$ is invertible in $D$.} at $\infty$, and set $M_B:=(rd)M_B'$. Then
\begin{equation}
\tau_M(\tau^*M_B)\subset r\,M_B[\fj_B^{-1}]. \nonumber
\end{equation}
Since $r$ is invertible in $D$, multiplication yields \emph{gluing} isomorphisms 
\begin{equation}\label{eq:glueing-M}
M_{\cO_F}\otimes_{\cO_F} D\stackrel{=}{\longrightarrow} M_D \stackrel{\sim}{\longleftarrow} M_B\otimes_B D. 
\end{equation}
So let $\cM$ be the unique (locally free) coherent sheaf extending $M_{\cO_F}$ such that $\cM(\Spec A\otimes B)=M_B$.

For an open $\Spec S\subset X$, set $N_S:=(M+\tau_M(\tau^*M))\cap M_S[\fj_S^{-1}]$. Then $N_S$ is a finitely generated, torsion-free $A\otimes S$-module containing $M_S$. We have gluing isomorphisms
\begin{equation}\label{eq:glueing-N}
N_{\cO_F}\otimes_{\cO_F} D \stackrel{\sim}{\longrightarrow} N_D \stackrel{\sim}{\longleftarrow} N_B\otimes_B D,
\end{equation}
so let $\cN$ be the (torsion free) coherent sheaf on $\Spec A\times X$ obtained from the Zariski gluing \eqref{eq:glueing-N}. Since $M_{\cO_F}\subset N_{\cO_F}$ and $M_B\subset N_B$, we get $\cM\subset \cN$. Moreover $M_{\cO_F}[\fj^{-1}]=N_{\cO_F}[\fj^{-1}]$ and $M_B[\fj_B^{-1}]=N_B[\fj_B^{-1}]$, hence the cokernel of $\cM\subset \cN$ is supported on $\Delta$.

Finally, since $\tau_M(\tau^*M_S)\subset N_S$ for all $S\in\{\cO_F,B,D\}$, the maps glue to a morphism of $\cO_{(\Spec A)\times X}$\nobreakdash-modules $\tau_{\cM}:\tau^*\cM\to \cN$. Since $\tau_M(\tau^*M_B)\subset rN_B$ and $r$ vanishes at $\infty$, we in fact have $\tau_{\cM}(\tau^*\cM)\subset \cN(-(0,\infty))$.
\end{proof}

\begin{Example}[Shtuka models of Carlitz twists]\label{ex:shtuka-moduel-carlitz-twists}
Consider the Carlitz $n$th twist of Notation \ref{not:t-setting}:
\[
\underline{A}(n):=\bigl(K[t],(t-\theta)^{-n}\tau\bigr).
\] 
Let $m:=\max(0,n)$ and choose $k\in \bZ$ such that $k\leq \lfloor \frac{n-1}{q-1}\rfloor$. Set $\cM:=\cO_{\Spec \bF[t]\times \bP^1}(k\cdot (0,\infty))$ and $\cN:=\cM(m\cdot \Delta)$. Then $\tau^*\cM\cong \cO_{\Spec \bF[t]\times \bP^1}(qk\cdot (0,\infty))$, and one checks that the condition on $k$ ensures that $(t-\theta)^{-n}\tau^*\cM\subset \cN(-(0,\infty))$. In particular, the triple $(\cN,\cM,\times (t-\theta)^{-n}\tau)$ defines an $\bP^1$-shtuka model for $\underline{A}(n)$.

In fact, all $\bP^1$-shtuka models for $\underline{A}(n)$ arise in this manner; in particular, there is a maximal one, namely the one for which $k:=\lfloor \tfrac{n-1}{q-1}\rfloor$. More generally, one can show that there is a unique maximal $X$-shtuka model for any $A$-motive (we will not use this fact).
\end{Example}

Let $\underline{\cM}=(\cN,\cM,\tau_{\cM})$ be an $X$-shtuka model for $\underline{M}$. The condition that the image of $\tau_{\cM}$ lands in $\cN(-(0,\infty))$ is crucial for the cohomological considerations to come. This is often materialized by the following lemma.

\begin{Lemma}\label{lem:isomorphism-at-infty}
Let $i:\Spec \cO_{\FI}\langle A \rangle \to \Spec(A\otimes \cO_{\FI}) \hookrightarrow (\Spec A)\times X$ be the canonical composition of $A$-schemes. Then $i^*\cM=i^*\cN$, and the induced morphism
\begin{equation}
\iota-\tau_{\cM}:i^*\cM(\Spec \cO_{\FI}\langle A\rangle)\longrightarrow i^*\cN(\Spec \cO_{\FI}\langle A\rangle) \nonumber
\end{equation}
is an isomorphism of $A$-modules.
\end{Lemma}

\begin{proof}
By Lemma \ref{lem:j-invertible-in-Tate-algebra}, we have $\fj\cO_{\FI}\langle A\rangle=\cO_{\FI}\langle A\rangle$. In particular, $i^*\Delta$ is the empty divisor on $\Spec \cO_{\FI}\langle A\rangle$. Hence $i^*\cM=i^*\cN$. 

Let $\pi_\infty$ be a uniformizer of $\cO_{\FI}$. Set $\Xi:=i^*\cM(\Spec \cO_{\FI}\langle A \rangle)$. Since $\tau_{\cM}(\tau^*\cM)\subset \cN(-(0,\infty))$, we have $\tau_{\cM}(\tau^*\Xi)\subset \pi_\infty \Xi$. In particular, for each $\xi\in \Xi$, the series
\begin{equation}
\psi:=\sum_{n=0}^{\infty}{\tau_{\cM}^n(\tau^{n*}\xi)} \nonumber
\end{equation}
converges in $\Xi$. The assignment $\xi\mapsto \psi$ defines an inverse to $\id-\tau_{\cM}$ on $\Xi$, hence $\iota-\tau_{\cM}$ is an isomorphism.
\end{proof}

\subsection{$C\times X$-shtuka models}\label{subsec:CxX-shtuka-models}
We now extend the construction of Proposition \ref{prop:existence-of-X-models} from $(\Spec A)\times X$ to $C\times X$.

\begin{Definition}[$C\times X$-shtuka model]\label{def:shtuka-model CC}
A \emph{$C\times X$-shtuka model} $\underline{\cM}$ for $\underline{M}$ is the datum $(\cN,\cM,\tau_{\cM})$ of
\begin{enumerate}[label=$(\alph*)$]
\item\label{item:CC1} a torsion free coherent sheaf $\cN$ on $C\times X$ such that $\cN(\Spec(A\otimes \cO_F))=N_{\cO_F}$,
\item\label{item:CC2} a locally free coherent subsheaf $\cM\subset \cN$ such that $\cM(\Spec(A\otimes \cO_F))=M_{\cO_F}$ and such that the cokernel of the inclusion $\iota:\cM\to \cN$ is supported on $\Delta$,
\item\label{item:CxC-shtuka-model-zero} a morphism of sheaves
\[
\tau_{\cM}:\tau^*\cM\longrightarrow \cN(-(0,\infty))
\]
which coincides with $\tau_M:\tau^*M_{\cO_F}\to N_{\cO_F}$ on $\Spec(A\otimes \cO_F)$,
\item\label{item:CC4} the localisation of $\tau_{\cM}$ at the divisor $(\infty,0)$ induces a map
\[
\tau_{\cM}:\tau^*\cM_{(\infty,0)}\longrightarrow \cM_{(\infty,0)},
\]
and there exists an integer $h>0$ such that
\[
\tau_{\cM}^h\bigl(\tau^{h*}\cM_{(\infty,0)}\bigr)\subset \cM_{(\infty,0)}\bigl(-(\infty,0)\bigr).
\]
\end{enumerate}
\end{Definition}

\begin{Remark}
Clearly, the restriction of a $C\times X$-shtuka model for $\underline{M}$ to $(\Spec A)\times X$ is an $X$-shtuka model for $\underline{M}$.
\end{Remark}

\begin{Example}[Shtuka models of Carlitz twists, resumed]\label{ex:shtuka-model-carlitz-resumed}
We continue Example \ref{ex:shtuka-moduel-carlitz-twists}. Among the $\bP^1$-shtuka models described there, only those with $n>0$ extend to $\bP^1\times \bP^1$-shtuka models: if $n<0$, the obstruction comes from the fact that $(t-\theta)^{-n}$ has a pole at $(\infty,0)$ of order $-n$. For $n=0$, the strict contraction condition at $(\infty,0)$ does not hold. 

If $n> 0$, the datum $\cM:=\cO_{\bP^1\times \bP^1}$, $\cN:=\cM(m\Delta)$ and $\tau_{\cM}:=\times (t-\theta)^{-n}\tau$ defines a $\bP^1\times \bP^1$-shtuka model for $\underline{A}(n)$.
\end{Example}

The main result of this subsection is the following.

\begin{Proposition}\label{prop:existence-CxX-shtuka-models}
If the weights of $\underline{M}$ are negative, then a $C\times X$-shtuka model for $\underline{M}$ exists.
\end{Proposition}

Before initiating the proof of Proposition \ref{prop:existence-CxX-shtuka-models}, we supply some ingredients on function field isocrystals with negative weights. The next statement follows from Dieudonn\'e theory as in \cite[\S 3.1]{gazda}.

\begin{Lemma}\label{lem:stable-lattice-if-negative-weights}
Let $\underline{M}$ be an $A$-motive over $F$ whose weights are all negative. There exist an $\cA_{\infty}(F)$-lattice $T$ in $\cB_\infty(M)$ which is stable under $\tau_M$ and two positive integers $d$ and $h$ such that
\[
\tau_M^h\bigl(\tau^{h*}T\bigr)\subset \fm_{\infty}^d\,T.
\]
\end{Lemma}

\begin{proof}[Proof of Proposition \ref{prop:existence-CxX-shtuka-models}]
The first part of the proof consists of constructing a suitable coherent sheaf on a neighborhood of the divisor $(\infty,0)=\{\infty\}\times X$. It splits into finitely many irreducible components $D_1\cup\cdots\cup D_\delta$, where $\delta\leq\deg(\infty)$, with respective generic points $\eta_1,\ldots,\eta_\delta$. Let $\Spec B\subset C\times X$ be an affine open containing all the $\eta_i$, \eg by
\href{https://stacks.math.columbia.edu/tag/09NJ}{[09NJ]}; let also $\fp_i\subset B$ be the height-one prime ideal corresponding to $\eta_i$ and put
\[
R:=S^{-1}B
\quad\text{where}\quad
S:=B\setminus\bigcup_{i=1}^{\delta}\fp_i.
\]
Since $B$ is normal, $R$ is a semilocal Dedekind domain whose maximal ideals are the $\fp_iR$. Choose a uniformizer $\pi_\infty\in K$ at $\infty$. Each $R_{\fp_iR}=\cO_{C\times X,\eta_i}$ is a discretely valued ring with uniformizer $\pi_\infty$ and residue field a finite constant field extension $F_i$ of $F$. The completion $\widehat R$ of $R$ at $\pi_\infty R=\bigcap_i\fp_iR$ is $\cA_\infty(F)$. The morphism $\tau$ permutes the points $\eta_i$ and hence induces an endomorphism of $R$, which fixes $\pi_\infty$.

Let $L_0$ be an $R$-lattice inside $M\otimes_{A\otimes F}\bF(C\times X)$. We claim that there exists an integer $N>0$ such that
\begin{equation}\label{eq:negative-weight-lattice-proof}
\tau_M^N(\tau^{N*}L_0)\subset\pi_\infty L_0.
\end{equation}
To see this, let $T$ be the $\cA_\infty(F)$-lattice provided by Lemma \ref{lem:stable-lattice-if-negative-weights}. The completion $\widehat L_0$ of $L_0$ at $\pi_\infty$ is an $\cA_\infty(F)$-lattice in $\cB_\infty(M)$ commensurable with $T$, and for some $a\geq0$ we have
\[
\pi_\infty^aT\subset\widehat L_0
\subset\pi_\infty^{-a}T.
\]
Thus, for $r\geq0$ large enough,
\[
\tau_M^{rh}(\tau^{rh*}\widehat L_0)
\subset\pi_\infty^{rd-2a}\widehat L_0
\subset\pi_\infty\widehat L_0.
\]
Faithful flatness of $R\to\widehat R$ gives \eqref{eq:negative-weight-lattice-proof} with $N=rh$. Then set
\[
L=L_0+\tau_M(\tau^*L_0)+\ldots+
\tau_M^{N-1}(\tau^{(N-1)*}L_0).
\]
Since $\tau$ fixes $\pi_\infty$, the lattice $L$ is $\tau_M$-stable and satisfies the contraction property
\begin{equation}\label{eq:contraction-lattice}
\tau_M^N(\tau^{N*}L)\subset\pi_\infty L.
\end{equation}

Because $R$ is semilocal and $L$ is a finite projective $R$-module, it is free.
Denote by $L_B$ the $B$-submodule of $L$ generated by a finite basis over $R$, and let $\cL_B=\tilde L_B$ denote the associated locally free coherent sheaf
on $\Spec B$. Choose an $X$-shtuka model $(\cN_0,\cM_0,\tau_0)$ for $\underline M$ and consider the open subscheme $U:=(C\times X)\setminus(\{\infty\}\times X)$ of $C\times X$. Because $\cL_B|_{U\cap\Spec B}$ and $\cM_0|_{U\cap\Spec B}$ agree generically, with common generic fiber $M\otimes_{A\otimes F}\bF(C\times X)$, there exists an open neighborhood $V\subset\Spec B$ of all the $\eta_i$ on which this generic identification
induces an isomorphism $\cL_B|_{U\cap V}\cong\cM_0|_{U\cap V}$. We set $\cL:=\cL_B|_V$.
Shrinking $V$ if necessary, we may assume that it meets neither $\Delta$ nor $C\times\{\infty\}$: neither of these divisors contains an irreducible component of $\{\infty\}\times X$.

We are now in a position to construct the $C\times X$-shtuka model for $\underline M$.
The isomorphism $\cL|_{U\cap V}\cong\cM_0|_{U\cap V}$ allows us to glue $\cL$ and $\cM_0$ into a locally free coherent sheaf $\cM_W$ on $W:=U\cup V$.
As $V$ does not meet $\Delta$, we also have $\cL|_{U\cap V}\cong\cN_0|_{U\cap V}$, allowing us to construct a torsion free coherent sheaf $\cN_W$, together with an inclusion map $\iota_W:\cM_W\to\cN_W$ whose cokernel is supported at $W\cap\Delta$.

Let $j:W\to C\times X$ denote the inclusion; its complement has codimension $2$. We set $\cM:=j_*\cM_W$ and $\cN:=j_*\cN_W$. By the Hartogs equivalence for reflexive sheaves
\href{https://stacks.math.columbia.edu/tag/0EBJ}{[0EBJ]}, $\cM$ is coherent and reflexive; since $C\times X$
is a regular surface, $\cM$ is locally free
\href{https://stacks.math.columbia.edu/tag/0B3N}{[0B3N]}.
The sheaf $\cN$ is coherent by
\href{https://stacks.math.columbia.edu/tag/0AWA}{[0AWA]},
since $C\times X$ is excellent, $\cN_W$ is torsion free,
and the complement of $W$ has codimension $2$.
It is also torsion free.
Applying $j_*$ to $\iota_W$ gives an inclusion
$\iota:\cM\to\cN$ and, since $\cM_W$ and $\cN_W$
agree on $W\setminus\Delta$, their pushforwards agree
on $(C\times X)\setminus\Delta$.

It remains to construct the Frobenius morphism.
The original map $\tau_M$ induces a rational morphism
\[
\tau^*\cM\dashrightarrow\cN(-(0,\infty)).
\]
Since $\tau^{-1}(U)=U$, this morphism is regular on $U$,
where it coincides with $\tau_0$.

On $V$, the target identifies with $\cL$, since $V$
meets neither $\Delta$ nor $C\times\{\infty\}$.
Both $(\tau^*\cM)|_V$ and $\cL$ are locally free.
The rational morphism is regular at every codimension-one
point of $V$.
Indeed, this is already known at points lying in $U$.
The remaining codimension-one points are the $\eta_i$,
where regularity follows from the $\tau_M$-stability
of the semilocal lattice $L$.
Since $V$ is normal, the locally free sheaf
$\mathcal{H}om((\tau^*\cM)|_V,\cL)$ is reflexive,
and regularity in codimension one extends the rational
morphism uniquely over $V$
\href{https://stacks.math.columbia.edu/tag/0AVB}{[0AVB]}.

The resulting morphisms on $U$ and $V$ agree on
$U\cap V$, since both induce $\tau_M$ generically.
They therefore glue to a morphism
\[
j^*\tau^*\cM\longrightarrow\cN_W(-(0,\infty)).
\]
By adjunction, this induces a morphism
\[
\tau_{\cM}:\tau^*\cM
\longrightarrow j_*(\cN_W(-(0,\infty)))
\cong\cN(-(0,\infty)),
\]
where the right-hand isomorphism follows from the
projection formula.

At the generic points of $\{\infty\}\times X$,
this morphism agrees with the original map on $L$.
Its stability and contraction properties therefore
follow from \eqref{eq:contraction-lattice}.
The remaining conditions follow from the construction,
so the triple $(\cN,\cM,\tau_{\cM})$ defines a
$C\times X$-shtuka model for $\underline M$.
\end{proof}

\begin{Remark}
It is noteworthy that the converse of Proposition \ref{prop:existence-CxX-shtuka-models} holds: if $\underline{M}$ admits a $C\times X$-shtuka model, then all the weights of $\underline{M}$ are negative. Indeed, one shows that the existence of such a model implies that the closure $T$ of $\cM_{(\infty,0)}$ is a stable $\cA_\infty(F)$-lattice in $\cB_\infty(M)$ such that
\[
\tau_M^h\bigl(\tau^{h*}T\bigr)\subset \fm_{\infty}T
\]
for some $h>0$. The existence of such a lattice occurs if and only if the slopes of the $\infty$-isocrystal attached to $\underline{M}$ are positive, \ie if and only if $\underline{M}$ has negative weights.
\end{Remark}

Let $\underline{M}$ be an $A$-motive with negative weights. Fix a $C\times X$-shtuka model $\underline{\cM}$ of $\underline{M}$, and let $i:\Spec \cA_{\infty}(\cO_F)\to C\times X$ be the canonical morphism. We set
\begin{align*}
L_{\cO_F}&:=(i^*\cM)(\Spf \cA_{\infty}(\cO_F))=(i^*\cN)(\Spf \cA_{\infty}(\cO_F)),\\
L&:=L_{\cO_F}\otimes_{\cA_{\infty}(\cO_F)}\cA_{\infty}(F),
\end{align*}
where the first equality uses that $\Delta$ is not supported on $\Spf\cA_{\infty}(\cO_F)$. The map $\tau_{\cM}$ induces an $\cO_{\infty}$-linear endomorphism of $L$ (resp.\ of $L_{\cO_F}$). The next lemma records an additional pleasant feature of $C\times X$-shtuka models, which will be used later in Subsection \ref{sec:proof}.

\begin{Lemma}\label{lem:id-tau-isomorphism-on-TO}
The morphism $\id-\tau_{\cM}$ induces an $\cO_{\infty}$-linear automorphism of $L$ and of $L_{\cO_F}$.
\end{Lemma}

\begin{proof}
It suffices to prove the claim for $L_{\cO_F}$, as tensoring with $\cA_{\infty}(F)$ yields the corresponding statement for $L$. Because $\cM$ is locally free of finite rank, $L_{\cO_F}$ is a finite projective $\cA_{\infty}(\cO_F)$\nobreakdash-module. In particular, for all integer $n\geq 0$ we have
\begin{equation}\label{eq:intersection-for-L}
L_{\cO_F}\cap \fm_{\infty}^n L = \fm_{\infty}^n L_{\cO_F}.
\end{equation}
By Lemma \ref{lem:stable-lattice-if-negative-weights}, there exist integers $h,d>0$ and an $\cA_{\infty}(F)$-lattice $T$ such that
\[
\tau_{\cM}^h\bigl(\tau^{h*}T\bigr)\subset \fm_\infty^d\,T.
\]
There also exists some integer $a\geq 0$ for which $\fm_{\infty}^{a} T\subset L \subset \fm_{\infty}^{-a} T$. Given an arbitrary integer $m\geq 0$, we deduce that 
\[
\tau_M^m(\tau^{m*}L) \subset \fm_{\infty}^{-a}\tau_M^m(\tau^{m*}T) \subset \fm_{\infty}^{d \lfloor\frac{m}{h}\rfloor-a} T \subset \fm_{\infty}^{d \lfloor\frac{m}{h}\rfloor-2a} L.
\] 
We deduce from \eqref{eq:intersection-for-L} that the same inclusion holds for $L_{\cO_F}$ in place of $L$:
\[
\text{for~all~}m\geq 0:\quad \tau_M^m(\tau^{m*}L_{\cO_F})\subset  \fm_{\infty}^{d \lfloor\frac{m}{h}\rfloor-2a} L_{\cO_F}.
\]
\smallskip\noindent
\emph{Injectivity.}
Let $x\in \ker(\id-\tau_{\cM}|_{L_{\cO_F}})$, so $x=\tau_{\cM}( \tau^*x)$. Iterating gives, for all $n\geq 1$,
\[
x=\tau_{\cM}^{nh}\bigl(\tau^{nh*}x\bigr)\in \fm_\infty^{nd-2a}\,L_{\cO_F}.
\]
Since $d>0$ and $\bigcap_{n\ge 1}\fm_\infty^{nd}L_{\cO_F}=\{0\}$, we get $x=0$.

\smallskip\noindent
\emph{Surjectivity.}
Let $y\in L_{\cO_F}$. For $n\geq 0$ and $u\in\{0,1,\dots,h-1\}$ we have
\[
\tau_{\cM}^{nh+u}\bigl(\tau^{(nh+u)*}y\bigr)\in \fm_\infty^{nd-2a}\,L_{\cO_F}.
\]
Hence the series
\[
f:=\sum_{t=0}^{\infty}\tau_{\cM}^t\bigl(\tau^{t*}y\bigr)
=\sum_{n=0}^{\infty}\left(\sum_{u=0}^{h-1}\tau_{\cM}^{nh+u}\bigl(\tau^{(nh+u)*}y\bigr)\right)
\]
converges in $L_{\cO_F}$, and it satisfies $f-\tau_{\cM}(\tau^*f)=y$. Thus $\id-\tau_{\cM}$ is surjective on $L_{\cO_F}$.
\end{proof}
 
\subsection{Shtuka models and extensions of Hodge-Pink structures}\label{subsec:shtuka-models-MHPS}
Let $\underline{M}$ be a rigid analytically trivial $A$-motive over $F$ whose weights are all negative. Let $(\cN,\cM,\tau_{\cM})$ be a $C\times X$-shtuka model for $\underline{M}$, whose existence is ensured by Proposition \ref{prop:existence-CxX-shtuka-models}. Let $\iota:\cM\to \cN$ be the inclusion of sheaves.

Recall that the ring $\cA_{\infty}(\cO_{\FI})$ is defined as the $\fm_{\infty}\otimes 1$-adic completion of $\OI\otimes \cO_{\FI}$. A trivial yet important remark is that $\cA_{\infty}(\cO_{\FI})$ is also $1\otimes \fm_{\FI}$-adically complete (hence $(\fm_{\infty}\otimes 1+1\otimes \fm_{\FI})$-adically complete). This follows from the identification \eqref{eq:identification-Binfty} which implies
\begin{equation}\label{eq:description-of-Ainfty}
\cA_{\infty}(\cO_{\FI})
=(\bF_{\infty}\otimes \cO_{\FI})[\![\pi_{\infty}\otimes 1]\!]
=(\bF_{\infty}\otimes k_{\FI})[\![1\otimes \varpi_{\infty},\pi_{\infty}\otimes 1]\!]
\end{equation}
where $k_{\FI}$ denotes the residue field of $\FI$ (it is a finite extension of $\bF_{\infty}$) and $\varpi_{\infty}\in \FI$ a uniformizer.

There is a map of ringed spaces
\begin{equation}\label{eq:morphism-of-ringed-space}
\Spf\cA_{\infty}(\cO_{\FI})\longrightarrow C\times X
\end{equation}
which, topologically, picks the formal neighborhood of $\{\infty\}\times \{\infty\}$ in $C\times X$. Denote by $\hat{\cN}$ and $\hat{\cM}$ the pullbacks of $\cN$ and $\cM$, respectively, along \eqref{eq:morphism-of-ringed-space}, and set
\[
\hat{\cN}_{\infty}:=\hat{\cN}\bigl(\Spf\cA_{\infty}(\cO_{\FI})\bigr), \qquad
\hat{\cM}_{\infty}:=\hat{\cM}\bigl(\Spf\cA_{\infty}(\cO_{\FI})\bigr),
\]
viewed as finitely generated $\cA_{\infty}(\cO_{\FI})$-modules.

The aim of this subsection is to prove that there is an exact sequence of $\KI$-vector spaces (Corollary \ref{cor:shtuka-model-to-hodge-pink-additive})
\[
0\longrightarrow \Betti{\underline{M}}^+\otimes_A \KI
\longrightarrow \KI\otimes_{\OI}\frac{\hat{\cN}_{\infty}}{(\iota-\tau_{\cM})(\hat{\cM}_{\infty})}
\longrightarrow \Ext^{1,\mathrm{ha}}_{\infty}(\mathbbm{1}^+,\sH^+(\underline{M}))
\longrightarrow 0.
\]
This exact sequence is a fundamental ingredient in the proof of Theorem \ref{thm:rank-dim}. A surprising feature is that it does not depend on the choice of shtuka model. The reader will have no trouble noticing how much this subsection relies on ideas from V.\ Lafforgue \cite[\S 4]{lafforgue}.

We start with a local comparison.
\begin{Proposition}\label{prop:locally-at-inftyxinfty}
There is an isomorphism of $\KI$-vector spaces
\[
\KI\otimes_{\OI}\frac{\hat{\cN}_{\infty}}{(\iota-\tau_{\cM})(\hat{\cM}_{\infty})}
\cong
\KI\otimes_{\OI}\frac{\hat{\cN}_{\infty}}{\iota(\hat{\cM}_{\infty})}.
\]
\end{Proposition}

We first show that a version of Proposition \ref{prop:locally-at-inftyxinfty} holds integrally  after rescaling $\hat{\cM}_{\infty}$ and $\hat{\cN}_{\infty}$ by some power of $(1\otimes \varpi_{\infty})$.
\begin{Lemma}\label{lem:iota-tau-and-iota-for-t-large}
For any large enough integer $t$, there exists an $\cO_{\infty}$-linear isomorphism 
\begin{equation}\label{eq:proof-first-step}
\frac{(1\otimes\varpi_{\infty})^t\hat{\cN}_{\infty}}{(\iota-\tau_{\cM})((1\otimes\varpi_{\infty})^t\hat{\cM}_{\infty})} \cong \frac{(1\otimes\varpi_{\infty})^t\hat{\cN}_{\infty}}{\iota((1\otimes\varpi_{\infty})^t\hat{\cM}_{\infty})}.
\end{equation}
\end{Lemma}
\begin{proof}
The core of the proof is to exhibit an isomorphism modulo $\fm_{\infty}$. As modules over the ring $\cA_{\infty}(\mathcal{O}_{F_\infty})/(\fm_{\infty}\otimes 1)=\bF_{\infty}\otimes \mathcal{O}_{F_\infty}$, set
\[
\overline{\cM} := \hat{\cM}_{\infty}/(\fm_{\infty}\otimes 1)\hat{\cM}_{\infty} \quad \text{and} \quad \overline{\cN} := \hat{\cN}_{\infty}/(\fm_{\infty}\otimes 1)\hat{\cN}_{\infty}
\] 
together with $\overline{\iota}:\overline{\cM}\to \overline{\cN}$ and $\overline{\tau}_{\cM}:\overline{\cM}\to \overline{\cN}$ (which is only semilinear). Since the cokernel of the inclusion $\iota:\cM\to \cN$ is supported at $\Delta$, there exists $s\geq 0$ such that $(\pi_{\infty}\otimes 1-1\otimes \pi_{\infty})^s v\in \iota(\hat{\cM}_{\infty})$ for all $v\in \hat{\cN}_{\infty}$. In particular, up to increasing $s$ if necessary, we see that there exists an integer $r$ such that
\[
(1\otimes \varpi_{\infty})^r \overline{\cN}\subset \overline{\iota}(\overline{\cM}).
\]

We claim that $\overline{\iota}$ remains injective; equivalently that $\hat{\cN}_{\infty}/\iota(\hat{\cM}_{\infty})$ is $(\pi_{\infty}\otimes 1)$-torsion free. To see this, let $n\in \hat{\cN}_{\infty}$ be such that $(\pi_{\infty}\otimes 1)n\in \iota(\hat{\cM}_{\infty})$. Let also $m\in \hat{\cM}_{\infty}$ be such that $\iota(m):=(\pi_{\infty}\otimes 1-1\otimes \pi_{\infty})^s n$. Then
\begin{equation}\label{eq:relation-among-m-and-n}
(\pi_{\infty}\otimes 1)\iota(m) = (\pi_{\infty}\otimes 1-1\otimes \pi_{\infty})^s (\pi_{\infty}\otimes 1)n \in (\pi_{\infty}\otimes 1-1\otimes \pi_{\infty})^s \iota(\hat{\cM}_{\infty}).
\end{equation}
Because $\hat{\cM}_{\infty}$ is a flat module over $\cA_{\infty}(\cO_{F_{\infty}})$ and $(\pi_{\infty}\otimes 1,\pi_{\infty}\otimes 1-1\otimes \pi_{\infty})$ forms a regular sequence in this ring, the quotient module $\hat{\cM}_{\infty}/(\pi_{\infty}\otimes 1-1\otimes \pi_{\infty})^s \hat{\cM}_{\infty}$ is $(\pi_{\infty}\otimes 1)$-torsion free. From \eqref{eq:relation-among-m-and-n}, we deduce that there exists $m_0\in \hat{\cM}_{\infty}$ such that $m=(\pi_{\infty}\otimes 1-1\otimes \pi_{\infty})^s m_0$, and thus that $n=\iota(m_0)\in \iota(\hat{\cM}_{\infty})$. This implies that $\hat{\cN}_{\infty}/\iota(\hat{\cM}_{\infty})$ is $(\pi_{\infty}\otimes 1)$-torsion free, hence that $\overline{\iota}$ is injective. 

Next, choose $t\geq 0$ such that $(q-1)t\geq r$. From $\tau_{\cM}(\cM)\subset \cN(-(0,\infty))$, we have
\begin{align*}
\overline{\tau}_{\cM}((1\otimes \varpi_{\infty})^t \overline{\cM}) &\subseteq (1\otimes \varpi_{\infty})^{qt+1} \overline{\cN} \\
&\subseteq (1\otimes \varpi_{\infty})^{t+r+1} \overline{\cN} \\
&\subseteq (1\otimes \varpi_{\infty})^{t+1}\overline{\iota}(\overline{\cM}).
\end{align*}
Since $\overline{\iota}$ is injective by the above claim, $\overline{\iota}^{-1}\overline{\tau}_{\cM}$ induces a topologically nilpotent endomorphism of $(1\otimes\varpi_{\infty})^t\overline{\cM}$. This implies that $\id-\overline{\iota}^{-1}\overline{\tau}_{\cM}$ is an automorphism of $(1\otimes\varpi_{\infty})^t\overline{\cM}$, and hence that
\begin{equation}\label{eq:equality-modulo-pi}
(\overline{\iota}-\overline{\tau}_{\cM})((1\otimes\varpi_{\infty})^t\overline{\cM})=\overline{\iota}((1\otimes\varpi_{\infty})^t\overline{\cM})
\end{equation}
as submodules of $(1\otimes\varpi_{\infty})^t\overline{\cN}$.

Now consider the quotient of $(1\otimes\varpi_{\infty})^t\overline{\cN}$ by the common submodule \eqref{eq:equality-modulo-pi}. This quotient is finite-dimensional over $\bF_{\infty}$, and we let $(n_1,\ldots,n_d)$ denote a family in $(1\otimes\varpi_{\infty})^t\hat{\cN}_{\infty}$ lifting a basis of it. The $\cO_{\infty}$-linear maps
\begin{align*}
f_\iota:(1\otimes \varpi_{\infty})^t \hat{\cM}_{\infty} \oplus \cO_{\infty}^{d} &\longrightarrow (1\otimes \varpi_{\infty})^t \hat{\cN}_{\infty}, \\
   (m,(c_j))&\longmapsto \iota(m)+\sum_j c_jn_j,
\end{align*}
and
\begin{align*}
f_{\iota-\tau}:(1\otimes \varpi_{\infty})^t \hat{\cM}_{\infty} \oplus \cO_{\infty}^{d}
   & \longrightarrow (1\otimes \varpi_{\infty})^t \hat{\cN}_{\infty}, \\
(m,(c_j)) &\longmapsto (\iota-\tau_{\cM})(m)+\sum_j c_jn_j
\end{align*}
are bijective modulo $\fm_{\infty}$. In addition, their source and target are $\fm_{\infty}$-adically complete, separated and torsion free. Consequently, both $f_{\iota}$ and $f_{\iota-\tau}$ are isomorphisms. This shows \eqref{eq:proof-first-step} as both sides are free of rank $d$ over $\cO_{\infty}$. 
\end{proof}

Two further lemmas complete the proof of Proposition \ref{prop:locally-at-inftyxinfty}. 
\begin{Lemma}\label{lem:preparation-un-lafforgue}
Let $t$ be a positive integer. Then $\iota-\tau_{\cM}$ and $\iota$ respectively induce isomorphisms of $\KI$-vector spaces
\[
\KI\otimes_{\OI}\left(\frac{\hat{\cM}_{\infty}}{(1\otimes\pi_{\infty})^t \hat{\cM}_{\infty}}\right)
\xrightarrow{\ \iota-\tau_{\cM}\ }
\KI\otimes_{\OI}\left(\frac{\hat{\cN}_{\infty}}{(1\otimes\pi_{\infty})^t \hat{\cN}_{\infty}}\right),
\]
\[
\KI\otimes_{\OI}\left(\frac{\hat{\cM}_{\infty}}{(1\otimes\pi_{\infty})^t \hat{\cM}_{\infty}}\right)
\xrightarrow{\ \iota\ }
\KI\otimes_{\OI}\left(\frac{\hat{\cN}_{\infty}}{(1\otimes\pi_{\infty})^t \hat{\cN}_{\infty}}\right).
\]
\end{Lemma}
\begin{proof}
Recall that there exists some integer $s\geq 0$ such that $(\pi_{\infty}\otimes 1-1\otimes \pi_{\infty})^s v\in \hat{\cM}_{\infty}$ for all $v\in \hat{\cN}_{\infty}$. We let
\[
\iota':\hat{\cN}_{\infty}\longrightarrow \hat{\cM}_{\infty}
\]
be multiplication by $(\pi_{\infty}\otimes 1-1\otimes \pi_{\infty})^s$. The following holds:
\begin{itemize}
\item $\iota'$ is an injective $\cA_{\infty}(\cO_{\FI})$-linear morphism,
\item $\iota'\circ \iota:\hat{\cM}_{\infty}\to \hat{\cM}_{\infty}$ is multiplication by $(\pi_{\infty}\otimes 1-1\otimes \pi_{\infty})^s$,
\item $\iota\circ \iota':\hat{\cN}_{\infty}\to \hat{\cN}_{\infty}$ is multiplication by $(\pi_{\infty}\otimes 1-1\otimes \pi_{\infty})^s$.
\end{itemize}
The multiplication by
\[
\left(\sum_{k=0}^{t-1}{\pi_{\infty}^{-(k+1)}\otimes \pi_{\infty}^k}\right)^s
=\left(\frac{1-\pi_{\infty}^{-t}\otimes \pi^{t}_{\infty}}{\pi_{\infty}\otimes 1-1\otimes \pi_{\infty}}\right)^s
\equiv (\pi_{\infty}\otimes 1-1\otimes \pi_{\infty})^{-s}
\]
on $\KI\otimes_{\OI}\bigl(\hat{\cM}_{\infty}/(1\otimes \pi_{\infty})^t \hat{\cM}_{\infty}\bigr)$ defines an inverse of $\iota'\iota$. The same argument shows that $\iota\iota'$ is an automorphism of $\KI\otimes_{\OI}\bigl(\hat{\cN}_{\infty}/(1\otimes \pi_{\infty})^t \hat{\cN}_{\infty}\bigr)$.

On the other hand, we have $\tau_{\cM}(\tau^* \hat{\cM}_{\infty})\subseteq (1\otimes \varpi_{\infty})\hat{\cN}_{\infty}$ by condition \ref{item:CxC-shtuka-model-zero} in Definition \ref{def:shtuka-model CC}, hence
$(\iota'\tau_{\cM})(\tau^*\hat{\cM}_{\infty})\subset (1\otimes \varpi_{\infty})\hat{\cM}_{\infty}$. Thus $\iota'\tau_{\cM}$ is nilpotent on
$\hat{\cM}_{\infty}/(1\otimes \pi_{\infty})^t \hat{\cM}_{\infty}$, and so is $(\iota'\iota)^{-1}(\iota'\tau_{\cM})$. In particular,
\[
\iota'(\iota-\tau_{\cM})
=(\iota'\iota)\Bigl(\id-(\iota'\iota)^{-1}(\iota'\tau_{\cM})\Bigr)
\]
is an isomorphism. It follows that $\iota-\tau_{\cM}$ is injective and $\iota'$ is surjective. Since $\iota\iota'$ is invertible, $\iota'$ is injective. We deduce that $\iota-\tau_{\cM}$, $\iota'$ and thus $\iota$ are isomorphisms.
\end{proof}

\begin{Lemma}
Let $t$ be a non-negative integer. Then the canonical maps
\begin{equation}\label{eq:quotient-iota-tau}
\KI\otimes_{\OI}\frac{(1\otimes\pi_{\infty})^t \hat{\cN}_{\infty}}{(\iota-\tau_{\cM})((1\otimes\pi_{\infty})^t \hat{\cM}_{\infty})}
\longrightarrow
\KI\otimes_{\OI}\frac{\hat{\cN}_{\infty}}{(\iota-\tau_{\cM})(\hat{\cM}_{\infty})},
\end{equation}
\begin{equation}\label{eq:quotient-iota}
\KI\otimes_{\OI}\frac{(1\otimes\pi_{\infty})^t \hat{\cN}_{\infty}}{\iota((1\otimes\pi_{\infty})^t \hat{\cM}_{\infty})}
\longrightarrow
\KI\otimes_{\OI}\frac{\hat{\cN}_{\infty}}{\iota(\hat{\cM}_{\infty})},
\end{equation}
are isomorphisms of $\KI$-vector spaces.
\end{Lemma}

\begin{proof}
In the category of $\OI$-modules, we have a diagram with exact rows and commutative squares:
\begin{equation}\label{diagram:from-t-to-general-lafforgue}
\begin{tikzcd}
0 \arrow[r] &
(1\otimes\pi_{\infty})^t\hat{\cM}_{\infty}
\arrow[r]\arrow[d,"\iota-\tau_{\cM}"] &
\hat{\cM}_{\infty}
\arrow[r]\arrow[d,"\iota-\tau_{\cM}"] &
\hat{\cM}_{\infty}/(1\otimes\pi_{\infty})^t\hat{\cM}_{\infty}
\arrow[r]\arrow[d,"\iota-\tau_{\cM}"] &
0 \\
0 \arrow[r] &
(1\otimes\pi_{\infty})^t\hat{\cN}_{\infty}
\arrow[r] &
\hat{\cN}_{\infty}
\arrow[r] &
\hat{\cN}_{\infty}/(1\otimes\pi_{\infty})^t\hat{\cN}_{\infty}
\arrow[r] &
0.
\end{tikzcd}
\end{equation}
By Lemma \ref{lem:preparation-un-lafforgue}, the third vertical arrow becomes an isomorphism after tensoring with $\KI$ over $\OI$. The first isomorphism then follows from the snake lemma. The second one follows from the same argument, with $\iota$ in place of $\iota-\tau_{\cM}$.
\end{proof}

\begin{proof}[Proof of Proposition \ref{prop:locally-at-inftyxinfty}]
Choose $t$ large enough so that Lemma \ref{lem:iota-tau-and-iota-for-t-large} applies. The desired isomorphism is the composition
\[
\begin{tikzcd}
\KI\otimes_{\OI}\displaystyle\frac{\hat{\cN}_{\infty}}{(\iota-\tau_{\cM})(\hat{\cM}_{\infty})} \arrow[r,"\eqref{eq:quotient-iota-tau}","\sim"']\arrow[d,dashed] & \KI\otimes_{\OI}\displaystyle\frac{(1\otimes \pi_{\infty})^t\hat{\cN}_{\infty}}{(\iota-\tau_{\cM})((1\otimes \pi_{\infty})^t\hat{\cM}_{\infty})}\arrow[d,"\text{Lemma}~\ref{lem:iota-tau-and-iota-for-t-large}","\wr"'] \\
\KI\otimes_{\OI}\displaystyle\frac{\hat{\cN}_{\infty}}{\iota(\hat{\cM}_{\infty})} & \KI\otimes_{\OI}\displaystyle\frac{(1\otimes \pi_{\infty})^t\hat{\cN}_{\infty}}{\iota((1\otimes \pi_{\infty})^t\hat{\cM}_{\infty})}\arrow[l,"\sim","\eqref{eq:quotient-iota}"']
\end{tikzcd} \nonumber
\]
This proves the claim.
\end{proof}

Recall that $\FI[\![\fj]\!]$ denotes the topological ring $(A\otimes \FI)^{\wedge}_{\fj}$ equipped with the product topology. We proved in Proposition \ref{prop:nu-extends-by-continuity} that the morphism $\nu:A\to \FI[\![\fj]\!]$, $a\mapsto a\otimes 1$, extends by continuity to a unique continuous map $\KI\to \FI[\![\fj]\!]$, also denoted by $\nu$. In particular, we get a continuous map $\nu\otimes \id:\KI\otimes \KI \to \FI[\![\fj]\!]$.

\begin{Lemma}\label{lem:B-completion-j-completion}
Let $\fj_{\infty}$ denote the kernel of the multiplication map $\cA_{\infty}(\cO_{\FI})\to \cO_{\FI}$. Then the map $\nu\otimes \id$ induces an isomorphism
\[
\cB_{\infty}(\cO_{\FI})^{\wedge}_{\fj_{\infty}} \stackrel{\sim}{\longrightarrow} \FI[\![\fj]\!].
\]
\end{Lemma}

\begin{proof}
Set $u:=\pi_{\infty}\otimes 1-1\otimes \pi_{\infty}$ in $\OI\otimes \cO_{\FI}$. By restriction, $\nu\otimes \id$ induces a map $\OI\otimes \cO_{\FI}\to \FI[\![\fj]\!]$, and the image of $u$ belongs to $\fj\setminus \fj^2$. This map is continuous for the $(\fm_{\infty}\otimes 1+1\otimes \fm_{\FI}+\fj_{\infty})$-adic topology on $\OI\otimes \cO_{\FI}$, and since the target is complete, it factors uniquely through its completion.

A quick look at the isomorphism $\OI\cong \bF_{\infty}[\![\pi_{\infty}]\!]$ shows that, as nontopological rings, the $(\fm_{\infty}\otimes 1+1\otimes \fm_{\FI})$-adic completion of $\OI\otimes \cO_{\FI}$ coincides with its $(\fm_{\infty}\otimes 1)$-adic completion, namely $\cA_{\infty}(\cO_{\FI})$. Therefore, $\nu\otimes \id$ extends uniquely to a map
\[
\cA_{\infty}(\cO_{\FI})^{\wedge}_{\fj_{\infty}} \longrightarrow \FI[\![\fj]\!],
\quad\text{hence to}\quad
\cB_{\infty}(\cO_{\FI})^{\wedge}_{\fj_{\infty}} \longrightarrow \FI[\![\fj]\!].
\]
It remains to show that the latter morphism is an isomorphism. Since $\FI[\![\fj]\!]$ is a DVR with residue field $\FI$ and uniformizer the image of $u$, it suffices to show that $\cB_{\infty}(\cO_{\FI})^{\wedge}_{\fj_{\infty}}$ is a DVR with residue field $\FI$ and uniformizer $u$.

Observe first that\footnote{We slightly abuse notation by also denoting by $\fj_{\infty}$ the kernel of the multiplication map $\OI\otimes_{\bF_{\infty}}\cO_{\FI}\to \cO_{\FI}$.}
\[
(\OI\otimes \cO_{\FI})^{\wedge}_{\fj_{\infty}}\cong (\OI\otimes_{\bF_{\infty}} \cO_{\FI})^{\wedge}_{\fj_{\infty}}.
\]
Indeed, if $f\in \bF_{\infty}$, then for all $n\geq 0$ we have
\[
f\otimes 1-1\otimes f=(f\otimes 1-1\otimes f)^{q^{d_\infty n}}\in \fj_{\infty}^{\,q^{d_\infty n}},
\]
where $d_\infty=[\bF_{\infty}:\bF]$, hence $f\otimes 1-1\otimes f$ vanishes in the $\fj_{\infty}$-adic completion.

Consequently,
\begin{align*}
(\OI\otimes \cO_{\FI})^{\wedge}_{\fm_{\infty}\otimes 1+1\otimes \fm_{\FI}+\fj_{\infty}}
&\cong (\OI\otimes_{\bF_{\infty}} \cO_{\FI})^{\wedge}_{\fm_{\infty}\otimes 1+1\otimes \fm_{\FI}+\fj_{\infty}} \\
&\cong \cO_{\FI}[\![\pi_{\infty}\otimes 1]\!]^{\wedge}_{\fj_{\infty}} \\
&\cong \cO_{\FI}[\![u]\!],
\end{align*}
where in the last isomorphism we used that $\fj_{\infty}$ is principal, generated by $u$, and that $\cO_{\FI}[\![\pi_{\infty}\otimes 1]\!]$ is a subring of $\cO_{\FI}[\![u]\!]$. Inverting $\pi_{\infty}\otimes 1$ (which amounts to inverting $1\otimes \pi_{\infty}$ as they coincide modulo $u$) yields
\[
\cB_{\infty}(\cO_{\FI})^{\wedge}_{\fj_{\infty}}\cong \FI[\![u]\!],
\]
which concludes the proof.
\end{proof}

We are now in a position to prove the main result of this subsection.

\begin{Theorem}\label{thm:shtuka-model-to-hodge-pink-additive}
Let $(\cN,\cM,\tau_{\cM})$ be a $C\times X$-shtuka model for $\underline{M}$. Then there is an isomorphism of $\KI$-vector spaces
\[
\KI\otimes_{\OI}\frac{\hat{\cN}_{\infty}}{(\iota-\tau_{\cM})(\hat{\cM}_{\infty})}
\stackrel{\sim}{\longrightarrow}
\frac{(M+\tau_M(\tau^* M))\otimes_{A\otimes F}\FI[\![\fj]\!]}{M\otimes_{A\otimes F}\FI[\![\fj]\!]},
\]
where the $\KI$-vector space structure on the right-hand side is given through $\nu$.
\end{Theorem}

\begin{proof}
Let $\Spec B\subset C$ (resp. $\Spec B'\subset X$) be an affine open subscheme containing the closed point $\infty$. Let $M_B$ (resp. $N_B$) be the value of $\cM$ (resp. $\cN$) on the Zariski open subscheme $\Spec B\otimes B'$ of $C\times X$. By the sheaf property and since the cokernel of $\iota:\cM\to \cN$ is supported on $\Delta$, we have
\[
\frac{\hat{\cN}_{\infty}}{\hat{\cM}_{\infty}}
\cong
\frac{N_B}{M_B}\otimes_{B\otimes B'}\cA_{\infty}(\cO_{\FI})
\cong
\frac{N_B}{M_B}\otimes_{B\otimes B'}\cA_{\infty}(\cO_{\FI})^{\wedge}_{\fj_{\infty}}.
\]
Extending scalars along $\OI\to \KI$ and invoking Lemma \ref{lem:B-completion-j-completion}, we obtain
\[
\KI\otimes_{\OI}\frac{\hat{\cN}_{\infty}}{\hat{\cM}_{\infty}}
\cong
\frac{N_B}{M_B}\otimes_{B\otimes B'}\cB_{\infty}(\cO_{\FI})^{\wedge}_{\fj_{\infty}}
\cong
\frac{N_B}{M_B}\otimes_{B\otimes B'}\FI[\![\fj]\!]
\cong
\frac{N}{M}\otimes_{A\otimes F}\FI[\![\fj]\!].
\]
Precomposing with the isomorphism of Proposition \ref{prop:locally-at-inftyxinfty} yields the desired isomorphism.
\end{proof}

As announced, we have the following consequence.

\begin{Corollary}\label{cor:shtuka-model-to-hodge-pink-additive}
There is an exact sequence of $\KI$-vector spaces
\[
0\longrightarrow \Betti{\underline{M}}^+\otimes_A\KI
\longrightarrow
\KI\otimes_{\OI}\frac{\hat{\cN}_{\infty}}{(\iota-\tau_{\cM})(\hat{\cM}_{\infty})}
\longrightarrow
\Ext^{1,\mathrm{ha}}_{\infty}(\mathbbm{1}^+,\sH^+(\underline{M}))
\longrightarrow 0.
\]
\end{Corollary}

\begin{proof}
This is a consequence of Theorem \ref{thm:shtuka-model-to-hodge-pink-additive} combined with Theorem \ref{thm:extensions-of-MHPS-arising-from-motives}.
\end{proof}

\section{Proof of the main theorems}\label{sec:reg-fin-thm}
Let $\underline{M}$ be a rigid analytically trivial $A$-motive over $F$ and $\underline{H}^+:=\sH^+(\underline{M})$ be the Hodge-Pink structure with infinite Frobenius associated to $\underline{M}$ in Section \ref{sec:rat-mix-A-mot}. This section proves Theorems A and B of the introduction through the following three statements. Theorem \ref{thm:finiteness-motcoh} restates Theorem~A. Theorem~B follows from Theorems \ref{thm:rank-dim} and \ref{thm:regulator-injective-intro}.

\begin{Theorem}\label{thm:finiteness-motcoh}
The $A$-modules $\Ext^{1,\operatorname{reg}}_{\cO_F,\infty}(\mathbbm{1},\underline{M})$ and $\operatorname{Cl}(\underline{M})$ are finitely generated. If the weights of $\underline{M}$ are all negative, then $\operatorname{Cl}(\underline{M})$ is finite. 
\end{Theorem}

\begin{Theorem}\label{thm:rank-dim}
Assume that the weights of $\underline{M}$ are all negative. Then the rank of $\Ext^{1,\operatorname{reg}}_{\cO_F,\infty}(\mathbbm{1},\underline{M})$ as an $A$-module equals the dimension of the $\KI$-vector space $\Ext^{1,\operatorname{ha}}_{\infty}(\mathbbm{1}^+,\underline{H}^+)$.
\end{Theorem}

\begin{Theorem}\label{thm:regulator-injective-intro}
Assume that the weights of $\underline{M}$ are all negative. Then the $\KI$-linearized regulator $\Reg(\underline{M})_{\KI}$ is injective.
\end{Theorem}

Theorems \ref{thm:finiteness-motcoh} and \ref{thm:rank-dim} will follow from computing the coherent cohomology of shtuka models on certain coverings of $C\times X$. Theorem \ref{thm:regulator-injective-intro} will be deduced from a careful rigid analytic study of the $\tau$-difference equations that arise from computing extension modules.

\subsection{Cohomological computations}\label{sec:cohomological-computations}
In this subsection, we recall general facts on coherent cohomology of schemes and how to compute it on coverings that are not necessarily Zariski. This will subsequently be applied to shtuka models in the next subsection to complete the proofs of Theorems \ref{thm:finiteness-motcoh} and \ref{thm:rank-dim}. Throughout this subsection, all schemes and morphisms are assumed to be quasi-compact and quasi-separated.

Let $A$ be a commutative ring\footnote{For the sake of this exposition, $A$ is not meant to be $\cO_C(C\setminus\{\infty\})$ here, as it generally is in this text.}. To avoid issues with functoriality of cones, we shall use the language of stable $\infty$-categories and the modern construction of the derived category of quasi-coherent sheaves on a scheme, \eg described at length in \cite[Part I]{gaitsgory-rozenblyum}. That is, given an $A$-algebra $R$, we shall denote by $\sD(R)$ the stable derived $\infty$-category of $R$-modules. More generally, given a scheme $X\xrightarrow{p} S$ over the affine base scheme $S=\Spec A$, we denote by $\sD_{\qc}(X)$ the stable derived $\infty$-category of quasi-coherent sheaves on $X$ over $S$: it can be recovered as the limit\footnote{For this to make sense, we treat $\sD$ as a functor from commutative rings to presentable stable $\infty$-categories.}
\[
\sD_{\qc}(X):= \lim_{\Spec R \to X} \sD(R)
\]
indexed over all affine $S$-schemes mapping to $X$. Under this description, an \emph{honest} quasi\nobreakdash-coherent sheaf $\cM$ amounts to the data, for all affines $\Spec R$ mapping to $X$, of an $R$-module $M_R$ that suitably glues on overlaps $\Spec R_1 \times_X \Spec R_2$; under this description, $\cM$ is seen as an object in $\sD_{\qc}(X)$ via the family $(M_R[0])_{\Spec R \to X}$ together with its gluing data. In particular, the $\infty$-category $\sD_{\qc}(X)$ comes with a canonical ring object $\cO_X$, obtained as the gluing of the $R[0]$'s, and any object of $\sD_{\qc}(X)$ is a module over $\cO_X$. 

We also recall that $\sD_{\qc}(X)$ is symmetric monoidal with respect to the derived tensor product $\otimes_{\cO_X}^L$. 

If $f:X\to Y$ is a morphism of $S$-schemes, the map on the indexing categories $\{\Spec R\to X\}\to \{\Spec R\to Y\}$ given by composition with $f$ induces the \emph{pullback functor} $f^*:\sD_{\qc}(Y)\to \sD_{\qc}(X)$. It is symmetric monoidal, and in particular satisfies $f^*\cO_Y\cong \cO_X$. Since it commutes with arbitrary colimits, it admits a right-adjoint by Lurie's theorem: this is the pushforward functor $f_*:\sD_{\qc}(X)\to \sD_{\qc}(Y)$. Given a quasi-coherent sheaf $\cF$ on $X$ we denote
\[
R\Gamma(X,\cF):=p_*\cF \in \sD(A)
\]
where $p:X\to S$ is the structure morphism. When $f:\Spec R\to \Spec Q$ is a morphism of affine $S$-schemes, the pushforward $f_*:\sD(R)\to \sD(Q)$ is simply given by restricting the module structure along the $A$-algebra map $Q\to R$. Given an object $\cF$ in $\sD_{\qc}(X)$ and a map $f:\Spec R\to X$, the module $\cF(\Spec R)\in \sD(A)$ \emph{of $R$-sections of $\cF$} is defined to be $(p_R)_*f^*\cF$ where $p_R:\Spec R\to S$ is the structure morphism; the $R$-sections of the \emph{honest} quasi-coherent sheaf $\cM$ are $M_R[0]$, as expected. 

Now let $T$ be a separated $S$-scheme. Let $U$, $V$ and $W$ be affine $S$-schemes which fit into a commutative diagram
\[
\begin{tikzcd}
U \arrow[r,"i"] & T \\
W \arrow[u]\arrow[ur,dashed,"k"] \arrow[r] & V \arrow[u,"j"']
\end{tikzcd}
\]
Given a quasi-coherent sheaf $\cF$ on $T$ (concentrated in degree zero), we denote by $S(\cF)$ the following sequence in $\sD_{\qc}(T)$:
\[
S(\cF):\quad \cF \longrightarrow i_*i^* \cF\oplus j_*j^* \cF \longrightarrow k_*k^*\cF
\]
where the morphisms are given by the adjunction unit (note that the data of $S(\cF)$ is functorial in $\cF$). The next lemma will allow us to compute the cohomology of shtuka models on coverings that are not necessarily Zariski:

\begin{Lemma}\label{lem:bhatt}
Assume that $S(\cO_T)$ is a fiber sequence. Then $S(\cF)$ is a fiber sequence for any $\cF\in \sD_{\qc}(T)$. In particular, the natural map
\begin{equation}\label{eq:computeRGamma}
R\Gamma(T,\cF)\longrightarrow \left[\cF(U)\oplus \cF(V)\longrightarrow \cF(W)\right],
\end{equation}
where the right-hand side is a complex sitting in degrees $0$ and $1$, is a quasi-isomorphism.
\end{Lemma}
\begin{proof}
Recall the projection formula: given a quasi-compact morphism of schemes $f:X\to Y$, $E$ an object in $\sD_{\qc}(Y)$ and $F$ in $\sD_{\qc}(X)$, we have
\[
f_*(f^*E\otimes^L_{\cO_X} F)\cong E\otimes^L_{\cO_Y}f_*F \quad \text{in}\quad \sD_{\qc}(Y)
\]
(see \cite[Lemma~I.III.3.2.4]{gaitsgory-rozenblyum}). If $f:X\to T$ is an affine morphism of schemes with target $T$, then it is quasi-compact and hence we obtain
\begin{equation}\label{eq:projection-formula-in-the-use}
f_*f^* \cF\cong \cF\otimes_{\cO_T}^L f_* f^* \cO_T
\end{equation}
where we plugged $E=\cF$ and $F=\cO_X=f^* \cO_Y$ in the projection formula. Applying formula \eqref{eq:projection-formula-in-the-use} to $f=i,j,k$, we obtain $S(\cF)\cong \cF\otimes^L_{\cO_T}S(\cO_T)$. By assumption, $S(\cO_T)$ is a fiber sequence, hence $S(\cF)$ is a fiber sequence. 

The second assertion follows from the former, after applying the pushforward functor $p_*=R\Gamma(T,-)$ to $S(\cF)$: because $i$ is an affine morphism, we get $(p\circ i)_*i^*\cF\cong R\Gamma(U,i^*\cF)\cong \cF(U)$, and similarly for $j$ and $k$.
\end{proof}

The language of \cite{gaitsgory-rozenblyum} is equally suited to describe categories of formal quasi-coherent sheaves: assume that $I\subset A$ is a finitely generated ideal in $A$. We consider 
\begin{equation}\label{eq:inclusion-of-complete-modules}
\sD(R)^{\wedge}_I \subset \sD(R)
\end{equation}
the full subcategory of derived $I$-adically complete objects. We may also denote this category by $\sD_{\qc}(\Spf R)$. The inclusion \eqref{eq:inclusion-of-complete-modules} admits a left-adjoint which is the derived $I$-completion functor. 

We globalize the situation for a scheme $X$ and define the category of quasi-coherent sheaves on $\mathfrak{X}:=X^{\wedge}_I$ as
\[
\sD_{\qc}(\mathfrak{X}):=\sD_{\qc}(X)^{\wedge}_{I}:= \lim_{\Spec R \subseteq X} \sD(R)^{\wedge}_{I}.
\]
It also comes with a canonical ring object $\cO_\mathfrak{X}=(\cO_X)^{\wedge}_I$, is symmetric monoidal with respect to the derived completed tensor product $\widehat{\otimes}_{\cO_{\mathfrak{X}}}^L$, and admits pushforwards and pullbacks. There is an inclusion functor $\sD_{\qc}(\mathfrak{X})\subset \sD_{\qc}(X)$ as well as a derived completion functor $\sD_{\qc}(X)\to \sD_{\qc}(\mathfrak{X})$. Given a quasi-coherent sheaf $\mathfrak{F}$ on $\mathfrak X$ we denote
\[
R\Gamma(\mathfrak{X},\mathfrak{F}):=p_*\mathfrak{F} \in \sD(A)
\]
where $p:X\to S$ is the structure morphism. 

We recall the derived version of Grothendieck's theorem on formal functions \href{https://stacks.math.columbia.edu/tag/0A0G}{[0A0G]}, applied to our situation:
\begin{Theorem}[Derived formal functions]\label{thm:GCT}
Let $f:X\to Y$ be a qcqs morphism over $S$. Then the functor $(f_*)^{\wedge}_I:\sD_{\qc}(X)^{\wedge}_I\to \sD_{\qc}(Y)^{\wedge}_I$ is right-adjoint to $(f^*)^{\wedge}_I$. In particular, if $X$ is a proper $S$-scheme, there is an equivalence
\[
R\Gamma(X,\cF)^{\wedge}_I\cong R\Gamma(\mathfrak{X},\cF^{\wedge}_I)
\]
in the derived category of complete modules $\sD(A)^{\wedge}_I$.
\end{Theorem}

In the setting of Lemma \ref{lem:bhatt}, where $\mathfrak{T}$, $\mathfrak{U}$, $\mathfrak{V}$ and $\mathfrak{W}$ denote the formal completions of the schemes $T$, $U$, $V$ and $W$, we consider the sequence in $\sD_{\qc}(\mathfrak{T})$:
\[
S(\mathfrak{F}):\quad \mathfrak{F}\longrightarrow i_*i^* \mathfrak{F}\oplus j_*j^* \mathfrak{F} \longrightarrow k_*k^* \mathfrak{F}
\]
The proof of Lemma \ref{lem:bhatt} applies verbatim to give the following formal analogue:
\begin{Lemma}\label{lem:formal-analog}
Assume that $S(\cO_{\mathfrak{T}})$ is a fiber sequence. Then $S(\mathfrak{F})$ is a fiber sequence for any $\mathfrak{F} \in \sD_{\qc}(\mathfrak{T})$. In particular, there is a quasi-isomorphism
\begin{equation}\label{eq:computeRGamma-formal}
R\Gamma(\mathfrak{T},\mathfrak{F})\longrightarrow \left[\mathfrak{F}(\mathfrak{U})\oplus \mathfrak{F}(\mathfrak{V})\longrightarrow \mathfrak{F}(\mathfrak{W})\right].
\end{equation}
\end{Lemma}

\subsection{Proof of Theorems \ref{thm:finiteness-motcoh} and \ref{thm:rank-dim}}\label{sec:proof}
Let $\underline{M}$ be a rigid analytically trivial $A$-motive over $F$. We now assemble the ingredients collected in the preceding subsections to prove Theorems \ref{thm:finiteness-motcoh} and \ref{thm:rank-dim}.  As announced, the complex $G_{\underline{M}}$ admits an interpretation in terms of the Zariski cohomology of $X$-shtuka models of $\underline{M}$ (Definition~\ref{def:shtuka-model-C}):
\begin{Proposition}\label{prop:ge-in-terms-of-cohomology}
Let $(\cN,\cM,\tau_{\cM})$ be an $X$-shtuka model for $\underline{M}$. Let $\iota$ denote the inclusion of $\cM$ in $\cN$. There is a quasi-isomorphism in $\sD(A)$:
\begin{equation}
G_{\underline{M}}\stackrel{\sim}{\longrightarrow} \cofib\left(R\Gamma(\Spec A \times X,\cM) \xrightarrow{\iota-\tau_{\cM}} R\Gamma(\Spec A \times X,\cN)\right). \nonumber
\end{equation}
\end{Proposition}

The first part of Theorem \ref{thm:finiteness-motcoh} follows from the above:
\begin{Corollary}
The $A$-modules $\Ext_{\cO_F,\infty}^{1,\operatorname{reg}}(\mathbbm{1},\underline{M})$ and $\operatorname{Cl}(\underline{M})$ are finitely generated.
\end{Corollary}
\begin{proof}
As $\Spec A\times X$ is proper over $\Spec A$, both $R\Gamma(\Spec A \times X,\cM)$ and $R\Gamma(\Spec A \times X,\cN)$ are perfect complexes of $A$-modules. By Proposition \ref{prop:ge-in-terms-of-cohomology}, so is $G_{\underline{M}}$. We conclude by Proposition~\ref{prop:motivic-exact-sequence}.
\end{proof}

\begin{proof}[Proof of Proposition \ref{prop:ge-in-terms-of-cohomology}]
The main ingredients are the cohomological preliminaries of Section~\ref{sec:cohomological-computations}. We consider the particular setting of $S=\Spec A$ and of the commutative diagram of $S$-schemes
\begin{equation}\label{eq:covering-of-AxC}
\begin{tikzcd}
\Spec \cO_{\FI}\langle A \rangle \arrow[r,"i"] & (\Spec A)\times X \\
\Spec \FI\langle A \rangle \arrow[u,"p"]\arrow[ur,dashed,"k"] \arrow[r,"q"] & \Spec A\otimes \cO_F \arrow[u,"j"']
\end{tikzcd}
\end{equation}
We want to show that $R\Gamma(\Spec A \times X,\cM)$ (and similarly for $\cN$) is computed by the two-term complex
\[
\left[M_{\cO_F}\oplus i^* \cM(\Spec \cO_{\FI}\langle A \rangle)\longrightarrow M\otimes_{A\otimes F}\FI\langle A \rangle \right]
\]
This does not follow from Zariski nor fpqc descent; while diagram \eqref{eq:covering-of-AxC} gives an fpqc cover of $(\Spec A)\times X$, fpqc-descent would require to contemplate its \v{C}ech nerve and thus have to deal with much sophisticated schemes, such as $\Spec F_{\infty}\langle A\rangle \otimes_{A\otimes \cO_F} F_{\infty}\langle A\rangle$ and its higher versions.

To deduce this claim from Lemma \ref{lem:bhatt}, we need the following:
\begin{Lemma}\label{lem:check-assumption-of-bhatt-lemma}
For $T:=(\Spec A) \times X$, the sequence $S(\cO_T):\cO_T\to i_*i^*\cO_T\oplus j_*j^*\cO_T\to k_*k^*\cO_T$ is a fiber sequence in $\sD_{\qc}(T)$.
\end{Lemma}
\begin{proof}
Write $B_0:=\cO_F$, $U=\Spec B_0$, and let $V=\Spec B_{\infty}$ be an affine open subscheme of $X$ that contains the point $\infty$, so that $(\Spec A)\times U$ and $(\Spec A)\times V$ together form a Zariski covering of $\Spec A\times X$. By Zariski descent for $\sD(-)$, the projection functor $f_i^*:\sD_{\qc}(\Spec A\times X)\to \sD(A\otimes B_i)$, for $i\in \{0,\infty\}$, are jointly conservative. In particular, it suffices to check that $f_i^*S(\cO_T)$ is a fiber sequence for $i\in \{0,\infty\}$. 

Applying (derived) base change \cite[Proposition I.III.2.2.2]{gaitsgory-rozenblyum}, at $i=0$ the sequence $f^*_0S(\cO_T)$ becomes
\[
A\otimes \cO_F \xrightarrow{x\mapsto (x,x)} \FI\langle A \rangle \oplus (A\otimes \cO_F) \xrightarrow{(x,y)\mapsto x-y} \FI\langle A \rangle.
\]
Similarly, at $i=\infty$, the sequence becomes
\[
A\otimes B_{\infty} \xrightarrow{x\mapsto (x,x)} \cO_{\FI}\langle A \rangle \oplus (A\otimes B_{0\infty}) \xrightarrow{(x,y)\mapsto x-y} \FI\langle A \rangle
\]
where $B_{0\infty}:=\cO_X(U\cap V)$ is the smallest subring of $F$ that contains both $B_0$ and $B_{\infty}$. Both sequences come from short exact sequences of $A$-modules, hence are fiber sequences.
\end{proof}

Now, let $\underline{\cM}=(\cN,\cM,\tau_{\cM})$ be an $X$-shtuka model for $\underline{M}$. We have
\begin{align*}
(j^*\cM)(\Spec A\otimes \cO_F)&=\cM(\Spec A\otimes \cO_F)=M_{\cO_F},\\
(j^*\cN)(\Spec A\otimes \cO_F)&=\cN(\Spec A\otimes \cO_F)=N_{\cO_F},\\
(k^*\cM)(\Spec \FI\langle A \rangle)&=(k^*\cN)(\Spec \FI\langle A \rangle)=M\otimes_{A\otimes F}\FI\langle A\rangle. \nonumber 
\end{align*}
By Lemma \ref{lem:check-assumption-of-bhatt-lemma} we may apply Lemma \ref{lem:bhatt} with the honest quasi-coherent sheaves $\cM$ and $\cN$. From Lemma \ref{lem:isomorphism-at-infty}, we deduce a diagram in $\sD(A)$ whose columns and first two rows are fiber sequences:
\begin{equation}
\begin{tikzcd}
R\Gamma(\Spec A\times X,\cM) \arrow[r]\arrow[d,"\iota-\tau_{\cM}"] & (i^*\cM)(\Spec \cO_{\FI}\langle A \rangle) \arrow[r]\arrow[d,"\id-\tau_{\cM}"] & \displaystyle\frac{M\otimes_{A\otimes F}F_{\infty}\langle A\rangle}{M_{\cO_F}} \arrow[d,"\id-\tau_{M}"]   \\
R\Gamma(\Spec A \times X,\cN) \arrow[r]\arrow[d] & (i^*\cN)(\Spec \cO_{\FI}\langle A \rangle) \arrow[r]\arrow[d] & \displaystyle\frac{M\otimes_{A\otimes F}F_{\infty}\langle A\rangle}{N_{\cO_F}} \arrow[d] \\
\cofib(\iota-\tau_{\cM}|\Spec A\times X) \arrow[r] & 0 \arrow[r] & G_{\underline{M}}\left[1\right]
\end{tikzcd}
\nonumber
\end{equation}
Note that because $i$ is flat, the classical pullback functor $i^*$ coincides with the derived pullback. The third row is then a fiber sequence, and the proposition follows.
\end{proof}

Theorems \ref{thm:finiteness-motcoh} (second part) and \ref{thm:rank-dim} will follow from the study of the cohomology of a $C\times X$\nobreakdash-shtuka model of $\underline{M}$ on $\Spf\cO_{\infty}\hat{\times} X$. The latter corresponds to the completion of the noetherian scheme $\Spec \OI\times X$ at the closed subscheme $\{\infty\}\times X$. The argument given here is a refinement of the one given in the proof of Proposition \ref{prop:ge-in-terms-of-cohomology} where we use $C\times X$-shtuka models instead of $X$-shtuka models. To ensure the existence of a $C\times X$-shtuka model, we now assume that all the weights of $\underline{M}$ are negative (Proposition \ref{prop:existence-CxX-shtuka-models}). 

We apply the results of Section \ref{sec:cohomological-computations} in the formal setting. Let $S=\Spec\OI$ and $I=\mathfrak{m}_{\infty}\subset \OI$.  We consider the commutative square of schemes over $S$:
\begin{equation}\label{eq:formal-covering}
\begin{tikzcd}
\Spec\OI\otimes \cO_F \arrow[r,"i"] & (\Spec\OI)\times X \\
\Spec\OI\otimes \FI \arrow[u]\arrow[ur,dashed,"k"] \arrow[r] & \Spec\OI\otimes\cO_{\FI}  \arrow[u,"j"']
\end{tikzcd}
\end{equation}
To ensure that Lemma \ref{lem:formal-analog} applies in our situation, we require:
\begin{Lemma}\label{lem:check-assumption-formal}
Let $\mathfrak{T}=\Spf \OI \hat{\times} X$ denote the formal completion of $(\Spec \OI)\times X$ at $I$. The sequence $S(\cO_{\mathfrak{T}}):\cO_{\mathfrak{T}}\to i_*i^* \cO_{\mathfrak{T}}\oplus j_*j^* \cO_{\mathfrak{T}}\to k_* k^* \cO_{\mathfrak{T}}$ is a fiber sequence in $\sD_{\qc}(\mathfrak{T})$.
\end{Lemma}
\begin{proof}
We reuse the argument and notation from the proof of Lemma \ref{lem:check-assumption-of-bhatt-lemma}: we need to prove that the sequences $f_i^*S(\cO_{\mathfrak{T}})$, for $i\in\{0,\infty\}$, are fiber sequences. This amounts to checking that the sequences
\[
\cA_{\infty}(\cO_F)\to \cA_{\infty}(\cO_F)\oplus \cA_{\infty}(\FI) \to \cA_{\infty}(\FI) \quad \text{and} \quad \cA_{\infty}(B_{\infty})\to \cA_{\infty}(B_{0\infty})\oplus \cA_{\infty}(\cO_{\FI}) \to \cA_{\infty}(\FI)
\]
are fiber sequences of complete $\OI$-modules. But this is clear, as both sequences arise as the complete base change along $\bF\to \OI$ of the short exact sequences
\[
0\to \cO_F\to \cO_F\oplus \FI \to \FI\to 0 \quad \text{and} \quad 0\to B_{\infty}\to B_{0\infty}\oplus \cO_{\FI} \to \FI\to 0.
\]
\end{proof}

We recall some notation introduced earlier:
\begin{align*}
\hat{\cN}_{\infty}&:=\hat{\cN}(\Spf \OI\hat{\otimes}\cO_{\FI}) & L &:=\hat{\cM}(\Spf \OI\hat{\otimes}F)=\hat{\cN}(\Spf \OI\hat{\otimes}F) \\
\hat{\cM}_{\infty}&:=\hat{\cM}(\Spf \OI\hat{\otimes}\cO_{\FI}) & L_{\cO_F} &:=\hat{\cM}(\Spf \OI\hat{\otimes}\cO_F)=\hat{\cN}(\Spf \OI\hat{\otimes}\cO_F).
\end{align*}
Note that $L$ is an $\cA_{\infty}(F)$-lattice in $\cB_{\infty}(M):=M\otimes_{A\otimes F}\cB_\infty(F)$ stable under $\tau_M$. \\

From Lemma \ref{lem:check-assumption-formal}, we may apply Lemma \ref{lem:formal-analog} to the ``covering'' \eqref{eq:formal-covering} on the quasi-coherent sheaves $\hat{\cM}$ and $\hat{\cN}$ of $\sD_{\qc}(\Spf \OI \hat{\times }X)$. We obtain a morphism of fiber sequences in $\sD_{\qc}(\OI)^{\wedge}_{\fm_{\infty}}$:
\begin{equation}\label{eq:diagram-completed-cohomology}
\begin{tikzcd}
R\Gamma(\Spf \cO_{\infty} \hat{\times} X,\hat{\cM}) \arrow[r]\arrow[d,"\iota-\tau_{\cM}"] & \hat{\cM}_{\infty}\arrow[r]\arrow[d,"\iota-\tau_{\cM}"] & \displaystyle\frac{L\otimes_{\cA_{\infty}(F)}\cA_{\infty}(F_{\infty})}{L_{\cO_F}} \arrow[d,"\id-\tau_{M}"] \\
R\Gamma(\Spf \cO_{\infty} \hat{\times} X,\hat{\cN}) \arrow[r] & \hat{\cN}_{\infty} \arrow[r] & \displaystyle\frac{L\otimes_{\cA_{\infty}(F)}\cA_{\infty}(F_{\infty})}{L_{\cO_F}} 
\end{tikzcd}
\end{equation}
The next lemma shows that the middle vertical arrow is injective once base changed along $\OI\to \KI$:
\begin{Lemma}\label{lem:injectivity-of-V-to-W}
The morphism $\iota-\tau_{\cM}:\hat{\cM}_{\infty}\otimes_{\OI}\KI\to \hat{\cN}_{\infty}\otimes_{\OI}\KI$ is injective.
\end{Lemma}
\begin{proof}
For $t>0$, $a\geq 0$ and $x\in \hat{\cM}_{\infty}$, we have
\begin{equation}
(\iota-\tau_{\cM})((1\otimes \varpi_{\infty})^{t+a}x)\equiv (1\otimes \varpi_{\infty})^{t+a}\iota(x) \pmod{(1\otimes \varpi_{\infty})^{t+a+1}\hat{\cN}_{\infty}} \nonumber
\end{equation}
where $\varpi_{\infty}\in \FI$ is a uniformizer. In particular, the first vertical arrow in diagram \ref{diagram:from-t-to-general-lafforgue} is injective. The lemma then follows from Lemma \ref{lem:preparation-un-lafforgue} together with the snake lemma.
\end{proof}

We return to diagram \eqref{eq:diagram-completed-cohomology} extended along the flat base change $\OI\to \KI$. By Lemma \ref{lem:injectivity-of-V-to-W}, the middle vertical map becomes injective and its cofiber ($=$cokernel) is computed by Theorem \ref{thm:shtuka-model-to-hodge-pink-additive}. The third vertical arrow is a (quasi-)isomorphism by Lemma \ref{lem:id-tau-isomorphism-on-TO}. We deduce that the cofibers of the first two vertical maps are equivalent and we get 
\begin{equation}\label{eq:formal-cohomology-hodge}
\cofib\left(R\Gamma(\Spf \cO_{\infty}\hat{\times} X,\hat{\cM}) \xrightarrow{\iota-\tau_{\cM}} R\Gamma(\Spf \cO_{\infty} \hat{\times} X,\hat{\cN})\right)_{\KI}\cong \displaystyle\frac{(M+\tau_M(\tau^*M))\otimes_{A\otimes F}\FI[\![\fj]\!]}{M\otimes_{A\otimes F}\FI[\![\fj]\!]}
\end{equation}
where the subscript $_{\KI}$ denotes the base change along $\OI\to \KI$. Since $\Spec\cO_\infty\times X\to\Spec\cO_\infty$ is proper and $\cM$ is coherent, $R\Gamma(\Spec\cO_\infty\times X,\cM)$ belongs to $D^b_{\mathrm{coh}}(\OI)$. Since $\OI$ is complete for the $\mathfrak m_\infty$-adic topology, this complex is derived $\mathfrak m_\infty$-complete \href{https://stacks.math.columbia.edu/tag/0A05}{[0A05]}. Thus the derived formal functions theorem \ref{thm:GCT} provides a quasi-isomorphism
\[
R\Gamma(\Spec\cO_\infty\times X,\cM) \cong R\Gamma(\Spf\cO_\infty\widehat{\times}X,\hat{\cM}).
\]
The same holds for $\cN$. Here and below, since $\Spec\OI\times X$ is not an open subscheme of $C\times X$, the notation ``$R\Gamma(\Spec\OI \times X,-)$'' stands for $R\Gamma(\Spec\OI \times X,q^*-)$ where $q:\Spec\OI\times X\to C\times X$. This allows us to rewrite \eqref{eq:formal-cohomology-hodge} as 
\[
\cofib\left(R\Gamma(\Spec\OI\times X,\cM) \xrightarrow{\iota-\tau_{\cM}} R\Gamma(\Spec\OI \times X,\cN)\right)_{\KI} \cong \displaystyle\frac{(M+\tau_M(\tau^*M))\otimes_{A\otimes F}\FI[\![\fj]\!]}{M\otimes_{A\otimes F}\FI[\![\fj]\!]}, 
\]
from which we deduce
\[
\cofib\left(R\Gamma(\Spec A\times X,\cM) \xrightarrow{\iota-\tau_{\cM}} R\Gamma(\Spec A \times X,\cN)\right)\otimes_{A}\KI \cong \displaystyle\frac{(M+\tau_M(\tau^*M))\otimes_{A\otimes F}\FI[\![\fj]\!]}{M\otimes_{A\otimes F}\FI[\![\fj]\!]}. 
\]
Since a $C\times X$-shtuka model restricts to an $X$-shtuka model, the left-hand side now identifies with the complex $G_{\underline{M}}\otimes_A\KI$ by Proposition \ref{prop:ge-in-terms-of-cohomology}, finally yielding:
\begin{equation} 
G_{\underline{M}}\otimes_A \KI\cong \displaystyle\frac{(M+\tau_M(\tau^*M))\otimes_{A\otimes F}\FI[\![\fj]\!]}{M\otimes_{A\otimes F}\FI[\![\fj]\!]}. \nonumber
\end{equation}
Therefore, the complex $G_{\underline{M}}\otimes_A \KI$ has cohomology concentrated in degree $0$ and the second part of Theorem \ref{thm:finiteness-motcoh} follows:
\begin{Proposition}
The $A$-module $\operatorname{H}^1(G_{\underline{M}})\cong \operatorname{Cl}(\underline{M})$ is torsion, and thus finite. 
\end{Proposition}

It remains to prove Theorem \ref{thm:rank-dim}. We first introduce an auxiliary comparison isomorphism; we will use this isomorphism to establish equality of dimensions. The injectivity of the Hodge-Pink regulator itself will be proved separately in the final subsection.
\begin{Definition}\label{def:conjectural-form-regulator}
Let $\underline{M}$ be a rigid analytically trivial $A$-motive whose weights are all negative and let $\underline{\cM}:=(\cM,\cN,\tau_{\cM})$ be a $C\times X$-shtuka model for $\underline{M}$. We call \emph{the mock regulator of $\underline{\cM}$}, and denote it by $\rho(\underline{\cM})$, the isomorphism of $\KI$-vector spaces
\begin{equation}
\displaystyle\frac{\left\{\xi\in M\otimes_{A\otimes F}\FI\langle A \rangle~|~\xi-\tau_M(\tau^*\xi)\in N_{\cO_F}\right\}}{M_{\cO_F}}\otimes_{A}\KI\stackrel{\sim}{\longrightarrow} \displaystyle\frac{(M+\tau_M(\tau^*M))\otimes_{A\otimes F}\FI[\![\fj]\!]}{M\otimes_{A\otimes F}\FI[\![\fj]\!]} \nonumber
\end{equation} 
obtained by the vertical composition of the quasi-isomorphisms in $\sD(\KI)$:
\begin{equation}
\begin{tikzcd}[column sep=0.2em]
\displaystyle\frac{\left\{\xi\in M\otimes_{A\otimes F}\FI\langle A \rangle~|~\xi-\tau_M(\tau^*\xi)\in N_{\cO_F}\right\}}{M_{\cO_F}}\otimes_{A}\KI  \\
G_{\underline{M}}\otimes_A \KI \arrow[u,"\wr~\operatorname{Cl}(\underline{M})\otimes_A\KI=0"']\\
\cofib\left[R\Gamma(\Spec A \times X,\cM)\xrightarrow{\iota-\tau_{\cM}} R\Gamma(\Spec A \times X,\cN)\right]\otimes_{A}\KI \arrow[d,"\wr"] \arrow[u,"\wr~\text{Proposition}~\ref{prop:ge-in-terms-of-cohomology}"'] \\
\cofib\left[R\Gamma(\Spec\OI \times X,\cM)\xrightarrow{\iota-\tau_{\cM}} R\Gamma(\Spec\OI \times X,\cN)\right]\otimes_{\OI}\KI \arrow[d,"\wr~\text{derived~formal~functions}"] \\
\cofib\left[R\Gamma(\Spf \OI \hat{\times} X,\hat{\cM})\xrightarrow{\iota-\tau_{\cM}} R\Gamma(\Spf \OI \hat{\times} X,\hat{\cN})\right]\otimes_{\OI}\KI \arrow[d,"\wr~\eqref{eq:diagram-completed-cohomology}~\text{and~Lemma}~\ref{lem:id-tau-isomorphism-on-TO}"] \\
\left[\hat{\cM}_{\infty}\xrightarrow{\id-\tau_M} \hat{\cN}_{\infty}\right]\otimes_{\OI}\KI \quad (\text{in~degrees~}[-1,0])\arrow[d,"\wr~\text{Lemma}~\ref{lem:injectivity-of-V-to-W}"] \\
\displaystyle\frac{\hat{\cN}_{\infty}}{(\iota-\tau_{\cM})(\hat{\cM}_{\infty})}\otimes_{\OI}\KI \arrow[d,"\wr~\text{Theorem}~\ref{thm:shtuka-model-to-hodge-pink-additive}"] \\
\displaystyle\frac{(M+\tau_M(\tau^*M))\otimes_{A\otimes F}\FI[\![\fj]\!]}{M\otimes_{A\otimes F}\FI[\![\fj]\!]}
\end{tikzcd}\nonumber
\end{equation}
\end{Definition}

\begin{Example}
Let $\underline{M}=\underline{A}(n)$ be the Carlitz $n$th twist over $K=\bF(\theta)$ (Examples \ref{ex:shtuka-moduel-carlitz-twists} \& \ref{ex:shtuka-model-carlitz-resumed}). Regardless of what shtuka model is chosen, one computes that the mock regulator is given by mapping the class of $\xi\in \KI\langle t \rangle$ to that of $e(t)\in (t-\theta)^{-n}\mathbb{F}[t,\theta]$, obtained as 
\[
e(t):=g-\frac{\tau(g)}{(t-\theta)^n}
\]
where $g\in \bF[t,\theta]$ denotes the projection of $\xi$ under the decomposition $\KI\langle t \rangle=\bF[t,\theta]\oplus \frac{1}{\theta}\bF[\![\frac{1}{\theta}]\!]\langle t \rangle$.

In general, we conjecture that the mock regulator only depends on $\underline{M}$ and not on the choice of the $C\times X$-shtuka model. 
\end{Example}

\begin{proof}[Proof of Theorem \ref{thm:rank-dim}]
For weight reasons, $\Hom_{\AMot_F}(\mathbbm{1},\underline{M})=0$. In particular, Proposition \ref{prop:motivic-exact-sequence} implies that we have an exact sequence of $\KI$-vector spaces:
\[
(\Betti{\underline{M}})^+_{\KI} \hookrightarrow \displaystyle\frac{\left\{\xi\in M\otimes_{A\otimes F}\FI\langle A \rangle~|~\xi-\tau_M(\tau^*\xi)\in N_{\cO_F}\right\}}{M_{\cO_F}}\otimes_{A}\KI \twoheadrightarrow \Ext^{1,\text{reg}}_{\cO_F,\infty}(\mathbbm{1},\underline{M})\otimes_A\KI. 
\]
On the other hand, by Theorem \ref{thm:extensions-of-MHPS-arising-from-motives}, we have an exact sequence of $\KI$-vector spaces:
\[
 (\Betti{\underline{M}})^+_{\KI}\hookrightarrow \displaystyle\frac{(M+\tau_M(\tau^*M))\otimes_{A\otimes F}\FI[\![\fj]\!]}{M\otimes_{A\otimes F}\FI[\![\fj]\!]} \twoheadrightarrow \Ext^{1,\text{ha}}_{\infty}(\mathbbm{1}^+,\sH^+(\underline{M})). 
\]
The middle terms have the same dimensions as they are isomorphic via $\rho(\underline{\cM})$. Consequently, the right-hand terms have the same dimensions as well.
\end{proof}

\subsection{The linearized Hodge-Pink regulator is injective}
Let $\underline{M}=(M,\tau_M)$ be a rigid analytically trivial $A$-motive over $F$ whose weights are all negative. Assume as above that $F$ admits a unique place above $\infty$. In this final section, we prove that the $\KI$-linearized Hodge-Pink regulator
\[
\Reg(\underline{M})_{\KI} : \Ext^{1,\mathrm{reg}}_{\mathcal{O}_F,\infty}(\mathbbm{1},\underline{M})\otimes_A \KI \longrightarrow \Ext^{1,\mathrm{ha}}_{\infty}(\mathbbm{1}^+,\sH^+(\underline{M}))
\]
is an isomorphism. Theorem \ref{thm:rank-dim} shows that the source and target of $\Reg(\underline{M})_{\KI}$ have the same dimension as vector spaces over $\KI$; it therefore suffices to show that $\Reg(\underline{M})_{\KI}$ is injective. 

As before, we denote by $M_{\mathcal{O}_F}$ the maximal integral model of $\underline{M}$ over $\mathcal{O}_F\subset F$ and by $N_{\cO_F}$ the intersection
\[
N_{\mathcal{O}_F}=M_{\mathcal{O}_F}[\mathfrak{j}^{-1}]\cap (M+\tau_M(\tau^*M)). 
\]
Let $L$ be a complete $\FI$-subfield of $\CI$ such that Theorem \ref{thm:analytic-continuation} holds. We also recall that we denoted by $\mathfrak{U}$ the analytification of $(C\setminus \{\infty\})\times \Spec L$ and by $\mathcal{M}_J(\mathfrak{U})$ the ring of meromorphic functions on $\mathfrak{U}$ whose poles are supported at $J=\{\fj,\fj_1,\fj_2,\ldots\}$.

By the explicit description of the regulator given in Proposition \ref{prop:explicit-computation-regulator}, a class in its kernel admits a representative $\psi$ whose Laurent expansion becomes regular at $\mathfrak j$ after subtracting an element of $(\underline M_B^+)_{K_\infty}$. The following proposition then shows that this representative defines the zero extension class and thus implies Theorem \ref{thm:regulator-injective-intro}.
\begin{Proposition}\label{prop:regulator-injective}
Let $s\geq 1$ and choose $m_1,\ldots,m_s\in N_{\mathcal{O}_F}$. For each $i\in\{1,\ldots,s\}$, let $\xi_i\in M\otimes_{A\otimes F}\mathcal{M}_J(\mathfrak{U})$ be a solution of $\xi_i-\tau_M(\tau^*\xi_i)=m_i$. Choose $b_1,\ldots,b_s\in\KI$ and set 
\[
\psi:=b_1\xi_1+\ldots+b_s\xi_s \in (M\otimes_{A\otimes F}\mathcal{M}_J(\mathfrak{U}))\otimes_A \KI.
\]
Assume that $\psi$ is regular at $\mathfrak{j}$; \emph{i.e.} that its image under the Laurent expansion map belongs to $M\otimes_{A\otimes F}L[\![\fj]\!]$. Then 
\[
\psi\in M_{\mathcal{O}_F}\otimes_A \KI.
\]
\end{Proposition}

The proof of Proposition \ref{prop:regulator-injective} has three steps. We first use the functional equation to remove the possible poles at the Frobenius translates of $\mathfrak j$. We then use the negativity of the weights to control growth at infinity. Finally, we combine the resulting formal integrality with boundedness of the coefficients to show that $\psi$ belongs to $M_{\mathcal O_F}\otimes_A K_\infty$.

We first establish the properties of the auxiliary algebra \(\mathcal M_J(D_\infty^\times)\) needed in the proof.

\subsubsection*{Meromorphic functions on $D_{\infty}^{\times}$}
The first step in the proof of Proposition \ref{prop:regulator-injective} is to embed $\KI$ and $\mathcal{M}_J(\mathfrak{U})$ into a common algebra $\mathcal{M}_J(D_{\infty}^{\times})$ whose definition we now give. 

\begin{Construction}[Analytic and meromorphic functions on $D_{\infty}^{\times}$.]
Let $L$ be a complete non-archimedean valued field (for the sake of this section, the reader is free to take a complete $\FI$-subfield of $\CI$ as above) equipped with a continuous map $\iota:F\hookrightarrow L$ with respect to the topology induced by $v_{\infty}$. We denote by $\mathcal{O}_L$ its valuation ring and by $k_L$ its residue field (which then contains a copy of $\bF_{\infty}$). We define $C_L^{\mathrm{an}}$ as the rigid analytification of $C_L:=C\times \Spec L$. The points of $C_L^{\mathrm{an}}$ are in one-to-one correspondence with the closed points of $C_L$ and, since the latter is proper over $L$, there is a well-defined \emph{reduction map} of the underlying sets of closed points
\[
\mathrm{red}:C_L^{\mathrm{an}}=|C_L| \longrightarrow |C_{k_L}|.
\]
The subset $\{\infty\}\times k_L$ of $|C_{k_L}|$ corresponds to a finite set of places $\{\overline{\infty}_\sigma\}_{\sigma}$, indexed by the $\bF$-linear embeddings $\bF_{\infty}\to L$, and we define $D_{\infty} \subset C_L^{\mathrm{an}}$ as its inverse image under $\mathrm{red}$. Geometrically, $\{\infty\}\times L$ decomposes into the same number of points $\{\infty_{\sigma}\}_{\sigma}$ in $C_L^{\mathrm{an}}$, and $D_{\infty}$ into the disjoint union of the rigid open unit disks $D_{\sigma}=\mathrm{red}^{-1}(\{\overline{\infty}_{\sigma}\})$ around these points. Note that there is a distinguished disk $D_{\iota}$, corresponding to the infinite point associated with $\iota$.\\
A standard computation in rigid analytic geometry shows that $D_{\infty}$ is an admissible subdomain of $C_L^{\mathrm{an}}$ and identifies the analytic functions on $D_{\infty}$ as a subring of $(\mathbb{F}_{\infty}\otimes L)[\![\pi_{\infty}]\!]$, where $\pi_{\infty}$ is any uniformizer of $C$ at $\infty$:
\begin{equation}\label{eq:ODinf}
\mathcal{O}(D_{\infty}) =\left\{\sum_{n\geq 0}{c_n \pi_{\infty}^n}~\Bigg |~\text{for all }0<r<1:~|c_n|r^n\xrightarrow[n\to \infty]{} 0\right\}.
\end{equation}
The morphism $\id\times \Frob_L$ of $C_L$ analytifies to an endomorphism of $C_L^{\mathrm{an}}$ stabilizing $D_{\infty}$, and therefore inducing an endomorphism of $\mathcal{O}(D_{\infty})$ which we denote by $f\mapsto f^{(1)}$. On analytic functions, it acts by
\[
f=\sum_{n\geq 0}{c_n \pi_{\infty}^n} \longmapsto f^{(1)}=\sum_{n\geq 0}{(\id_{\bF_{\infty}}\otimes \Frob_q)(c_n) \pi_{\infty}^n}.
\]
The algebra $\mathbb{F}_{\infty}\otimes L$ decomposes as a finite product of copies of $L$, indexed over the embeddings $\sigma:\bF_{\infty}\to L$, resulting in the decomposition 
\[
\mathcal{O}(D_{\infty})=\prod_{\sigma}\mathcal{O}(D_{\sigma}), \quad \mathcal{O}(D_{\sigma})=\left\{\sum_{n\geq 0}{c_{n,\sigma} \pi_{\infty}^n}\in L[\![\pi_{\infty}]\!]~\Bigg|~\text{for all }0<r<1:~|c_{n,\sigma}|r^n\xrightarrow[n\to \infty]{} 0\right\}.
\]
We define $D_{\infty}^{\times}$ as the subdomain of $D_{\infty}$ obtained by removing the centers $\{\infty_{\sigma}\}_{\sigma}$. The computation of the analytic functions on $D_{\infty}^{\times}$ is standard as well, and gives
\begin{equation}\label{eq:ODinfx}
\mathcal{O}(D_{\infty}^{\times}) =\left\{\sum_{n\in \mathbb{Z}}{c_n \pi_{\infty}^n}~\Bigg|~\text{for all }0<r<1:~|c_n|r^n\xrightarrow[n\to \pm \infty]{} 0\right\}.
\end{equation}
Note that $f\mapsto f^{(1)}$ extends to $\mathcal{O}(D_{\infty}^{\times})$ and that we have a similar decomposition $\mathcal{O}(D_{\infty}^{\times})=\prod_{\sigma}{\mathcal{O}(D_{\sigma}^{\times})}$. \\
The algebra $\mathcal{M}(D_{\infty}^{\times})$ of meromorphic functions on $D_{\infty}^{\times}$ is then defined as the localization of $\mathcal{O}(D_{\infty}^{\times})$ at regular elements. Algebraically, it corresponds to 
\begin{equation}\label{eq:MDinfx}
\mathcal{M}(D_{\infty}^{\times})=\prod_{\sigma}{\mathcal{M}(D_{\sigma}^{\times})}, \quad \mathcal{M}(D_{\sigma}^{\times})=Q(\mathcal{O}(D_{\sigma}^{\times})).
\end{equation}
The reader not familiar with rigid analytic geometry can very well take \eqref{eq:ODinf}, \eqref{eq:ODinfx} and \eqref{eq:MDinfx} as definitions for the algebras $\mathcal{O}(D_{\infty})$, $\mathcal{O}(D_{\infty}^{\times})$ and $\mathcal{M}(D_{\infty}^{\times})$.
\end{Construction}

As in Subsection \ref{subsection:analytic continuation}, we let $\fj_n\subset A\otimes L$ denote the maximal ideal given by the kernel of the map induced from $A\times L\to L$, $(a,l)\mapsto \iota(a)^{q^n}l$, and write $J$ for the infinite set $\{\fj_n~|~n\geq 0\}$. Given a subset $I$ of $J$, we may consider the ring $\mathcal{M}_I(D_{\infty}^{\times})$ of meromorphic functions on $D_{\infty}^{\times}$ whose poles are supported at $I$. Explicitly, it is given by the inverse limit
\begin{equation}\label{eq:MJDinfx}
\mathcal{M}_I(D_{\infty}^{\times})=\lim_{0<r\leq R<1} \left\{\sum_{n\in \mathbb{Z}}{c_n \pi_{\infty}^n}~\Bigg|~|c_n|R^n\xrightarrow[n\to \infty]{} 0,~|c_n|r^n\xrightarrow[n\to -\infty]{} 0, \right\}\left[\mathfrak{i}^{-1}~|~\mathfrak{i}\in I\right]
\end{equation}
indexed over pairs of real numbers $(r,R)$ satisfying $0<r\leq R<1$. For each pair $(r,R)$, the corresponding algebra consists of meromorphic functions on the closed annuli of inner radius $r$ and outer radius $R$, with poles of finite order supported at the points of $I$. Only finitely many ideals $\mathfrak{j}_n$ have to be inverted for each fixed pair.

\begin{Construction}[Mapping $\mathcal{M}_I(\mathfrak{U})$ to $\mathcal{M}_I(D_{\infty}^{\times})$]
Elements in $\mathcal{O}(\mathfrak{U})$ embed in $\mathcal{O}(D_{\infty}^{\times})$. To see this algebraically, recall that given a basis $(a_i)_{i\geq 0}$ of $A$ over $\bF$ ordered by increasing degree, any element $f\in\mathcal{O}(\mathfrak{U})$ can be written uniquely as a series
\[
f=\sum_{i\geq 0}{a_i\otimes \ell_i} \quad \text{with}~\ell_i\in L
\]
where, for all $\rho>1$, the sequence $(|\ell_i|\rho^{\deg(a_i)})_i$ is a null-sequence. We may therefore express $f$ as 
\begin{equation}\label{eq:well-defined-element-in-O(Dinfx)}
\sum_{n\in \mathbb{Z}}{\left(\sum_{\substack{i\geq 0 \\ \deg(a_i)\geq - d_{\infty} n}}{\lambda_{in}\otimes \ell_i}\right)\pi_{\infty}^n}
\end{equation}
where $\lambda_{in}\in \bF_{\infty}$ denotes the $\pi_{\infty}^n$-coefficient of $a_i$ seen in $\KI=\bF_{\infty}(\!(\pi_{\infty})\!)$ and where $d_{\infty}$ is the degree of the place $\infty$. Denote by $c_n$ the inner sum in parentheses. Let $0<r<1$. Set $\rho_0:=r^{-1/d_{\infty}}>1$, so that
\[
|c_n|r^n \leq \sum_{\substack{i\geq 0 \\ \deg(a_i)\geq -d_{\infty}n}}{|\ell_i|\rho_0^{-d_{\infty}n}}\leq \sum_{\substack{i\geq 0 \\ \deg(a_i)\geq -d_{\infty}n}}{|\ell_i|\rho_0^{\deg(a_i)}}.
\]
Since the summand of the right-hand series tends to zero as $i\to \infty$, as $n\to -\infty$ we obtain $|c_n|r^n\to 0$. As $n\to \infty$, we also have
\[
|c_n|r^n\leq r^{n/2}\left(\sum_{i=0}^{\infty}{|\ell_i| \rho_1^{\deg(a_i)}}\right)
\] 
where $\rho_1:=r^{-1/2d_{\infty}}>1$. The series in parentheses converges by assumption, and we also get $|c_n|r^n\to 0$. \\
This shows that the assignment $f\mapsto \eqref{eq:well-defined-element-in-O(Dinfx)}$ defines a map $\mathcal{O}(\mathfrak{U})\to \mathcal{O}(D_{\infty}^{\times})$. A simple verification shows that \eqref{eq:well-defined-element-in-O(Dinfx)} determines $f$, and so the map is injective. Inverting the ideals $\mathfrak{i}$ of $A\otimes L$ for $\mathfrak{i}\in I$ gives the injective map  $\mathcal{M}_I(\mathfrak{U})\to \mathcal{M}_I(D_{\infty}^{\times})$.
\end{Construction}

The two ring homomorphisms $\KI\to \mathcal{M}_J(D_{\infty}^{\times})$ and $\mathcal{M}_J(\mathfrak{U}) \to \mathcal{M}_J(D_{\infty}^{\times})$ assemble into a canonical map
\begin{equation}\label{eq:tensor-product-injects}
\mathcal{M}_J(\mathfrak{U})\otimes_A \KI\longrightarrow \mathcal{M}_J(D_{\infty}^{\times}).
\end{equation}
This common ambient algebra allows us to combine meromorphic functions on $\mathfrak U$ with scalars in $K_\infty$. The following lemma ensures that no additional relations are introduced.
\begin{Lemma}\label{lem:tensor-injects}
The map \eqref{eq:tensor-product-injects} is injective.
\end{Lemma}
\begin{proof}
Suppose that a nonzero element $\sum_i f_i\otimes h_i$ maps to zero in $\mathcal{M}_J(D_{\infty}^{\times})$. Choose such a relation with the number of terms minimal. Then the elements $(h_i)_i$ of $\KI$ are linearly independent over $K$. Since the possible poles lie outside the closed unit affinoid, the finitely many $f_i$ are analytic on a common slightly larger affinoid. Fix an $\bF$-basis $(a_j)_{j\geq 0}$ of $A$ ordered by increasing degree. Each $f_i$ has a unique expansion
\[
f_i=\sum_{j}{a_j\otimes \ell_{ij}}, \quad |\ell_{ij}|s^j\xrightarrow[j\to \infty]{} 0
\]
for some real number $s>1$. Up to multiplying by some well-chosen constant in $L^{\times}$, we may therefore assume that $\ell_{ij}\in \mathcal{O}_L$ for all $i$ and $j$ and that at least one of them is a unit. We may then reduce $f_i$ to $\bar{f}_i\in A\otimes k_L$.

Choose $r$ and $R$ such that $|\iota(\pi_{\infty})|<r\leq R<1$ and consider the map
\begin{equation}\label{eq:restriction-to-annulus}
\mathcal{M}_J(D_{\infty}^{\times})\longrightarrow \mathcal{O}(D_{r,R}).
\end{equation}
Through this map, any element of $\mathcal{M}_J(D_{\infty}^{\times})$ admits an expansion
\[
\sum_{n\in \mathbb{Z}}{c_n \pi_{\infty}^n} \in (\mathbb{F}_{\infty}\otimes L)[\![\pi_{\infty}^{\pm 1}]\!], \quad |c_n|R^n\xrightarrow[n\to \infty]{} 0,~|c_n|r^n\xrightarrow[n\to -\infty]{} 0.
\]
One verifies easily that the expansion of each $f_i$ belongs to $(\mathbb{F}_{\infty}\otimes \mathcal{O}_L)[\![\pi_{\infty}^{\pm 1}]\!]$ and that its reduction modulo $\mathfrak{m}_L$ is the expansion of $\bar{f}_i$ in $(\mathbb{F}_{\infty}\otimes k_L)[\![\pi_{\infty}^{\pm 1}]\!]$. It follows from this observation that the reduction modulo $\mathfrak{m}_L$ of the expansion of $\sum_i{f_ih_i}$ is
\[
\sum_i{\bar{f}_ih_i} \in \KI\otimes k_L.
\]
Choose an $\mathbb{F}$-linear functional $k_L\to\mathbb{F}$ that is nonzero on at least one coefficient of the family $(\bar f_i)_i$. Applying it to the relation $\sum_i\bar f_i h_i=0$ gives a nonzero $A$-linear relation among the $h_i$, contradicting their $K$-linear independence.
\end{proof}

Now that we have proven that both rings embed compatibly in $\mathcal{M}_J(D_{\infty}^{\times})$, we establish two further properties of this ring. First, given $n\geq 0$, there is a Laurent expansion map $\mathcal{M}_J(D_{\infty}^{\times})\to L(\!(\fj_n)\!)$ which extends Construction \ref{cons:Laurent-series-expansion-map}. We have the classical fact:
\begin{Lemma}[Removable singularities]\label{lem:regular-at-J}
Let $I'\subset I\subset J$ be subsets. Let $f\in \mathcal{M}_{I}(D_{\infty}^{\times})$ be such that its image via $\mathcal{M}_J(D_{\infty}^{\times})\to L(\!(\mathfrak{i})\!)$ lands in $L[\![\mathfrak{i}]\!]$ for all $\mathfrak{i}\in I'$. Then $f\in \mathcal{M}_{I\setminus I'}(D_{\infty}^{\times})$.
\end{Lemma}

Frobenius twists shift the allowed poles as follows. Given $I\subset J$, we define $I^{(1)}$ or more generally $I^{(h)}$ as $\{\mathfrak{j}_{n+h}~|~\mathfrak{j}_n\in I\}$. The next lemma is immediate.
\begin{Lemma}\label{lem:twist-polar-part}
Let $f\in \mathcal{M}_{I}(D_{\infty}^{\times})$. Then $f^{(1)}\in \mathcal{M}_{I^{(1)}}(D_{\infty}^{\times})$.  
\end{Lemma}

The second important feature of $\mathcal{M}_J(D_{\infty}^{\times})$ is that of \emph{growth at $\infty$}. This is measured via families of Gauss norms. 

\begin{Construction}[Norms on meromorphic functions]
Consider a real number $0<\rho<1$. Each algebra $\mathcal{O}(D_{\sigma}^{\times})$ carries a multiplicative Gauss norm $\|\cdot\|_{\rho,\sigma}$, given by the formula
\[
\|f\|_{\rho,\sigma}:=\sup_{n\in \bZ} |c_{n,\sigma}|\rho^n, \quad \text{with} \quad f=\sum_{n\in \mathbb{Z}}{c_{n,\sigma}\pi_{\infty}^n}\in \mathcal{O}(D_{\sigma}^{\times}).
\]
%Note that $\mathcal{O}(D_{\sigma}^{\times})$ is a Fréchet space, with family of norms $(\|\cdot\|_{\rho})_{0<\rho<1}$.
Because each $\|\cdot \|_{\rho,\sigma}$ is multiplicative, it extends to the field of fractions $\mathcal{M}(D_{\sigma}^{\times})$ of $\mathcal{O}(D_{\sigma}^{\times})$. By taking the maximum of the $\|\cdot \|_{\rho,\sigma}$ over $\sigma$, we obtain a norm $\|\cdot \|_{\rho}$ on $\mathcal{M}(D_{\infty}^{\times})$; be aware that this procedure cuts out the multiplicativity property. \\
In addition, one verifies easily that 
\begin{equation}\label{eq:behaviour-of-norms-under-twist}
\|f^{(1)}\|_{\rho}=\|f\|_{\rho^{1/q}}^q.
\end{equation}
\end{Construction}

The next lemma controls the growth of rational functions near infinity.
\begin{Lemma}\label{lem:rational-function-growth}
Let $f\in Q(A\otimes L)$ and let $s=\min_{\sigma}\{\mathrm{ord}_{\infty_{\sigma}}(f)\}$. As $\rho$ tends to zero,
\[
\|f\|_{\rho} =O(\rho^{s}).
\]
\end{Lemma}
\begin{proof}
Write $f=\frac{g}{h}$ with $g,h\in A\otimes L$ and $h\neq 0$. We may assume $g\neq 0$ since the assertion is immediate for $f=0$. Write the expansions of $g$ and $h$ at $\infty_{\sigma}$:
\[
g_{\sigma}=\sum_{n\geq u_{\sigma}}{g_{n,\sigma}\pi_{\infty}^n}, \quad h_{\sigma}=\sum_{n\geq v_{\sigma}}{h_{n,\sigma}\pi_{\infty}^n} \quad \text{in }L(\!(\pi_{\infty})\!)
\]
where $u_{\sigma}$ and $v_{\sigma}$ are chosen so that $g_{u_{\sigma},\sigma},h_{v_{\sigma},\sigma}\neq 0$. Since $g\in A\otimes L$, there exists a real constant $M_{\sigma}> 0$ such that $|g_{n,\sigma}|\leq M_{\sigma}$ and, in particular,
\[
\text{for }n>u_{\sigma}:\quad |g_{n,\sigma}|\rho^n \leq M_{\sigma}\rho^{u_{\sigma}+1}
\]
Choosing $\rho$ so that $M_{\sigma}\rho <|g_{u_{\sigma},\sigma}|$, we obtain 
\[
\|g_{\sigma}\|_{\rho,\sigma}=|g_{u_{\sigma},\sigma}|\rho^{u_{\sigma}}.
\]
Applying the same argument to $h_\sigma$, for each $\sigma$, gives the following equality for all sufficiently small $\rho$:
\[
\|f\|_{\rho}=\max_{\sigma}\left\{\frac{|g_{u_{\sigma},\sigma}|}{|h_{v_{\sigma},\sigma}|}\rho^{u_{\sigma}-v_{\sigma}}\right\}.
\]
This establishes the assertion with $s=\min_{\sigma}(u_{\sigma}-v_{\sigma})=\min_{\sigma}\{\mathrm{ord}_{\infty_{\sigma}}(f)\}$. 
\end{proof}

We conclude these preliminaries with a boundedness property. Any element $f$ of $\mathcal{M}_J(D_{\infty}^{\times})$, by restriction to the annulus $D_{r,R}$ with $|\iota(\pi_{\infty})|<r\leq R<1$, admits a unique expansion
\[
f=\sum_{n\in \mathbb{Z}}{c_n\pi_{\infty}^n} \quad c_n\in \mathbb{F}_{\infty}\otimes L.
\]
We say that $f$ has \emph{bounded coefficients} if there exists a positive constant $C$ such that $|c_n|\leq C$ for all $n\in \mathbb{Z}$. 

\begin{Lemma}\label{lem:bounded-coefficients}
Any element in the image of \eqref{eq:tensor-product-injects} has bounded coefficients.
\end{Lemma}
\begin{proof}
The property of having bounded coefficients, inside the ring $\mathcal{M}_J(D_{\infty}^{\times})$, is preserved under sums and products. Since elements in $K_{\infty}$ have bounded coefficients, it suffices to show that the same holds for any element $g$ of $\mathcal{M}_J(\mathfrak{U})$. But this is clear as $g$ restricts to an element of $\mathcal{O}(D_{r,1})$, hence its coefficients are bounded by $\|g\|_1$. 
\end{proof}

\subsubsection*{Proof of Proposition \ref{prop:regulator-injective}}
We keep the field $L$ fixed above and use the notation of Proposition \ref{prop:regulator-injective}. In particular, set 
\[
\psi:=b_1\xi_1+\ldots+b_s\xi_s \in (M\otimes_{A\otimes F}\mathcal{M}_J(\mathfrak{U}))\otimes_A \KI.
\]
By Lemma \ref{lem:tensor-injects} together with the fact that $M$ is a projective $A\otimes F$-module, we may consider $\psi$ as an element in $M\otimes_{A\otimes F} \mathcal{M}_J(D_{\infty}^{\times})$. We assume that, after applying the Laurent expansion map, we have 
\[
\psi \in M\otimes_{A\otimes F}L[\![\fj]\!].
\]
Our aim is to prove that $\psi\in M_{\mathcal{O}_F}\otimes_A \KI$.

By Lemma \ref{lem:stable-lattice-if-negative-weights}, choose a $\tau_M$-stable $\mathcal{A}_{\infty}(F)$-lattice $\Lambda$ in $\mathcal{B}_{\infty}(M)$ and positive integers $h,d$ such that 
\begin{equation}\label{eq:contractible-lattice}
\tau_M^h(\tau^{h*}\Lambda) \subset \mathfrak{m}_{\infty}^d \Lambda.
\end{equation}
Using the Beauville-Laszlo lemma, the gluing of $M$ and $\Lambda$ produces a locally free sheaf $\mathcal{M}$ on $C_F:=C\times \Spec F$ which restricts to $M$ on the affine subscheme $\Spec A\otimes F$, together with a morphism $\tau_{\mathcal{M}}:\tau^* \mathcal{M}|_{C_F\setminus \mathfrak{j}}\longrightarrow \mathcal{M}|_{C_F\setminus \mathfrak{j}}$ extending $\tau_M$. Let $V=\Spec B_F$ be an affine open neighbourhood of $\infty_F$ in $C_F$ on which $\mathcal{M}$ is free. Let $(e_1,\ldots,e_r)$ be a basis of $\mathcal{M}(V)$, so that
\[
\Lambda:=\bigoplus_{i=1}^r \mathcal{A}_{\infty}(F)e_i.
\]
Using this basis, we define a norm on $M\otimes_{A\otimes F}\mathcal{M}(D_{\infty}^{\times})$ by
\[
\|f\|_{\Lambda,\rho}:=\max_i\{\|f_i\|_{\rho}\} \quad \text{where} \quad f=\sum_i{f_i e_i}.
\]
A different choice of basis will produce a different, but equivalent, norm. 
\begin{Lemma}\label{lem:tau-norm-estimates}
Let $f\in M\otimes_{A\otimes F}\mathcal{M}(D_{\infty}^{\times})$. There exists a constant $C>0$ such that, as $\rho$ tends to zero, we have
\[
\|\tau_M^h(\tau^{h*}f)\|_{\Lambda,\rho}\leq C\rho^d \|f\|^{q^h}_{\Lambda,\rho^{q^{-h}}}.
\]
In addition, suppose that $f\in(M\otimes_{A\otimes F}\mathcal{M}_J(\mathfrak{U}))\otimes_A\KI$ and that its image in $M\otimes_{A\otimes F}\mathcal{M}(D_{\infty}^{\times})$ belongs to $M\otimes_{A\otimes F}\mathcal{O}(D_{\infty}^{\times})$. Then there exists a constant $C>0$ such that, as $\rho$ tends to zero,
\[
\|\tau_M^h(\tau^{h*}f)\|_{\Lambda,\rho}
\leq C\rho^d\|f\|_{\Lambda,\rho}.
\]
\end{Lemma}
\begin{proof}
Condition \eqref{eq:contractible-lattice} amounts to saying that the matrix representing the family $(\tau_M^h(\tau^{h*}e_i))_{i}$ written in the basis $(e_i)_i$ is of the form $\pi_{\infty}^d H_h$ where $H_h$ is an $r\times r$-matrix whose coefficients belong both to $Q(A\otimes F)$ and to $\mathcal{A}_{\infty}(F)$. In particular, its coefficients are regular at $\infty_F$, and by Lemma \ref{lem:rational-function-growth} we deduce that $\|H_h\|_{\rho}$ is bounded by some positive constant $C$ as $\rho$ goes to zero. It follows that 
\[
\text{for all }i\in \{1,\ldots,r\}:\quad \|\tau_M^h(\tau^{h*} e_i)\|_{\Lambda,\rho}\leq C\rho^d.
\]
Combining \eqref{eq:behaviour-of-norms-under-twist} with the preceding bound gives the first assertion.

The second assertion follows from Lemma \ref{lem:bounded-coefficients} which implies that $\|f\|^{q^h}_{\Lambda,\rho^{q^{-h}}}\leq \|f\|_{\Lambda,1}^{q^h-1}\|f\|_{\Lambda,\rho}$.
\end{proof}

\begin{proof}[Proof of Proposition \ref{prop:regulator-injective}]
Set $\mu:=b_1m_1+\ldots+b_sm_s$, viewed as an element of $M\otimes_{A\otimes F}\mathcal{M}_{\{\mathfrak{j}\}}(D_{\infty}^{\times})$. Then $\psi$ satisfies the functional equation:
\begin{equation}\label{eq:functional-relation}
\psi-\tau_M(\tau^* \psi)=\mu.
\end{equation}
Let us first prove from \eqref{eq:functional-relation} that $\psi$ is analytic on $D_{\infty}^{\times}$; \emph{i.e.} that $\psi$ belongs to $M\otimes_{A\otimes F}\mathcal{O}(D_{\infty}^{\times})$. We shall proceed by induction on $i\geq 0$ to prove that $\psi$ belongs to $M\otimes_{A\otimes F}\mathcal{M}_{J^{(i)}}(D_{\infty}^{\times})$; this suffices because $M$ is a finitely generated projective $A\otimes F$-module, so tensoring with $M$ commutes with arbitrary intersections:
\[
\bigcap_{i\geq 0}\left(M\otimes_{A\otimes F}\mathcal{M}_{J^{(i)}}(D_{\infty}^{\times})\right) =M\otimes_{A\otimes F}\left(\bigcap_{i\geq 0}\mathcal{M}_{J^{(i)}}(D_{\infty}^{\times})\right),
\] 
and then we may use $\bigcap_{i\geq 0}\mathcal{M}_{J^{(i)}}(D_{\infty}^{\times})=\mathcal{O}(D_{\infty}^{\times})$ which follows from Lemma \ref{lem:regular-at-J}. 

This holds for $i=1$. If $\psi\in M\otimes_{A\otimes F}\mathcal{M}_{J^{(i)}}(D_{\infty}^{\times})$---which holds for $i=1$---then we have
\[
\tau_M(\tau^*\psi)\in M\otimes_{A\otimes F}\mathcal{M}_{J^{(i+1)}\cup \{\mathfrak{j}\}}(D_{\infty}^{\times})
\]
as follows from Lemma \ref{lem:twist-polar-part}. Equation \eqref{eq:functional-relation} then implies that 
\[
\psi=\mu+\tau_M(\tau^*\psi)\in M\otimes_{A\otimes F}\left(\mathcal{M}_{J^{(i)}}(D_{\infty}^{\times})\cap \mathcal{M}_{J^{(i+1)}\cup \{\mathfrak{j}\}}(D_{\infty}^{\times})\right) = M\otimes_{A\otimes F}\mathcal{M}_{J^{(i+1)}}(D_{\infty}^{\times}),
\]
completing the induction step.

Next, we show that $\psi$ belongs to $\mathcal{B}_{\infty}(M_L):=M\otimes_{A\otimes F} \mathcal{B}_{\infty}(L)$. Let $h\geq 1$ be a positive integer and let $\mu_h:=\sum_{h>k\geq 0}{\tau_M^k(\tau^{k*}\mu)}$ so that 
\begin{equation}\label{eq:functional-relation-h}
\psi-\tau_M^h(\tau^{h*}\psi)=\mu_h. 
\end{equation}
If $h$ is chosen as above, then Lemma \ref{lem:tau-norm-estimates} implies that choosing $\rho>0$ small enough gives $\|\tau_M^h(\tau^{h*}\psi)\|_{\Lambda,\rho}<\|\psi\|_{\Lambda,\rho}$. We obtain
\begin{equation}\label{eq:estimate}
\|\psi\|_{\Lambda,\rho}=\|\psi-\tau_M^h(\tau^{h*}\psi)\|_{\Lambda,\rho}=\|\mu_h\|_{\Lambda,\rho}=O(\rho^s)
\end{equation}
for some integer $s$. The coordinates $(\psi_i)_i$ of $\psi$ in the basis $(e_i)_{i}$ are elements of the subring generated by $Q(A\otimes F)$ and $\mathcal{O}(D_{\infty}^{\times})$. On sufficiently small punctured disks, these coordinates are analytic Laurent series. The bound \eqref{eq:estimate} forces their coefficients in degrees below $s$ to vanish. Thus each coordinate has only finitely many negative Laurent coefficients, and its expansion belongs to $\mathcal{B}_{\infty}(L)$. We obtain 
\[
\psi\in \mathcal{B}_{\infty}(M_L).
\]

We may now consider the relation \eqref{eq:functional-relation-h} in $\mathcal{B}_{\infty}(M_L)$ which we equip with the $\mathfrak{m}_{\infty}$-adic topology induced by any lattice. Given an element $b$ therein, the sequence $(\tau_M^k(\tau^{k*}b))_{k\geq 0}$ tends to zero by \eqref{eq:contractible-lattice}. Iterating the functional equation and using this convergence, we obtain
\[
\psi=\sum_{k\geq 0}{\tau_M^k(\tau^{k*}\mu)} \in \mathcal{B}_{\infty}(M_{\mathcal{O}_F}).
\]
Because $M_{\mathcal{O}_F}$ is flat over $A\otimes \mathcal{O}_F$, this gives 
\[
\psi \in (M_{\mathcal{O}_F}\otimes_{A\otimes \mathcal{O}_F} \mathcal{B}_{\infty}(\mathcal{O}_F))\cap (M_{\mathcal{O}_F}\otimes_{A\otimes \mathcal{O}_F} \mathcal{O}(D_{\infty}^{\times}))=M_{\mathcal{O}_F}\otimes_{A\otimes \mathcal{O}_F} (\mathcal{B}_{\infty}(\mathcal{O}_F)\cap \mathcal{O}(D_{\infty}^{\times}))
\]
where the latter intersection is understood as an intersection of submodules of $(\mathbb{F}_{\infty}\otimes L)[\![\pi_{\infty}^{\pm 1}]\!]$ (be aware that this is not a ring). We interpret $M_{\mathcal{O}_F}$ as a direct summand in a finite free $A\otimes \mathcal{O}_F$\nobreakdash-module $(A\otimes \mathcal{O}_F)^g$. Because $\psi$ comes from an element in $M\otimes_{A\otimes F}(\mathcal{M}_J(\mathfrak{U})\otimes_A\KI)$, it has bounded coefficients by Lemma \ref{lem:bounded-coefficients}, and thus belongs to $(\mathbb{F}_{\infty}\otimes V)(\!(\pi_{\infty})\!)^g$ for some finite dimensional $\bF$-vector space $V\subset \mathcal{O}_F$ and, in particular, belongs to $(\KI\otimes \mathcal{O}_F)^g$. Projecting back onto the $M_{\mathcal{O}_F}$-component, we obtain
\[
\psi \in M_{\mathcal{O}_F}\otimes_{A\otimes \mathcal{O}_F} (\KI\otimes \mathcal{O}_F)= M_{\mathcal{O}_F}\otimes_A \KI
\]
as desired.
\end{proof}

This proves Proposition \ref{prop:regulator-injective} and hence Theorem \ref{thm:regulator-injective-intro}. Together with Theorem \ref{thm:rank-dim}, this establishes Theorem B.

\end{document}